\documentclass{amsart}
\usepackage{amscd,amssymb,amsxtra}
\usepackage[mathscr]{eucal}
\usepackage[all]{xy}

\title[Quadratic functions]
{Quadratic functions in geometry, topology,
and $M$-theory}

\author{M.~J.~Hopkins}
\address{Department of Mathematics \\ Massachusetts Institute of
Technology\\Cambridge, MA 02139-4307}
\email{mjh@math.mit.edu}
\thanks{The first author would like to acknowledge support from NSF
grant \#DMS-9803428}
\author{I.~M.~Singer}
\address{Department of Mathematics \\ Massachusetts Institute of
Technology\\Cambridge, MA 02139-4307}
\email{ims@math.mit.edu}
\thanks{The second author would like to acknowledge support from DOE
grant \#DE-FG02-ER25066}
     
\newtheorem{thm}[equation]{Theorem}
\newtheorem{cor}[equation]{Corollary}
\newtheorem{lem}[equation]{Lemma}
\newtheorem{prop}[equation]{Proposition}

\ifx\undefined\theoremstyle
	\newtheorem{definition}[equation]{Definition}
	\newenvironment{defin}{\begin{definition}\rm}{\end{definition}}
	\newtheorem{conj}[equation]{Conjecture}
	\newtheorem{rem}[equation]{Remark}
	\newtheorem{rems}[equation]{Remarks}

\else
	\theoremstyle{definition}

	\newtheorem{defin}[equation]{Definition}

	\theoremstyle{remark}

	\newtheorem{rem}[equation]{Remark}
	
\ifx\undefined\prem	\newtheorem{prem}{\bf Working Remark} \fi

\fi

\newtheorem{eg}[equation]{Example}

\ifx\undefined\pf
\newenvironment{pf}{\bigskip{\em Proof:\/}}{\qed\medskip}
\newenvironment{pf*}[1]{\bigskip{\em #1:\/}}{\qed\medskip}
\fi

\ifx\undefined\numberwithin
\def\numberwithin#1#2{\makeatletter\@ifundefined{c@#1}{\@nocnterrr}{%
  \@ifundefined{c@#2}{\@nocnterr}{%
  \@addtoreset{#1}{#2}%
  \toks@\expandafter\expandafter\expandafter{\csname the#1\endcsname}%
  \expandafter\xdef\csname the#1\endcsname
    {\expandafter\noexpand\csname the#2\endcsname
     .\the\toks@}}}\makeatother}\fi

\ifx\undefined\qed\newcommand{\qed}{\hfil\rule{4pt}{6pt}\bigskip}\fi
\ifx\undefined\\operatorname\newcommand{\operatorname}[1]{\mathop{\mbox{\rm #1}}}\fi

\newcommand{\Aut}{\operatorname{Aut}}
\newcommand{\ext}{\operatorname{Ext}}

\DeclareMathOperator{\coker}{coker}

\newcommand{\Q}{{\mathbb Q}}
\newcommand{\Z}{{\mathbb Z}}
\newcommand{\R}{{\mathbb R}}

\newcommand{\zerowidth}[1]{\hbox to 0pt{\hss$\displaystyle #1$\hss}}
\newcommand{\LL}{[\mkern-2mu[}
\newcommand{\RR}{]\mkern-2mu]}

\newcommand{\C}{\mathbb C}

\ifx\undefined\eqref\newcommand{\eqref}[1]{\rm (\ref{#1})}\fi

\newcounter{thmItem}

\newcommand{\thmItemref}[1]
 	 {{\rm \ref{#1})}}

\newenvironment{thmList}{\begin{list}%
{\rm \roman{thmItem})}{\usecounter{thmItem}
\setlength{\labelwidth}{2em}
\setlength{\itemindent}{2em}
\setlength{\leftmargin}{0pt}
\setlength{\listparindent}{0pt}
\setlength{\parsep}{0pt}
\setlength{\partopsep}{0pt}
\setlength{\itemsep}{\medskipamount}
\setlength{\topsep}{\medskipamount}
}}{\end{list}}

\newcounter{textItem}

\newcommand{\textItemref}[1]
 	 {{\rm (\ref{#1})}}

\newenvironment{textList}{\begin{list}%
{\rm (\arabic{textItem})}{\usecounter{textItem}
\setlength{\labelwidth}{2em}
\setlength{\itemindent}{2em}
\setlength{\leftmargin}{0pt}
\setlength{\listparindent}{0pt}
\setlength{\parsep}{0pt}
\setlength{\partopsep}{0pt}
\setlength{\itemsep}{\medskipamount}
\setlength{\topsep}{\medskipamount}
}}{\end{list}}

\newcounter{condItem}

\ifx\undefined\slot
\newcommand{\slot}{\,-\,}
\fi

\DeclareMathOperator{\su}{\mathit{SU}}
\DeclareMathOperator{\bsu}{\mathit{BSU}}
\DeclareMathOperator{\bu}{\mathit{BU}}

\DeclareMathOperator{\spin}{Spin}
\DeclareMathOperator{\spinc}{Spin^{c}}

\DeclareMathOperator{\harm}{harm}
\DeclareMathOperator{\HA}{\mathit{HA}}
\DeclareMathOperator{\HZ}{\mathit{H}\Z}
\DeclareMathOperator{\HQ}{\mathit{H}\Q}

\newcommand{\hz}{\HZ}
\newcommand{\hq}{\HQ}

\newcommand{\chs}{\check H}

\newcommand{\mok}[1]{\mo\left(#1\right)}
\DeclareMathOperator{\mo}{\mathit{MO}}
\DeclareMathOperator{\mso}{\mathit{MSO}}
\DeclareMathOperator{\bso}{\mathit{BSO}}
\DeclareMathOperator{\mspin}{MSpin}
\DeclareMathOperator{\bspin}{BSpin}

\newcommand{\rz}{\R/\Z}
\newcommand{\qz}{\Q/\Z}

\newcommand{\tpi}{2\pi i}

\DeclareMathOperator{\Sq}{Sq}

\newcommand{\cp}{{CP^\infty}}
\newcommand{\rp}{{RP^\infty}}

\newcommand{\dbar}{\bar \partial}

\newcommand{\hcat}{\mathcal{H}}

\newcommand{\chcat}{\check\hcat}

\newcommand{\chcochain}[1]{\Check C(#1)}

\newcommand{\chzero}{\Check C}
\newcommand{\chzerococycle}{\Check Z}
\newcommand{\chcocycle}[1]{\Check Z(#1)}

\newcommand{\chh}[1]{\Check H(#1)}
\newcommand{\chhs}[1]{\chh{#1}^{#1}}

\newcommand{\wustruct}{\lambda}
\newcommand{\lineb}{\mathcal L}

\newcommand{\simp}[1]{#1_{\bullet}}
\DeclareMathOperator{\sing}{sing}
\newcommand{\gsing}{\sing^{\text{geom}}}

\newcommand{\sets}{\underline{\text{Sets}}}

\newcommand{\spi}{\pi^{\text{simp}}}

\newcommand{\Omegasimp}{\Omega^{\text{simp}}}

\newcommand{\orientation}{$\chs$-orientation}
\newcommand{\oriented}[1]{$\chs$-oriented{#1}}

\newcommand{\rnbar}[1]{\bar{\R}^{#1}}

\newcommand{\slant}{/}

\DeclareMathOperator{\ch}{ch}

\makeatletter
\newcommand{\dirac}{\mathrel{\mathpalette\c@ncel\partial}}
\makeatother

\DeclareMathOperator{\determinant}{det}
\renewcommand{\det}{\determinant}

\newcommand{\bsowubar}{\bso\langle\bar\nu_{2k} \rangle}
\newcommand{\msowubar}{\mso\langle\bar\nu_{2k} \rangle}
\newcommand{\msowubarh}{\mso\langle\bar\nu_{2k} \rangle_{h\Z/2}}
\newcommand{\msobwu}{\mso\langle\beta\nu_{n+1} \rangle}
\newcommand{\msobwuh}{\mso\langle\beta\nu_{n+1} \rangle_{h\Z/2}}
\newcommand{\msobwuk}{\msobwu\wedge K(\Z,2k)}

\newcommand{\msobwukkh}{\msobwu\wedge K(\Z(1),2k)\wedge K(\Z(1),2k)_{h\Z/2}}

\newcommand{\msobwukh}{\msobwu\wedge K(\Z(1),2k)_{h\Z/2}}
\newcommand{\msobwukhp}{\left(\msobwu\wedge K(\Z(1),2k)_{+} \right)_{h\Z/2}}
\newcommand{\bsobwu}{\bso\langle\beta\nu_{n+1} \rangle}

\DeclareMathOperator{\fil}{filt}
\newcommand{\grass}[2]{{Gr}_{#1}\left(#2\right)}
\newcommand{\sgrass}[2]{\widetilde{Gr}_{#1}\left(#2\right)}
\newcommand{\spingrass}[2]{Gr^{\text{Spin}}_{#1}\left(#2 \right)}

\newcommand{\cat}[1]{\mathcal{#1}}

\DeclareMathOperator{\codim}{codim}

\DeclareMathOperator{\thom}{Thom}

\newcommand{\fred}{\mathcal F}

\newcommand{\grv}{\mathcal V}
\newcommand{\grw}{\mathcal W}

\newcommand{\nubar}{\bar{\nu}}
\newcommand{\zbar}{\bar{z}}
\newcommand{\nuhat}{\Check{\nu}}

\newcommand{\iotai}{\tilde{\iota}}

\numberwithin{equation}{section}

\begin{document}

\maketitle

\setcounter{tocdepth}{2}
\tableofcontents

\nocite{alvarez85:_cohom}
\nocite{brylinski93:_loop}
\nocite{deligne71:_theor_hodge}
\nocite{alvarez01:_beyon_ellip_genus}

\section{Introduction}
Quadratic refinements of the intersection pairing on a Riemann surface
appear to have two mathematical origins: one in complex
function theory dating back to Riemann in the 1870's, and one in
topology stemming from the work of Pontryagin in the 1930's.

Pontryagin's ideas were taken up and generalized by
Kervaire~\cite{kervaire60} in the late 1950's, who, among other things
used them to produce an example of a topological manifold of dimension
$10$ which does not admit a smooth structure.  Analogous invariants
for manifolds of other dimensions were investigated by many
topologists, most notably Kervaire, Browder~\cite{browder69:_kervair},
Brown-Peterson~\cite{brown66:_kervair,brown65:_kervair} and
Brown~\cite{brown72:kervair}, and play an important role in the
surgery classification of manifolds and in the homotopy groups of
spheres.

Riemann's quadratic function occurred in his theory of
$\vartheta$-func\-tions and while its topological aspects were
clarified in 1971 by Atiyah~\cite[Proposition~4.1]{atiyah71:_rieman}
and Mumford~\cite{mumford71:_theta}, until recently quadratic
functions in higher dimensions have remained in the province of
topology.  Our purpose in this paper is to bring to quadratic
functions in higher dimensions more of the geometry present in
Riemann's work.  There are two issues involved, one having to do with
constructing quadratic functions, and the other having to do with the
mathematical language in which to describe them.  As we explain below,
our motivation for doing this came from theoretical physics, and our
theory owes much to the papers~\cite{witten97:_five_m,witten97:_m} of
Witten.

In the case of Riemann surfaces, links between the
topological approach of Pontryagin and the analytic approach of
Riemann can be made using index theory.  In~\cite{atiyah71:_rieman}
Atiyah interprets Riemann's quadratic function in terms of the mod $2$
index of the Dirac operator.  It is also possible to deduce Riemann's
results from the theory of the determinant of the $\dbar$-operator.
Though this point of view seems relatively modern, it is arguably the
closest to Riemann's original analysis.  Riemann's quadratic function
occurred in the functional equation for his $\vartheta$-function.  The
$\vartheta$-function is (up to scale) the unique holomorphic section of
the determinant of the $\dbar$-operator, and it's functional equation
can studied from the symmetries of the determinant line.  In
\S\ref{subsec-determinants} we will give a proof of Riemann's results
along these lines.  Our proof also works in the algebraic setting.
While apparently new, it is related to that of Mumford, and gives
another approach to his results~\cite{mumford71:_theta}.

While it doesn't seem possible to construct quadratic functions in
higher dimensions using index theory alone, there is a lot to be
learned from the example of determinant line bundles on Riemann
surfaces.  Rather than trying to use the index of an operator, our
approach will be to generalize the index {\em formula}, i.e., the
topological index.  The index formula relates the determinant of the
$\dbar$ operator in dimension $2$ to the index of the $\spinc$ Dirac
operator in dimension $4$, and ultimately, quadratic functions in
dimension $2$ to the signature of $4$-manifolds.  Now on a
$4$-manifold $M$ the relation between $\spinc$-structures and
quadratic refinements of the intersection pairing has a simple
algebraic interpretation.  The first Chern class $\lambda$ of the
$\spinc$-structure is a characteristic element of the bilinear form on
$H^{2}(M)$:
\[
\int_{M}x\cup x \equiv \int_{M}x\cup \lambda \mod 2.
\]
The expression 
\begin{equation}\label{eq:71}
q(x)=\frac12 \int_{M}\left(x^{2}-x\lambda \right)
\end{equation}
is then a quadratic refinement of the intersection pairing.  It is
useful to compare this with the formula for the index of the
$\spinc$ Dirac operator:
\begin{equation}\label{eq:72}
\kappa(\lambda)=\frac18\int_{M}(\lambda^{2}-L(M)),
\end{equation}
where $L(M)$ is the characteristic class which gives the signature
when integrated over $M$.  Formula~\eqref{eq:71} gives the change of
$\kappa(\lambda)$ resulting from the change of $\spinc$-structure
$\lambda\mapsto \lambda-2x$.

The fact that~\eqref{eq:72} is an integer also has an algebraic
explanation: the square of the norm of a characteristic element of a
non-degen\-erate symmetric bilinear form over $\Z$ is always congruent
to the signature mod $8$.  This points the way to a generalization in
higher dimensions.  For manifolds of dimension $4k$, the
characteristic elements for the intersection pairing in the middle
dimension are the integer lifts $\lambda$ of the Wu-class $\nu_{2k}$.
The expression~\eqref{eq:72} is then an integer, and its variation
under to $\lambda\mapsto\lambda-2x$ gives a quadratic refinement of
the intersection pairing.  This can almost be described in terms of
index theory.  A $\spinc$-structure on a manifold of dimension $4k$
determines an integral Wu-structure, and the integer $\kappa(\lambda)$
turns out to be the index of an operator.  But we haven't found
good analytic way to understand the variation $\lambda\mapsto
\lambda-2x$.  In dimension $4$ this variation can implemented by
coupling the $\spinc$ Dirac operator to a $U(1)$-bundle with first
Chern class $-x$.  In higher dimensions one would need to couple the
operator to something manufactured out of a cohomology class of degree
$2k$.

In this paper we refine the expression~\eqref{eq:72} to a cobordism
invariant for families of manifolds.  The cobordism theory is the one
built from families $E/S$ of manifolds equipped with an integer
cocycle $\lambda\in Z^{2k}(E;\Z)$ whose mod $2$ reduction represents
the Wu class $\nu_{2k}$ of the relative normal bundle.  If $E/S$ has
relative dimension $(4k-i)$, then our topological interpretation
of~\eqref{eq:72} will produce an element 
\[
\kappa(\lambda)\in \tilde I^{i}(S)
\]
of a certain generalized cohomology group.  The cohomology theory
$\tilde I$ is a generalized cohomology theory known as the Anderson
dual of the sphere.  It is the dualizing object in the category of
cohomology theories (spectra).  When $i=2$, the group $\tilde
I^{2}(S)$ classifies graded line bundles.  By analogy with the case of
Riemann surfaces, we think of $\kappa(\lambda)$ the determinant line
of the Wu-structure $\lambda$ on $E/S$.  This ``generalized
determinant'' $\kappa(\lambda)$ can be coupled to cocycles of degree
$2k$, and can be used to construct quadratic functions.  

The relationship between Wu-structures, quadratic functions, and the
Kervaire invariant goes back to the early work on the Kervaire
invariant in~\cite{brown65:_kervair,brown65:_kervair,
brown66:_kervair,anderson66:_su_ko_kervair}, and most notably to the
paper of Browder~\cite{browder69:_kervair}.  It was further clarified
by Brown~\cite{brown72:kervair}.  The relationship between the
signature in dimension $4k$ and the Kervaire invariant in dimension
$(4k-2)$ was discovered by Milgram~\cite{milgram74:_surger} and
Morgan--Sullivan~\cite{morgan74}.  Our construction of
$\kappa(\lambda)$ is derived
from~\cite{brown72:kervair,milgram74:_surger,morgan74}, though our
situation is somewhat different.
In~\cite{brown72:kervair,milgram74:_surger,morgan74} the emphasis is
on surgery problems, and the class $\lambda$ is necessarily $0$.  In
this work it is the variation of $\lambda$ that is important.  Our
main technical innovation involves a systematic exploitation of
duality.

Even though $\kappa(\lambda)$ generalizes the determinant line, as
described so far our cobordism approach produces objects which are
essentially topological.  To enrich these objects with more geometric
content we introduce the language of {\em differential functions} and
{\em differential cohomology theories}.  Let $X$ be a topological
space, equipped with a real cocycle $\iota\in Z^{n}(X;\R)$, and $M$ a
smooth manifold.  A {\em differential function} from $M$ to
$(X,\iota)$ is a triple
\[
(c,h,\omega)
\]
with $c:M\to X$ a continuous function, $\omega\in \Omega^{n}(M)$ a
closed $n$-form, and $h\in C^{n-1}(M;\R)$ a cochain satisfying
\[
\delta h=\omega-c^{\ast}\iota.
\]
Using differential functions, we then revisit the basic constructions
of algebraic topology and introduce differential cobordism groups and
other differential cohomology theories.  It works out, for instance,
that the differential version of the group $\tilde I^{1}(S)$ is the
group of smooth maps from $S$ to $U(1)$.  Using the differential
version of $\tilde I^{2}(S)$ one can recover the category of
$U(1)$-bundles with connection.  In this way, by using {\em
differential} rather than continuous functions, our topological
construction refines to a something richer in its geometric aspects.

The differential version of $H^{k}(M;\Z)$ is the group of
Cheeger-Simons differential characters
$\chs^{k-1}(M)$~\cite{cheeger85:_differ}, and in some sense our theory
of differential functions is a non-linear generalization.  We began
this project intending to work entirely with differential characters.
But they turned out not robust enough for our purposes.  

The bulk of this paper is devoted to working out the theory of
differential function spaces.  To make them into {\em spaces} at all
we need to consider differential functions on the products $M\times
\Delta^{n}$ of $M$ with a varying simplex.  This forces us at the
outset to work with manifolds with corners (see
Appendix~\ref{sec:manif-with-corn}).  Throughout this paper, the term
{\em manifold} will mean {\em manifold with corners}.  The term {\em
manifold with boundary} will have its usual meaning, as will the term
{\em closed manifold}.

Using the language of differential function spaces and differential
cohomology theories the construction of quadratic functions in higher
dimension can be made to arise very much in the way they did for
Riemann in dimension $2$.  A differential integral Wu-structure is the
analogue of the canonical bundle $\lambda$, a choice of $\lambda/2$ is
the analogue of a $\vartheta$-characteristic (or $\spin$-structure),
and $\kappa(\lambda-2x)$ the analogue of the determinant of the
$\dbar$-operator coupled to a holomorphic line bundle.  For a more
detailed discussion, see \S\ref{sec-statement-of-results}.

Our interest in this project originated in a discussion with Witten.
It turns out that quadratic functions in dimension $6$ appear in the
ring of ``topological modular forms''~\cite{hopkins95:_topol_witten},
as a topologically defined mod $2$ invariant.  Modulo torsion, the
invariants coming from topological modular forms are accounted for by
index theory on loop space.  This suggested that it might be possible
to generalize Atiyah's interpretation of Riemann's quadratic function
to dimension $6$ using some kind of mod $2$ index on loop space.  We
asked Witten about this and he pointed out that he had used quadratic
refinements of the intersection pairing on certain $6$-manifolds in
describing the fivebrane partition function in
$M$-theory~\cite{witten97:_five_m,witten97:_m}.  The fivebrane
partition function is computed as the unique (up to scale) holomorphic
section of a certain holomorphic line bundle, and the quadratic
function is used to construct this line bundle.  We then realized that
an analogue of determinant lines could be used instead of mod $2$
indices to generalize Riemann's quadratic functions to higher
dimensions.

The organization of this paper is as follows.
Section~\ref{sec-statement-of-results} is devoted to background
material and the statement of our main result.  More specifically,
\S\ref{subsec-history} recalls the results of Riemann and Pontryagin.
In \S\ref{subsec-determinants} we give a proof of Riemann's results
using determinants.  Sections~\ref{subsec-differential-cocycles}
and~\ref{subsec-wustruct} introduce differential cocycles, which play
a role in higher dimensions analogous to the one played by line
bundles in dimension $2$.  We state our main result
(Theorem~\ref{thm-main}) in \S\ref{subsec-main-result}, and in
\S\ref{subsec-witten} relate it to Witten's construction.  In
\S\ref{sec-db-cs} we review Cheeger-Simons cohomology.  In
\S\ref{sec:gener-diff-cohom} we lay out the foundations of
differential function complexes and differential cohomology theories.
Section~\ref{sec-topology} contains the proof of
Theorem~\ref{thm-main}.

We had originally included in this paper an expository discussion,
primarily for physicists, describing the role of quadratic and
differential functions in the construction of certain partition
functions.  In the end we felt that the subject matter deserved a
separate treatment, which we hope to complete soon.

Our theory of differential function spaces provides a variation of
algebraic topology more suited to the needs of mathematical physics.
It has already proved useful in anomaly
cancellation~\cite{freed:_ramon_ramon_k}, and it appears to be a
natural language for describing fields and their action functionals.
We have many examples in mind, and hope to develop this point of view
in a later paper.

This work has benefited from discussions with many people.  We wish
to thank in particular Orlando Alvarez, Dan Freed, Jack Morava, Tom
Mrowka, C.~R.~Nappi, T.~Ramadas, and Edward Witten.  We especially want to
thank Dolan and Nappi for making available to us an early version of
their paper~\cite{dolan98}.

\section{Determinants, differential cocycles and statement of results}
\label{sec-statement-of-results}

\subsection{Background}\label{subsec-history} 
Let $X$ be a Riemann surface of genus $g$.  A {\em theta
characteristic} of $X$ is a square root of the line bundle
$\omega$ of holomorphic $1$-forms.  There are $2^g$ theta
characteristics of $X$, and they naturally form a principal
homogeneous space for the group of square roots of the trivial line
bundle.  Riemann associated to each theta characteristic $x$ a parity
$q(x)$, defined to be the dimension mod 2 of the space of its
holomorphic sections.  He showed that $q(x)$ is invariant under
holomorphic deformations and has remarkable algebraic
properties---namely that
\begin{equation}\label{eq-quadratic}
q(x\otimes a)-q(x) 
\end{equation}
is a quadratic function of $a$, and that the associated bilinear form
\[
e(a,b) = q(x\otimes a\otimes b)-q(x\otimes a)-q(x\otimes b)+q(x)
\]
is non-degenerate and independent of $x$.  
     In these equations, $a$ and $b$ are square roots of the trivial
bundle.  They are classified by elements $H^1(X;\Z/2)$, and the form
$e(a,b)$ corresponds to the cup product.  One can check that the
expression~\eqref{eq-quadratic} depends only on the cohomology class
underlying $a$, and so a choice of theta characteristic thus gives
rise, via Riemann's parity, to a quadratic refinement of the
intersection pairing.  Riemann derived the algebraic properties of the
function $q$ using his $\vartheta$-function and the Riemann
singularities theorem.  In the next section we will deduce these
results from properties of determinant line bundles, in a way
that can be generalized to higher dimensions.

Quadratic functions in topology arise from a famous error, an
unwitting testimony to the depth of the invariants derived from
quadratic functions.  In the 1930's, Pontryagin introduced a geometric
technique for investigating the homotopy groups of
spheres~\cite{pontryagin38_1, pontryagin38_2}.  His method was to
study a map between spheres in terms of the geometry of the inverse
image of a small disk centered on a regular value.  It led eventually
to a remarkable relationship between homotopy theory and differential
topology, and one can find in these papers the beginnings of both
bordism and surgery theories.

Pontryagin's first results concerned the homotopy groups $\pi_{n+1}
S^n$ and $\pi_{n+2} S^n$, for $n\ge 2$
(\cite{pontryagin38_1,pontryagin38_2}).  Because of a subtle error, he
was led to the conclusion%
\footnote{They are actually cyclic of order 2}
\[
\pi_{n+2} S^n =0\qquad n\ge 2.
\]
The groups were later correctly determined by George
Whitehead~\cite{whitehead50} and by Pontryagin
himself~\cite{pontryagin50:_homot} (see
also~\cite{pontryagin59:_smoot} and~\cite{kervaire62:_la_pontr}).
               
The argument went as follows.  A map $f:S^{n+2}\to S^n$ is homotopic
to a generic map.  Choose a small open disk $D\subset S^n$ not
containing any singular values, and let $a\in D$ be a point.  The
subspace $X=f^{-1}(a)\subset S^{n+2}$ is a Riemann surface, the
neighborhood $f^{-1}D$ is diffeomorphic to the normal bundle to the
embedding $X\subset S^{n+2}$, and the map $f$ gives this normal bundle
a framing.  The homotopy class of the map $f$ is determined by this
data (via what is now known as the Pontryagin--Thom construction).  If
$X$ has genus $0$, then the map $f$ can be shown to be null homotopic.
Pontryagin sketched a geometric procedure for modifying $f$ in such a
way as to reduce the genus of $X$ by $1$.  This involved choosing a
simple closed curve $C$ on $X$, finding a disk $D\subset S^{n+2}$
bounding $C$, whose interior is disjoint from $X$, and choosing a
suitable coordinate system in a neighborhood of of $D$.  Pontryagin's
procedure is the basic manipulation of framed surgery.  It is not
needed for the correct evaluation of the groups $\pi_{n+2}S^n$ and,
except for the dimension $0$ analogue, seems not to appear again
until~\cite{milnor61}.

The choice of coordinate system is not automatic, and there is an
obstruction
\[
\phi:H_1(X;\Z/2)\to \Z/2
\]
to its existence.  As indicated, it takes its values in $\Z/2$ and
depends only on the mod $2$ homology class represented by $C$.  As
long $\phi$ takes the value $0$ on a non-zero homology class, the
genus of $X$ can be reduced by $1$.  Pontryagin's error concerned the
algebraic nature of $\phi$, and in~\cite{pontryagin38_2} it was
claimed to be linear.  He later determined that $\phi$ is
quadratic~\cite{pontryagin50:_homot}, and in fact
\[
\phi(C_1+C_2) =
\phi(C_1)+\phi(C_2)+I(C_1,C_2),
\]
where $I(C_1,C_2)$ is the intersection number of $C_1$ and $C_2$.
Thus $\phi$ is a {\em quadratic refinement of the intersection
pairing}.  The Arf invariant of $\phi$ can be used to detect the
non-trivial element of $\pi_{n+2}(S^n)$.  This was the missing
invariant.

Around 1970 Mumford called attention to Riemann's parity, and raised
the question of finding a modern proof of its key properties.  Both
he~\cite{mumford71:_theta} and Atiyah~\cite{atiyah71:_rieman} provided
answers.  It is Atiyah's~\cite[Proposition~4.1]{atiyah71:_rieman} that
relates the geometric and topological quadratic functions.  Atiyah
identifies the set of theta-characteristics with the set of Spin
structures, and Riemann's parity with the mod 2 index of the Dirac
operator.  This gives immediately that Riemann's parity is a
Spin-cobordism invariant, and that the association
\[
(X,x) \mapsto q(x)
\]
is a surjective homomorphism from the cobordism group $\mspin_2$ of
Spin-manifolds of dimension $2$ to $\Z/2$.  Now the map from the
cobordism group $\Omega^{\text fr}_2=\pi_2 S^0$ of (stably) framed
manifolds of dimension $2$ to $\mspin_2$ is an isomorphism, and both
groups are cyclic of order $2$.  It follows that Riemann's parity
gives an {\em isomorphism} $\mspin_2\to\Z/2$, and that the restriction
of this invariant to $\Omega^{\text fr}_2$ coincides with the
invariant of Pontryagin.   The key properties of $q$ can be derived
from this fact.

In~\cite{mumford71:_theta}, Mumford describes an algebraic proof of
these results, and generalizes them to more general sheaves on
non-singular algebraic curves.  In~\cite{harris82:_theta} Harris
extends Mumford's results to the case of singular curves.
                              
\subsection{Determinants and the Riemann parity}\label{subsec-determinants}

In this section we will indicate how the key properties of the
Riemann's $q(x)$ can be deduced using determinants.  Let $E/S$ be a
holomorphic family\footnote{In the language of algebraic geometry,
$E/S$ needs to be a smooth proper morphism of complex analytic
spaces. The manifold $S$ can have singularities.  In this paper, we
will not be dealing with holomorphic structures, and $S$ will
typically have corners.  The notion ``smooth morphism'' will be
replaced with notion of ``neat map'' (see~\ref{sec:manif-with-corn}).}
of Riemann surfaces, and $L$ a holomorphic line bundle over $E$.
Denote by
\[
\det L
\]
the determinant line bundle of the $\dbar$-operator coupled to $L$.
The fiber of $\det L$ over a point $s\in S$ can be identified with
\[
\det H^0(L_s)\otimes \det H^1(L_s)^\ast.
\]
If $K=K_{E/S}$ is the line bundle of relative Kahler
differentials (holomorphic $1$-forms along the fibers), then by Serre
duality this equation can be re-written as
\begin{equation}\label{eq-symmetric-det}
\det H^0(L_s)\otimes \det H^0(K\otimes L_s^{-1}).
\end{equation}
This leads to an isomorphism
\[
\det L \approx \det \left(K\otimes L^{-1}\right),
\]
which, fiberwise, is given by switching the factors
in~\eqref{eq-symmetric-det}.  An isomorphism
$L^2\approx K$ (if one exists) then gives rise to automorphism
\[
\det L\to \det L.
\]
This automorphism squares to the identity, and so is given by a 
holomorphic map
\[
q:S\to\{\pm1\}.
\]
To compute the value of $q$ at a point $s\in S$, write 
\[
(\det L)_s = \Lambda^{2d}\left(H^0(L_s)\oplus(H^0(L_s)\right)
\]
($d=\dim H^0(L_s)$).  The sign encountered in switching the factors is
\[
(-1)^d.
\]
It follows that 
\[
q(s) = (-1)^{\dim H^0(L_s)}
\]
and so we recover Riemann's parity in terms of the symmetry of $\det
L$.  This key point will guide us in higher dimensions.  Note that our
approach shows that {\em Riemann's} $q:S\to \{\pm1 \}$ is holomorphic,
and so invariant under holomorphic deformations.

Riemann's algebraic properties are derived
from the quadratic nature of $\det$, i.e., the fact that the line bundle
\[
B(L_1,L_2) = 
\frac{\det(1)\det(L_1\otimes L_2)}
{\det(L_1)\det(L_2)}
\]
is bilinear in $L_1$ and $L_2$ with respect to tensor product.  We
will not give a proof of this property (see
Deligne~\cite{deligne87:_det} in which the notation $\langle
L_1,L_2\rangle$ is used), but note that it is suggested by
the formula for the first Chern class of $\det L$
\begin{equation}\label{eq-c1-det}
2\,c_1(\det L) = \int_{E/S} (x^2 - x c_1)
\end{equation}
where $x$ is the first Chern class of $L$, and $c_1$ is the first
Chern class of the relative tangent bundle of $E/S$.  

\subsection{Differential cocycles}
\label{subsec-differential-cocycles}

Now suppose $E/S$ has relative dimension $2n$.  In order for the
quadratic term in~\eqref{eq-c1-det} to contribute to the first Chern
class of a line bundle over $S$, $x$ must be an element of
$H^{n+1}(E;\Z)$.  This motivates looking for mathematical objects
classified by $H^{n+1}(E;\Z)$ in the way complex line bundles are
classified by $H^{2}(E;\Z)$.  In the discussion in
\S\ref{subsec-determinants} it was crucial to work with line bundles
and the isomorphisms between them, rather than with isomorphism
classes of line bundles; we need to construct a {\em category} whose
isomorphism classes of objects are classified by $H^{n+1}(E;\Z)$.

\begin{defin}
Let $M$ be a manifold and $n\ge0$ an integer.  The category of {\em
$n$-cocycles}, $\hcat^n(M)$ is the category whose objects are smooth
$n$-cocycles
\[
c \in Z^n(M;\Z)
\]
and in which a morphism from $c_1$ to $c_2$
is an element
\[
b\in C^{n-1}(M;\Z)/\delta\,C^{n-2}(M;\Z)
\]
such that
\[
c_1 + \delta b = c_2 
\]
\end{defin}

There is an important variation involving forms. 

\begin{defin}
Let $M$ be a manifold, and $n\ge0$ an integer.  The category of {\em
differential $n$-cocycles}, $\chcat^n(M)$ is the category whose
objects are triples
\[
(c,h,\omega)\subset C^n(M;\Z)\times C^{n-1}(M;\R)\times \Omega^n(M),
\]
satisfying
\begin{equation}\label{eq-closed}
\begin{gathered}
\delta c=0 \\
d\omega =0\\
\delta h = \omega -c.
\end{gathered}
\end{equation}
A morphism from $(c_1,h_1,\omega_1)$ to $(c_2,h_2,\omega_2)$
is an equivalence class of pairs
\[
(b,k)\in C^{n-1}(M;\Z)\times C^{n-2}(M;\R)
\]
satisfying
\begin{gather*}
c_1 - \delta b = c_2 \\
h_1 + \delta k + b = h_2.
\end{gather*}
The equivalence relation is generated by
\[
(b,k)\sim(b-\delta a, k+\delta k' + a).
\]
\end{defin}

For later purposes it will be convenient to write
\[
d(c,h,\omega) = (\delta c, \omega-c-\delta h, d\,\omega).
\]
Note that $d^2=0$, and that
condition~\eqref{eq-closed} 
says $d(c,h,\omega)=0$.  
We will refer to general triples 
\[
(c,h,\omega)\subset C^n(M;\Z)\times C^{n-1}(M;\R)\times \Omega^n(M),
\]
as {\em differential cochains} (of degree $n$), and those which are
differential cocycles as {\em closed}.  As will be explained in more
detail in \S\ref{subsec-differential-character}, these are the
$n$-cochains and cocycles in the cochain complex
$\chcochain{n}^{\ast}(M)$ with
\[
\chcochain{n}^{k}(M) = \left\{(c,h,\omega)\mid \omega=0 \text{ if }
k < n\right\} 
\subseteq C^{k}(M;\Z)\times C^{k-1}(M;\R)\times \Omega^{n}(M)
\]
and differential given by
\[
d(c,h,\omega)= (\delta c, \omega - c -\delta h, d\omega).
\]
The $k^{\text{th}}$ cohomology group of $\chcochain{n}^{\ast}(M)$ is
denoted 
\[
\chh{n}^{k}(M),
\]
and $\chh{n}^{n}(M),$ can be identified with the group of differential
characters $\chs^{n-1}(M)$ of
Cheeger--Simons~\cite{cheeger85:_differ}.

The operations of addition of cochains and forms define abelian group
structures on the categories $\chcat^n$.  The set of isomorphism
classes of objects in $\hcat^n(M)$ is the group $H^n(M;\Z)$, and the
automorphism group of the trivial object $0$, is $H^{n-1}(M;\Z)$.  The
set of isomorphism classes of objects in $\chcat^n(M)$ is the group
$\chh{n}^{n}(M)$.  The automorphism of the trivial object $(0,0,0)$ is
naturally isomorphic to the group $H^{n-2}(M;\rz)$.  There is a
natural functor $\chcat^n(M)\to\hcat^n(M)$ which is compatible with
the abelian group structures.  On isomorphism classes of objects it
corresponds to the natural map $\chh{n}^{n}(M)\to H^n(M;\Z)$, and on
the automorphism group of the trivial object it is the connecting
homomorphism
\[
H^{n-2}(M;\rz)\to H^{n-1}(M;\Z)
\]
of the long exact sequence associated to 
\[
0\to \Z\to\R\to\rz\to 0.
\]

\begin{eg}\label{eg:6}
The category $\chcat^2(M)$ is equivalent to the category of
$U(1)$-bundles with connection, with the group structure of tensor
product.  One way to present a $U(1)$ bundle is to give for each open
set $U$ of $M$ a principal homogeneous space $\Gamma(U)$ (possibly
empty) for the group of smooth maps $U\to U(1)$.  The points of
$\Gamma(U)$ correspond to local sections of a principal bundle, and
must restrict along inclusions and patch over intersections.  To add a
connection to such a bundle comes down to giving maps
\[
\nabla:\Gamma(U)\to \Omega^1(U)
\]
which are ``equivariant'' in the sense that
\begin{equation}\label{eq-gauge-transformation}
\nabla(g\cdot s) = \nabla(s)+g^{-1}dg.
\end{equation}
An object $x=(c,h,\omega)\in \chcat^2(M)$ gives a $U(1)$-bundle
with connection as follows:  the space $\Gamma(U)$ is the
quotient of the space of
\[
s=(c^1,h^0,\theta^1)\in C^1(U;\Z)\times C^{0}(U;\R)\times \Omega^1(U)
\]
satisfying
\[
ds=x
\]
by the equivalence relation
\[
s\sim s + dt,\qquad t\in C^0(U;\Z)\times \{0\}\times \Omega^0(U).
\]
Any two sections in $\Gamma(U)$ differ by an 
\[
\alpha	\in C^1(U;\Z)\times C^{0}(U;\R)\times \Omega^1(U)
\]
which is closed---in other words, an object in the category
$\chcat^1(M)$.  The equivalence relation among sections corresponds to
the isomorphisms in $\chcat^1(M)$, and so the space $\Gamma(U)$ is a
principal homogeneous space for the group of isomorphism classes in
$\chcat^1(U)$ ie, the group of smooth maps from $U$ to $\rz$.  The
function $\nabla$ associates to $s=(c^1,h^0,\theta^1)$ the $1$-form
$\theta^1$.  The equivariance condition~\ref{eq-gauge-transformation}
is obvious, once one identifies $\rz$ with $U(1)$ and writes the
action multiplicatively.
\end{eg}

Quite a bit of useful terminology is derived by reference
to the above example.  Given a differential cocycle
\[
x=(c,h,\omega)\in \chcat^n(M),
\]
we will refer to $\omega=\omega(x)$ as the {\em curvature} of $x$
, $c=c(x)$
as the {\em characteristic cocycle} (it is a cocycle representing
the $1^{\text{st}}$ Chern class when $n=2$), and 
\[
e^{2\pi i\,h}=e^{2\pi i\,h(x)}
\]
as the {\em monodromy}, regarded as a homomorphism from the group of
$(n-1)$-chains into $U(1)$.  The cohomology class $[x]\in H^n(M;\Z)$
represented by $c$ will be called the {\em characteristic class} of
$x$.  The set of differential cochains
\[
s\in C^{n-1}(M;\Z)\times C^{n-2}(M;\R)\times \Omega^{n-1}(M)
\]
satisfying
\[
ds=x
\]
will be called the space of {\em trivializations} of $x$.  Note that any
differential cochain $s$ is a trivialization of $ds$.

The reduction of $h$ modulo $\Z$
is also known as the {\em differential character} of $x$, so that
\[
\text{monodromy}=e^{2\pi i\,(\text{differential character})}.
\]
The curvature form, the differential character (equivalently the
monodromy) and the characteristic class are invariants of the
isomorphism class of $x$, while the characteristic cocycle is not.  In
fact the differential character (equivalently the monodromy)
determines $x$ up to isomorphism in $\chcat^n(M)$.

It is tempting to refer to the form component of a trivialization
$s=(c,h,\theta)$ as the ``curvature,'' but this
does not reduce to standard terminology.  By analogy with the case in
which $s$ has degree 1, we will refer to $\theta=\nabla(s)$ as the
{\em connection form} associated to $s$.

It is useful to spell this out in a couple of other cases.  The category
$\chcat^1(M)$ is equivalent to the category whose objects are smooth
maps from $M$ to $\rz$, and with morphisms, only the identity maps.
The correspondence associates to $(c^1,h^0,\omega^1)$ its differential
character.  From the point of view of smooth maps to $\rz$, the
curvature is given by the derivative, and the characteristic cocycle
is gotten by pulling pack a fixed choice of cocycle representing the
generator of $H^1(\rz;\Z)$.  The characteristic class describes the
effect of $f:M\to\rz$ in cohomology, and determines $f$ up to
homotopy.  A {\em trivialization} works out to be a lift of $f$ to $\R$.  If
the trivialization is represented by $(c^0,0,\theta^0)$, then the lift is
simply given by $\theta^0$----the connection form of the trivialization.
Finally, as the reader will easily check, $\chcat^0(M)$ is equivalent
to the category whose objects are maps from $M$ to $\Z$, and morphisms
the identity maps.

\subsection{Integration and \orientation{s}}
\label{sec:integr-orient}

Let $M$ be a smooth compact manifold and $V\to M$ a (real) vector
bundle over $M$ of dimension $k$.  A {\em differential Thom cocycle}
on $V$ is a ({\em compactly supported}) cocycle
\[
U=(c,h,\omega)\in \chcocycle{k}^k_{c}(V)
\]
with the property that for each $m\in M$,
\[
\int_{V_{m}}\omega= \pm 1.
\]
A choice of differential Thom cocycle determines an orientation of each $V_m$
by requiring that the sign of the above integral be $+1$.

The integral cohomology class underlying a differential Thom cocycle
is a Thom class $[U]$ in $H^k_c(V;\Z)$.  There is a unique $[U]$
compatible with a fixed orientation, and so any two choices of $U$
differ by a cocycle of the form $\delta b$, with
\[
b=(b,k,\eta)\in \chcochain {k-1}^{k-1},
\]

Using the ideas of 
Mathai--Quillen~\cite{quillen85:_super_chern,mathai86:_super_thom},
a differential Thom cocycle%
\footnote{The Thom form constructed by Mathai--Quillen is not
compactly supported, but as they remark
in~\cite{mathai86:_super_thom,quillen85:_super_chern} a minor
modification makes it so.}
can be associated to a metric and connection on $V$, up to addition of
a term $d\,(b,k,0)$.

\begin{defin}\label{def:5}
A {\em \orientation{} of $p:E\to S$} consists of the
following data
\begin{textList}
\item  A smooth embedding $E\subset S\times \R^N$ for some $N$;
\item  A tubular neighborhood\footnote{\label{fn:1}A {\em tubular neighborhood} of
$p:E\hookrightarrow S$ is a vector bundle $W$ over $E$, and an
extension of $p$ to a diffeomorphism of $W$ with a neighborhood of
$p(E)$.  The derivative of the embedding $W\to M$ gives a vector
bundle isomorphism $W\to \nu_{E/S}$ of $W$ with the normal bundle to
$E$ in $S$.}
\[
W\subset S\times \R^{N};
\]
\item A differential Thom cocycle $U$ on $W$.  
\end{textList}

An {\em \oriented{} map} is a map $p:E\to S$ together with
a choice of \orientation.
\end{defin}

While every map of compact manifolds $p:E\to S$ factors through an
embedding $E\subset S\times\R^{N}$, not every embedding $E\subset
S\times\R^{N}$ admits a tubular neighborhood.  A necessary and
sufficient condition for the existence of a tubular neighborhood is
that $p:E\to S$ be {\em neat}.  A {\em neat map} of manifolds with
corners is a map carrying corner points of codimension $j$ to corner
points of codimension $j$, and which is transverse (to the corner) at
these points.  A neat map $p:E\to S$ of compact manifolds factors
through a neat embedding $E\subset S\times\R^{N}$, and a neat
embedding has a normal bundle and admits a tubular neighborhood.
Every smooth map of closed manifolds is neat, and a map $p:E\to S$ of
manifolds with boundary is neat if $f\circ p$ is a defining function
for the boundary of $E$ whenever $f$ is a defining function for the
boundary of $S$.  See Appendix~\ref{sec:manif-with-corn}.

Let $E/S$ be an \oriented{} map of relative dimension $k$.  In
\S\ref{subsec-integration}
we will define natural, additive integration functors
\begin{gather*}
\int_{E/S}:\hcat^n(E)\to\hcat^{n-k}(S)\\
\int_{E/S}:\chcat^n(E)\to\chcat^{n-k}(S),
\end{gather*}
compatible with the formation of the ``connection form,'' and hence
curvature:
\begin{align*}
\omega\left(\int_{E/S}x\right) &= \int_{E/S}\omega\left(x\right) \\
\nabla\left(\int_{E/S}s\right) &= \int_{E/S}\nabla\left(s\right).
\end{align*}
In the special case when $E/S$ is a fibration over an open dense set,
the integration functors arising from different choices of
\orientation{} are naturally isomorphic (up to the usual sign).

When $E$ is a manifold with boundary and $S$ is a closed manifold, a
map $p:E\to S$ cannot be neat, and so cannot admit an \orientation{},
even if $E/S$ is orientable in the usual sense.  Integration can be
constructed in this case by choosing a defining function $f:E\to[0,1]$
for the boundary of $E$.  An orientation of $E/S$ {\em can} be
extended to an \orientation{} of
\[
E\to [0,1]\times S.
\]
Integration along this map is then defined, and gives a functor
\[
\int_{E/[0,1]\times S}:\chcat^n(E)\to\chcat^{n-k+1}
([0,1]\times S),\qquad (k=\dim E/S),
\]
which commutes with the differential $d$, and is compatible with
restriction to the boundary:
\[
\begin{CD}
\chcat^n(E)@>\int_{E/S}>>\chcat^{n-k}([0,1]\times S) \\
@V\text{res.}VV @VV\text{res.}V \\
\chcat^n(\partial E)@>>\int_{\partial E/S}>\chcat^{n-k}(\{0,1\}\times S).
\end{CD}
\]
The usual notion of integration over manifolds with boundary
is the iterated integral
\[
\int_{E/S}x=\int_0^{1}\int_{E/[0,1]\times S}x.
\]
The second integral does not commute with $d$ (and hence 
does not give rise to a functor), but rather satisfies Stokes
theorem:
\begin{equation}\label{eq-fundamental-calculus}
d\int_0^{1}s=
\int_0^{1} ds - (-1)^{|s|}
\left(s\vert_{\{1\}\times S}-s\vert_{\{0\}\times S}\right).
\end{equation}
So, for example, if $x$ is closed, 
\begin{equation}\label{eq-integral-boundary}
d\,\int_0^1\int_{E/[0,1]\times S}x = (-1)^{n-k}\int_{\partial E/S}x;
\end{equation}
in other words
\[
(-1)^{n-k}\int_0^1\int_{E/[0,1]\times S}x
\]
is a trivialization of
\[
\int_{\partial E/S}x.
\]

In low dimensions these integration functors work out to be familiar
constructions.  If $E\to S$ is an oriented family of $1$-manifolds
without boundary, and $x\in\chcat^{2}(E)$ corresponds to a line bundle
then
\[
\int_{E/S}x
\]
represents the function sending $s\in S$ to the monodromy of $x$
computed around the fiber $E_s$.  

The cup product of cocycles, and the wedge product of forms combine to
give bilinear functors
\[
\chcat^{n}(M)\times
\chcat^{m}(M)\to
\chcat^{n+m}(M)
\]
which are compatible with formation of the connection form (hence
curvature) and characteristic cocycle.  See
\S\ref{subsec-differential-character}.

\subsection{Integral Wu-structures}\label{subsec-wustruct} 
As we remarked in the introduction, the role of the canonical bundle
is played in higher dimensions by an integral Wu-structure.   

\begin{defin}
Let $p:E\to S$ be a smooth map, and fix a cocycle $\nu\in
Z^{2k}(E;\Z/2)$ representing the Wu-class $\nu_{2k}$ of the relative
normal bundle.  A {\em differential integral Wu-structure of degree
$2k$} on $E/S$ is a differential cocycle
\[
\wustruct=(c,h,\omega)\in\chcocycle{2k}^{2k}(E)
\]
with the property that $c\equiv\nu\mod 2$.
\end{defin}

We will usually refer to a differential integral Wu-structure of
degree $2k$ as simply an {\em integral Wu-structure}.  If $\wustruct$
and $\wustruct'$ are integral Wu-structures, then there is a unique
\[
\eta\in\chcochain{2k}^{2k}(E)
\]
with the property that 
\[
\wustruct'=\wustruct+2\eta.
\]

We will tend to overuse the symbol $\lambda$ when referencing integral
Wu-structures.  At times it will refer to a differential cocycle, and
at times merely the underlying topological cocycle (the $c$
component).  In all cases it should  be clear from context which meaning
we have in mind.  

For a cocycle $\nu\in Z^{2k}(E;\Z/2)$, let $\chcat^{2k}_{\nu}(E)$
denote the category whose objects are differential cocycles
$x=(c,h,\omega)$ with $c\equiv\nu\mod 2$, and in which a morphism
from $x$ to $x'$ is an equivalence class of differential cochains
\[
\tau = (b,k)\in \chcochain{2k}^{2k-1}(E)/d\left(\chcochain{2k}^{2k-2}(E)\right)
\]
for which $x'=x+2d\,\tau$.  The set of isomorphism classes in
$\chcat^{2k}_{\nu}(E)$ is a torsor (principal homogeneous space) for
the group $\chh{2k}^{2k}(M)$, and the automorphism group of any object
is $H^{2k-2}(E;\rz)$.  We will write the action of
$\mu\in\chh{2k}^{2k}(M)$ on $x\in\chcat^{2k}_{\nu}(E)$ as
\[
x\mapsto x + (2)\mu,
\]
with the parentheses serving as a reminder that the multiplication by
$2$ is formal;  even if $\mu\in\chh{2k}^{2k}(M)$ has order $2$, the
object $x + (2)\mu$ is not necessarily isomorphic to $x$.%
\footnote{It could happen that $2\mu$ can be written as $d(b',k')$,
but not as $2d(b,k)$.}
A differential cochain 
\[
y=(c',h',\omega')\in \chcochain{2k}^{2k}
\]
gives a functor
\begin{align*}
\chcat^{2k}_{\nu}(E) & \to
\chcat^{2k}_{\nu+d\bar y}(E)\\
x &\mapsto x + dy
\end{align*}
in which $\bar y$ denotes the mod $2$ reduction of $c'$.  Up to natural
equivalence, this functor depends only on the value of $\bar y$.  In
this sense the category $\chcat^{2k}_{\nu}(E)$ depends only on the
cohomology class of $\nu$.

With this terminology, a differential integral Wu-structures is an
object of the category $\chcat^{2k}_{\nu}(E)$, with $\nu$ a cocycle
representing the Wu-class $\nu_{2k}$.  An {\em isomorphism} of
integral Wu-structures is an isomorphism in $\chcat^{2k}_{\nu}(E)$.

We show in Appendix~\ref{sec:integer-wu-classes} that it is possible to
associate an integral Wu-structure $\nubar_{2k}(E/S)$ to a
$\spin$ structure on the relative normal bundle of $E/S$.
Furthermore, as we describe in
\S\ref{sec:char-class}, a connection on $\nu_{E/S}$ gives a refinement
of $\nubar_{2k}(E/S)$ to a differential integral Wu-structure
$\nuhat_{2k}(E/S)\in\chh{2k}_{\nu}^{2k}(E)$.  
We'll write
$\wustruct(s,\nabla)$ (or just $\wustruct(s)$) for the integral
Wu-structure associated to a $\spin$-structure $s$ and connection
$\nabla$.  By Proposition~\ref{thm:4}, if the $\spin$-structure is
changed by an element $\alpha\in H^1(E;\Z/2)$ then, up to isomorphism,
the integral Wu-structure changes according to the rule
\begin{equation}\label{eq-change-of-spin}
\wustruct(s+\alpha) \equiv \wustruct(s)+(2)\beta\,\sum_{\ell\ge
0}\alpha^{2^\ell-1}\nu_{2k-2^{\ell}}
\,\in\,\chcat^{2k}_{\nu}(E),
\end{equation}
where $\nu_j$ is the $j^{\text{th}}$ Wu class of the relative normal
bundle, and $\beta$ denotes the map 
\[ 
H^{2k-1}(E;\Z/2)\to H^{2k-1}(E;\rz)\hookrightarrow\chh{2k}^{2k}(E)
\]
described in \S\ref{subsec-differential-character}.

We now reformulate the above in terms of ``twisted differential
characters.''  These appear in the physics
literature~\cite{witten01:_overv_k} as Chern-Simons terms associated
to characteristic classes which do not necessarily take integer values
on closed manifolds.  We will not need this material in the rest of
the paper.

\begin{defin}
Let $M$ be a manifold,
\[
\nu\in Z^{2k}(M;\rz)
\]
a smooth cocycle, and and $2k\ge0$ an integer.  The category of {\em
$\nu$-twisted differential $2k$-cocycles}, $\chcat^{2k}_{\nu}(M)$, is the
category whose objects are triples
\[
(c,h,\omega)\subset C^{2k}(M;\R)\times C^{2k-1}(M;\R)\times \Omega^{2k}(M),
\]
satisfying
\[
\begin{aligned}
c&\equiv\nu\mod\Z \\
\delta c&=0 \\
d\omega &=0\\
\delta h &= \omega -c
\end{aligned}
\]
(we do not distinguish in notation between a form and the cochain
represented by integration of the form over chains).  A morphism from
$(c_1,h_1,\omega_1)$ to $(c_2,h_2,\omega_2)$ is an equivalence class
of pairs
\[
(b,k)\in C^{2k-1}(M;\Z)\times C^{2k-2}(M;\R)
\]
satisfying
\begin{gather*}
c_1 - \delta b = c_2 \\
h_1 + \delta k + b = h_2.
\end{gather*}
The equivalence relation is generated by
\[
(b,k)\sim(b-\delta a, k+\delta k' + a).
\]
\end{defin}
The category $\chcat^{2k}_{\nu}(M)$ is a torsor for the category
$\chcat^{2k}(M)$.  We will write the translate of
$v\in\chcat^{2k}_{\nu}(M)$ by $x\in\chcat^{2k}(M)$ as $v+x$.

\begin{rem}
We have given two meanings to the symbol $\chcat^{2k}_{\nu}(M)$; one
when $\nu$ is a cocycle with values in $\Z/2$, and one when $\nu$ is a
cocycle taking values in $\rz$.  For a cocycle $\nu\in
Z^{2k}(M;\Z/2)$ there is an {\em isomorphism} of categories
\begin{align*}
\chcat^{2k}_{\nu}(M) &\to \chcat^{2k}_{\frac12\nu}(M) \\
(c,h,\omega)&\mapsto \frac12(c,h,\omega),
\end{align*}
in which $\frac12\nu\in Z^{2k}(M;\rz)$ is the composite of $\nu$ with
the inclusion 
\begin{align*}
\Z/2 &\hookrightarrow \rz \\
t &\mapsto \frac12t
\end{align*}
\end{rem}

The isomorphism class of a $\nu$-twisted differential cocycle
$(c,h,\omega)$ is determined by $\omega$ and the value of $h$ modulo
$\Z$.
\begin{defin}
Let $M$ be a smooth manifold, and
\[
\nu\in Z^{2k}(M;\rz)
\]
a smooth cocycle.  A {\em $\nu$-twisted differential character} is
a pair
\[
(\chi,\omega)
\]
consisting of a character
\[
\chi:Z_{2k-1}(M)\to\rz
\]
of the group of smooth $(2k-1)$-cycles, and an $2k$-form $\omega$
with the property that for every smooth $2k$-chain $B$, 
\[
\int_{B}\omega - \chi(\partial B) \equiv \nu(B) \mod\Z.
\]
The set of $\nu$-twisted differential characters will be denoted
\[
\chs^{2k-1}_{\nu}(M).
\]
It is torsor for the group $\chs^{2k-1}(M)$.
\end{defin}

\subsection{The main theorem}\label{subsec-main-result}

\begin{thm}\label{thm-main}
Let $E/S$ be an \oriented{} map of manifolds of relative dimension
$4k-i$, with $i\le2$.  Fix a differential cocycle $L_{4k}$ refining
the degree $4k$ component of the Hirzebruch $L$ polynomial, and fix a
cocycle $\nu\in Z^{2k}(E;\Z/2)$ representing the Wu-class $\nu_{2k}$.
There is a functor
\[
\kappa_{E/S}:\chcat^{2k}_{\nu}(E) \to \chcat^{i}(S)
\]
with the following properties:
\begin{thmList}
\item\label{item:2} (Normalization) modulo torsion,
\[
\kappa(\wustruct)\approx \frac18
\int_{E/S}\wustruct\cup\wustruct - \check{L}_{4k};
\]
\item\label{item:3} (Symmetry) there is an isomorphism 
\[
\tau(\wustruct):\kappa(-\wustruct)\xrightarrow{\approx}{}\kappa(\wustruct)
\]
satisfying 
\[
\tau(\wustruct)\circ\tau(-\wustruct)=\text{identity map of
$\kappa(\wustruct)$};
\]

\item\label{item:4} (Base change) Suppose that $S'\subset S$ is a
closed submanifold and $p:E\to S$ is transverse to $S'$.  Then the
\orientation of $E/S$ induces an \orientation of
\[
E'=p^{-1}(S')\xrightarrow{p'}{} S'
\]
and, with $\tilde f$ denoting the map $E'\to E$, there is an isomorphism
\[
\tilde f^{\ast} \kappa_{E/S}(\wustruct)\approx
\kappa_{E'/S'}(f^{\ast}\wustruct);
\]

\item\label{item:5} (Transitivity) let
\[
E\to B\to S
\]
be a composition of \oriented{} maps, and suppose given a framing of
the stable normal bundle of $B/S$ which is compatible with its
\orientation.  The \orientation{s} of $E/B$ and $B/S$ combine to give
an \orientation{} of $E/S$, the differential cocycle $\check L_{4k}$
represents the Hirzebruch $L$-class for $E/S$, and a differential
integral Wu--structure $\wustruct$ on $E/B$ can be regarded as a
differential integral Wu--structure on $E/S$.  In this situation there
is an isomorphism
\[
\kappa_{E/S}(\wustruct) \approx \int_{B/S} \kappa_{E/B}(\wustruct);
\]
\end{thmList}
\end{thm}

We use $\kappa$ and a fixed differential integral Wu-structure
$\lambda$ to construct a quadratic functor
\[
q_{E/S}=q_{E/S}^{\wustruct}:\chcat^{2k}(E)\to \chcat^i(S)
\]
by
\[
q_{E/S}^{\wustruct}(x)=\kappa(\wustruct-2x).
\]

\begin{cor}\label{thm:28}
\begin{thmList}
\item\label{item:6} The functor $q_{E/S}$ is a quadratic refinement of
\[
B(x,y)=\int_{E/S}x\cup y
\]
ie, there is an isomorphism
\[
t(x,y):q_{E/S}(x+y)-q_{E/S}(x)-q_{E/S}(y)+q_{E/S}(0)
\xrightarrow{\approx}{} B(x,y);
\]
\item\label{item:7} (Symmetry) There is an
isomorphism
\[
s^{\wustruct}(x):q_{E/S}^{\wustruct}(\lambda-x) \approx q_{E/S}^{\wustruct}(x)
\]
satisfying 
\[
s^{\wustruct}(x)\circ s^{\wustruct}(\lambda-x)=\text{identity map};
\]
\item\label{item:8} (Base change) The functor $q$ satisfies the base change
property~\thmItemref{item:4} of Theorem~\ref{thm-main}:
\[
q_{E'/S'}(\tilde f^\ast x)\approx
f^\ast q_{E/S}(x);
\]

\item\label{item:9} (Transitivity) The functor $q$ satisfies the transitivity
property~\thmItemref{item:5} of 
Theorem~\ref{thm-main}:
\[
q_{E/S}(x) \approx \int_{B/S} q_{E/B}(x).
\]
\end{thmList}
\end{cor}

\begin{rem}
All the isomorphisms in the above are natural isomorphisms, ie
they are isomorphisms of {\em functors}.
\end{rem}

\begin{rem}\label{rem:7}
Suppose that $E/S$ is a fibration of relative dimension $n$, with Spin
manifolds $M_{s}$ as fibers, and $E\subset \R^{N}\times S$ is an
embedding compatible with the metric along the fibers.  Associated to
this data is an \orientation{}, differential integral Wu-structure,
and a differential $L$-cocycle.  To describe these we need the results
of \S\ref{sec:char-class} where we show that a connection on a
principal $G$-bundle and a classifying map determine differential
cocycle representatives for characteristic classes, and a differential
Thom cocycle for vector bundles associated to oriented orthogonal
representations of $G$.

Let $W$ be the normal bundle to
the embedding $E\subset S\times\R^{N}$.  The fiber at $x\in
M_{s}\subset E$ of the normal bundle $W$ of $E\subset S\times\R^{N}$
can be identified with the orthogonal complement in $\R^{N}$ of the
tangent space to $M_{s}$ at $x$.  Because of this, $W$ comes equipped
with a metric, a connection, a $\spin$-structure, and a classifying map
\[
\begin{CD}
W @>>> \xi^{\text{Spin}}_{N-n} \\
@VVV @VVV \\
E @>>> \spingrass{N-n}{\R^{N}},
\end{CD}
\]
where $\spingrass{N-n}{\R^{N}}$ is defined by the homotopy pullback
square
\[
\begin{CD}
\spingrass{N-n}{\R^{N}}@>>> \bspin(N-n) \\
@VVV @VVV \\
\sgrass{N-n}{\R^{N}} @>>> BSO(N-n).
\end{CD}
\]
In the above, $\xi^{\text{Spin}}_{N-n}$ is the universal $(N-n)$-plane
bundle, and $\sgrass{N-n}{\R^{N}}$ refers to the oriented Grassmannian.
The metric on $W$ can be used to construct a tubular neighborhood 
\[
W\subset S\times\R^{N},
\]
and the differential Thom cocycle and $L$-cocycle are the ones constructed
in \S\ref{sec:char-class}.

Fix once and for all a cocycle
\[
\tilde c \in Z^{2k}\left(\spingrass{N-n}{\R^{N}};\Z \right)
\]
representing $\nubar_{2k}$.  The the mod $2$ reduction of $\tilde c$
represents the universal Wu-class $\nu_{2k}$.  We set 
\begin{align*}
c &= f^{\ast}(\tilde c) \\
\nu &=  c \mod 2 \in Z^{2k}(E;\Z/2).
\end{align*}
We take  $\lambda$ to be the differential cocycle
associated to characteristic class $\nubar_{2k}$ and
the $\spin$-connection on the normal bundle to $E\subset
\R^{N}\times S$.
\end{rem}

\begin{rem}
To be completely explicit about the \orientation{} constructed in
Remark~\ref{rem:7} a convention would need to be chosen for
associating a tubular neighborhood to a (neat) embedding equipped with
a metric on the relative normal bundle.  It seems best to let the
geometry of a particular situation dictate this choice.  Different
choices can easily be compared by introducing auxiliary parameters.
\end{rem}

\begin{rem}
Another important situation arises when $K\to S$ is a map from a
compact manifold with boundary to a compact manifold without boundary.
As it stands, Theorem~\ref{thm-main} doesn't apply
since $K\to S$ is not a neat map.  To handle this case choose a
defining function $f$ for $\partial K=E$, and suppose we are given an
orientation and integral Wu--structure $\wustruct$ of $K/[0,1]\times
S$.  Then by~\eqref{eq-fundamental-calculus}, for each differential
cocycle $y\in\chcat^{2k}(K)$
\begin{equation}\label{eq-a}
(-1)^{i+1}\, d\int_0^1q_{K/[0,1]\times S}^{\wustruct}(y) 
= q_{K/[0,1]\times S}^{\wustruct}(y)\vert_{S\times\{1\}}-
q_{K/[0,1]\times S}^{\wustruct}(y)\vert_{S\times\{0\}}.
\end{equation}
Since
\[
\begin{CD}
E @>>> K\\
@VVV @VVV \\
S\times\{0,1\}@>>> [0,1]\times S
\end{CD}
\]
is a transverse pullback square, the base change property of $q$
gives a canonical isomorphism
\begin{equation}\label{eq-b}
q_{K/[0,1]\times S}^{\wustruct}(y)\vert_{S\times\{1\}} 
-q_{K/[0,1]\times S}^{\wustruct}(y)\vert_{S\times\{0\}} \approx
q_{E/S}^{\wustruct}(y\vert_{E}).
\end{equation}
Combining~\eqref{eq-a} and~\eqref{eq-b}, and writing $x=y\vert_{E}$, gives
\begin{equation}\label{eq-q-boundary}
(-1)^{i+1}\,d\int_0^1q_{K/[0,1]\times S}^{\wustruct}(y) = 
q_{E/S}^{\wustruct}(x).
\end{equation}
In other words, {\em writing the pair $(E/S,x)$ as a boundary gives a
trivialization of $q_{E/S}^{\wustruct}(x)$}.
\end{rem}

Theorem~\ref{thm-main} remains true when $i>2$, but with the
modification that $q_{E/S}^\wustruct$ is not a differential cocycle,
but a differential function to $(\tilde I_{i};\iota)$.  See
\S\ref{sec-topology} for more details.

The properties of $q^{\wustruct}_{E/S}(x)$ lead to an explicit formula
for $2q^{\wustruct}_{E/S}(x)$.  First of all, note that 
\[
q^{\wustruct}_{E/S}(\lambda)=\kappa(-\lambda)=\kappa(\lambda).
\]
It follows that
\begin{align*}
q^{\wustruct}_{E/S}(0) &=
q^{\wustruct}_{E/S}(\lambda) =q^{\wustruct}_{E/S}(x+(\lambda-x)) \\
&=q^{\wustruct}_{E/S}(x)+q^{\wustruct}_{E/S}(\lambda-x)-q^{\wustruct}_{E/S}(0)+B(x,\lambda-x) \\
\end{align*}
Using the symmetry
$q_{E/S}^{\wustruct}(\lambda-x)=q_{E/S}^{\wustruct}(x)$, and
collecting terms gives
\begin{cor}\label{cor-2q-formula}	
Let $E/S$ be as in Theorem~\ref{thm-main}.  There is an isomorphism
\[
2q_{E/S}^{\wustruct}(x)\approx 2\kappa_{E/S}(\lambda) + \int_{E/S}x\cup
x-x\cup\lambda.\qed
\]
\end{cor}

We now exhibit the effect of $q$ on automorphisms in the first
non-trivial case, $i=2$.  We start with
\[
x=0\in\chcat^{2k}(E).
\]
The automorphism group of $x$ is $H^{2k-2}(E;\rz)$, and the
automorphism group of $q(x)=\kappa(\lambda)\in\chcat^2(S)$ is
$H^0(S;\rz)$.  The map
\[
q:H^{2k-2}(E;\rz)\to H^0(S;\rz)
\]
is determined by its effect on fibers
\[
q:H^{2k-2}(E_s;\rz)\to H^0(\{s\};\rz)=\rz\qquad s\in S.
\]
By Poincar\'e duality, such a homomorphism is given by
\[
q(\alpha)=\int_{E_s}\alpha\cup \mu \qquad \alpha\in\Aut(0),
\]
for some $\mu\in H^{2k}(E_s)$.  After some work (involving either
chasing through the definitions, or replacing automorphisms with their
mapping tori in the usual way), the quadratic nature of $q$ implies
that for general $x\in\chcat^{2k}(E_s)$ one has
\begin{equation}\label{eq-aut-q}
q(\alpha)=\int_{E_s}\alpha\cup([x]+\mu)\qquad \alpha\in\Aut(x).
\end{equation}
Again, with a little work, Corollary~\ref{cor-2q-formula} gives the
identity
\[
2\mu=-[\lambda].
\]

This is a somewhat surprising formula.  The theory of Wu-classes shows
that for each $2n$-manifold $M$, the Wu class $\nu_{2k}$ of the
normal bundle vanishes, and hence the cohomology class $[\lambda]$ is divisible
by $2$.  Our theory shows that the choice of integral Wu-structure
actually leads to a canonical way of dividing $[\lambda]$ by $2$,
namely $-\mu$.

The cohomology class $\mu$ arises in another context.  Suppose that
$E/S$ is the product family $E=M\times S$, and the integral
Wu-structure is gotten by change of base from an integral
Wu-structure on $M\to\text{\{pt\}}$.  If $x\in\chcat^{2k}(E)$ is
arbitrary, and $y\in\chcat^{2k}(E)$ is of the form $(0,h,0)$, then
\begin{equation}\label{eq-mu-again}
q(x+y) \approx q(x) + \int_{E/S}y\cup(x+\mu\otimes 1),
\end{equation}
where $\mu\in H^{2k}(M;\Z)$ is the class described above.

\subsection{The fivebrane partition function} \label{subsec-witten} We
will now show how Corollary~\ref{thm:28} can be used to construct the
holomorphic line bundle described by Witten
in~\cite{witten97:_five_m}.  We will need to use some of the results
and notation of \S\ref{subsec-differential-character}, especially the
exact sequences~\eqref{eq-sequences}.

Suppose that $M$ is a Riemannian $\spin$ manifold of dimension $6$,
and let $J$ be the torus
\[
J=H^3(M)\otimes\rz.
\]
The Hodge $\ast$-operator on $3$-forms determines a complex structure
on $H^3(M;\R)$, making $J$ into a complex torus.  In fact $J$ is a
polarized abelian variety---the $(1,1)$-form $\omega$ is given by the
cup product, and the Riemann positivity conditions follow easily from
Poincar\'e duality.  These observations are made by Witten
in~\cite{witten97:_five_m}, and he raises the question of finding a
symmetric, holomorphic line bundle $\lineb$ on $J$ with first Chern
form $\omega$.  There is, up to scale, a unique holomorphic section of
such an $\lineb$ which is the main ingredient in forming the fivebrane
partition function in $M$-theory.

Witten outlines a construction of the line bundle $\lineb$ in case
$H^4(M;\Z)$ has no torsion.  We construct the line bundle $\lineb$ in
general.  As was explained in the introduction this was one motivation
for our work.

The group $J$ is naturally a subgroup of the group
\[
\chh 4^4(M)=\chs^{4-1}(M)
\]
of differential characters (\S\ref{subsec-differential-character}) of
degree $3$ on $M$.  In principle, the fiber of $\lineb$ over a point
$z\in J$ is simply the line
\[
q^{\wustruct}_{M/\text{pt}}(x)=q^{\wustruct}(x),
\]
where $x$ is a differential cocycle representing $z$, and $\wustruct$
is the integral Wu-structure corresponding to the $\spin$-structure on
$M$.  The key properties of $\lineb$ are a consequence of
Theorem~\ref{thm-main}.  One subtlety is that the line
$q^{\wustruct}(x)$ need not be independent of the choice of
representing cocycle.  Another is that the symmetry is not a symmetry
about $0$.  

To investigate the dependence on the choice of cocycle representing $z$,
note that any two representatives $x$ and $x'$ are isomorphic, and
so the lines $q^{\wustruct}(x)$ and $q^{\wustruct}(x')$ are also
isomorphic.  Any two isomorphisms of $x$ and $x'$ differ by an
automorphism $\alpha\in\Aut(x)$, and by~\eqref{eq-aut-q} the
corresponding automorphism $q^{\wustruct}(\alpha)$ of
$q^{\wustruct}(x)$ is given by
\[
q^{\wustruct}(\alpha)=\int_{M}\alpha\cup([z]+\mu).
\]
This means that as long as $[z]=-\mu$ the isomorphism
\[
q^{\wustruct}(x)\approx q^{\wustruct}(x')
\]
is independent of the choice of isomorphism $x\approx x'$, and the
lines $q^{\wustruct}(x)$ and $q^{\wustruct}(x')$ can be canonically
identified.  The desired $\lineb$ is therefore naturally found over
the ``shifted'' torus $J_t$, where $t\in A^4(M)$ is any element whose
underlying cohomology class is $-\mu$, and
\[
J_t\subset \chh{4}^{4}(M)
\]
denotes the inverse image of $t$ under the second map of the sequence
\[
H^3(M)\otimes \rz \to \chh{4}^{4}(M)\to A^4(M).
\]
The sequence shows that each space $J_t$ is a principal homogeneous
space for $J$.

To actually construct the line bundle $\lineb$ over $J_t$, we need to
choose a universal $Q\in \chh{4}^{4}(M\times J_t)$ (an analogue of a
choice of Poincar\'e line bundle).  We require $Q$ to have the
property that for each $p\in J_t$, the restriction map
\[
(1\times\{p\})^\ast:\chh{4}^{4}(M\times J_t)\to
\chh{4}^{4}(M\times \{p\})
\]
satisfies
\begin{equation}\label{eq-poincare-property}
(1\times\{p\})^\ast Q= p.
\end{equation}
The line bundle $\lineb$ will then turn out to be
\[
q(x)=q^{\wustruct}_{M\times J_t/J_t}(x)
\]
where $x=(c,h,G)$ is a choice of differential cocycle representing
$Q$.
We first construct
\[
Q_0\in\chh{4}^{4}(M\times J)
\]
and then use a point of $J_t$ to translate it to 
\[
Q\in\chh{4}^{4}(M\times J_t).
\]

Write 
\begin{align*}
V&=H^3(M;\R)\\
L&=H^3(M;\Z)
\end{align*}
so that $V$ is the tangent space to $J$ at the origin, and there is a
canonical isomorphism
\[
V/L=J.
\]
Let
\[
\theta\in\Omega^3(M\times V)
\]
be the unique $3$-form which is vertical with respect to the
projection
\[
M\times V\to V
\]
and whose fiber over $v\in V$ is $\harm(v)$---the unique harmonic form
whose deRham cohomology class is $v$.  The form $\theta$ gives rise to
an element
\[
\bar\theta\in\chh{4}^{4}(M\times V)
\]
via the embedding
\[
\Omega^3/\Omega^3_0 \to \chh{4}^{4}(M\times V).
\]
The restriction map
\[
(1\times\{v\})^\ast:\chh{4}^{4}(M\times V)\to
\chh{4}^{4}(M\times \{v\})
\]
has the property that for each $v\in V$,
\begin{equation}\label{eq-poincare-property-v}
(1\times\{v\})^\ast \bar\theta= v.
\end{equation}
This guarantees that~\eqref{eq-poincare-property} will hold for
the class $Q$ we are constructing.

We will now show that there is a class $Q_{0}\in\chh{4}^{4}(M\times J)$
whose image in $\chh{4}^{4}(M\times V)$ is $\bar\theta$.  Note that the
curvature ($d\,\bar\theta$) of $\bar\theta$ is translation invariant,
and is the pullback of a translation invariant four form
\[
\omega(Q_{0})\in\Omega^4_0(M\times J)
\]
with integral periods.  Consider the following diagram in which
the rows are short exact:
\[
\begin{CD}
H^3(M\times V;\R/\Z) @>>> \chh{4}^{4}(M\times V) @>>>
\Omega^4_0(M\times V)\\
@AAA @AAA @AAA \\
H^3(M\times J;\R/\Z) @>>> \chh{4}^{4}(M\times J) @>>>
\Omega^4_0(M\times J).
\end{CD}
\]
The leftmost vertical map is surjective.  We've already seen that
the curvature of $\bar\theta$ descends to $M\times J$.
An easy diagram chase 
gives the existence of the desired $Q_0\in\chh{4}^{4}(M\times J)$.
Such a class $Q_0$ is unique up to the addition of an element
\[
a\in\ker H^3(M\times J;\R/Z) \to H^3(M\times V;\R/\Z) = 
H^3(M;\R/\Z).
\]

When $t\in A^4(M\times J)$ is an element whose underlying cohomology
class is $-\mu\otimes 1$ the line bundle $q(x)$ is independent of both
the choice of $Q$ and the choice of $x$ representing $Q$.  The
independence from the choice of $x$ amounts, as described above, to
checking that $q$ sends automorphisms of $x$ to the identity map of
$q(x)$.  If $\alpha\in H^{2}(M\times J)$ is an automorphism of $x$,
then by~\eqref{eq-aut-q}
\[
q(\alpha)=\int_{M\times J_t/J_t}\alpha\cup([x]+\mu\otimes 1).
\]
The integral vanishes since the K\"unneth component of $([x]+\mu\otimes
1)$ in $H^4(M)\otimes H^0(J_t)$ is zero.    If $Q'=Q+Y$ is another
choice of ``Poincar\'e bundle'' then, by~\eqref{eq-mu-again}
\[
q(Q')= q(Q) + \int_{M\times J/J}Y\cup(Q+\mu\otimes 1).
\]
Once again, one can check that the integral vanishes by looking
at K\"unneth components.

Witten constructs $\lineb$ (up to isomorphism) by giving a formula for
monodromy around loops in $J$.  We now show that the monodromy of
$q(x)$ can be computed by the same formula.

Let $\gamma:S^1\to J$ be a loop, and consider the following diagram:
\[
\begin{CD}
M\times S^1 @>\tilde\gamma>> M\times J \\
@VVV @VVV \\
S^1@>>\gamma > J \\
@VVV @. \\
\text{pt} @. .
\end{CD}
\]  
The monodromy of $q(x)$ around
$\gamma$ given by 
\[
\exp\left(2\pi i\int_{S^1}\gamma^\ast q(x)\right).
\]
By the base change property of $q$
(Corollary~\ref{thm:28}, \thmItemref{item:8}), the integral above can
be computed as
\[
\int_{S^1} q_{M\times S^1/S^1}(\tilde\gamma^{\ast} x).
\]
Choose any framing of the stable normal bundle of $S^1$.  By
transitivity (Corollary~\ref{thm:28}, \thmItemref{item:9}), this,
in turn, is given by
\[
q_{M\times S^1/\text{pt}}(\tilde\gamma^{\ast}(x)).
\]
We can compute this by finding a section, which we do by writing
$M\times S^1$ as a boundary, extending $\tilde\gamma^{\ast}(x))$, and
using~\eqref{eq-q-boundary}.

The $7$-dimensional $\spin$-manifold $M\times S^1$, together with the
cohomology class represented by $\tilde\gamma^\ast c$ defines an
element of
\[
M\spin_7 K(\Z,4).
\]
Since this group is zero, there is a $\spin$-manifold $N$ with
$\partial N=M\times S^1$ and a cocycle $\tilde c$ on $N$ extending
$\tilde\gamma^\ast c$.  We can then find
\[
y\in \chcat^4(N)
\]
whose restriction to $M\times S^1$ is $\tilde\gamma^\ast(x)$.
By~\eqref{eq-q-boundary}
\[
q_{M\times S^1/\text{pt}}(\tilde\gamma^{\ast}(x))
= -d\,q_N(y).
\]
This means that $q_N(y)$ is a real number whose reduction modulo $\Z$
is $\frac1{2\pi i}$ times the log of the monodromy.  Note, by the same
reasoning, that $q_N(0)$ must be an integer since it gives the
monodromy of the constant line bundle $q^{\lambda}(0)$.  So the
monodromy of $q^{\lambda}(x)$ (divided by $2\pi i$) is also given by
\[
q_N(y)-q_{N}(0) \mod \Z.
\]
Using Corollary~\ref{cor-2q-formula} we find
\[
q_N(y)-q_{N}(0)=\frac{-1}2 \int_N G\wedge G - G\wedge F,
\]
where $G$ is the $4$-form $\omega(y)$, and $F$ is $\frac12$ the first
Pontryagin form of $N$---that is $\omega(\lambda(N))$.  This is
(up to sign) Witten's formula.

Several comments are in order.
\begin{textList}
\item Our construction of the line bundle $\lineb$ works for any
$\spin$ manifold $M$ of dimension $(4k-2)$.  The computation of
monodromy we have given only applies when the reduced bordism group
\[
\mspin_{4k-1} K(\Z,2k)
\]
vanishes.  This group can be identified with 
\[
H_{2k-1}(\bspin;\Z)
\]
and does not, in general, vanish.
\item In the presence of torsion, as the manifold $M$ moves through
bordisms, the different components of the shifted Jacobians can come
together, suggesting that the ``correct'' shifted Jacobian to use is
not connected.  We have constructed a line bundle over the entire
group $\chh{4}^{4}(M)$, but it is independent of the choices only on
the component indexed by $-\mu$.  This suggests that it is better
to work with a line bundle over the {\em category $\chcat^{4}(M)$.}
Of course these remarks also apply in dimension $4k-2$.
\item In our formulation the quadratic functions arises from a
symmetry of the bundle.  Its value on the points of order two in the
connected component of the identity of $J$ can be computed from the
monodromy of the connection, but not in general.  This suggests that
it is important to remember automorphisms of objects in $\chcat^4(M)$,
and once again places priority on $\chcat^4(M)$ over its set of
isomorphism classes.
\item Formula~\eqref{eq-change-of-spin} gives the effect of a change
of $\spin$--structure on $\lineb$, and the resulting quadratic
function.
\end{textList}

\section{Cheeger--Simons cohomology}\label{sec-db-cs}

\subsection{Introduction}
\label{subsec-db-cs-intro}
We wish to define a cohomology theory of smooth manifolds $M$, which
encodes the notion of a closed $q$-forms with integral periods.  More
specifically, we are looking for a theory to put in the corner of the
square
\begin{equation}\label{eq-chs-pullback}
\begin{CD}
? @>>> \Omega^q_{\text{closed}}(M) \\
@VVV @VVV \\
H^q(M;\Z) @>>> H^q(M;\R),
\end{CD}
\end{equation}
which will fit into a long exact Mayer-Vietoris type of sequence.
Here $\Omega^q_{\text{closed}}(M)$ is the space of closed $q$-forms on
$M$.  There are two approaches to doing this, and they lead to
equivalent results.  The first is via the complex of smooth cochains
and is due to Cheeger and Simons.  The second is via complexes of
sheaves, and is due to Deligne and Beilinson.  In this paper we follow
the approach of Cheeger--Simons, and ultimately generalize it to the
theory of differential function spaces
(\S\ref{sec:gener-diff-cohom}).%
\footnote{Since writing this paper we have learned of the work of Harvey and
Lawson (see for example~\cite{harvey:_from_spark_grund})
which contains another treatment of Cheeger--Simons cohomology.}

\subsection{Differential Characters}
\label{subsec-differential-character}

To construct the cohomology theory sketched in the previous
section, we refine~\eqref{eq-chs-pullback}
to a diagram of cochain complexes, and define a complex
$\chcochain q^{\ast}(M)$ by the homotopy cartesian square
\begin{equation}\label{eq:44}
\begin{CD}
\chcochain q^{\ast}(M)  @>>> \Omega^{\ast\ge q}(M) \\
@VVV @VVV \\
C^{\ast}(M;\Z) @>>> C^{\ast}(M;\R).
\end{CD}
\end{equation}
More explicitly, the complex $\chcochain q^{\ast}(M)$ is given by
\[
\chcochain q^{n}(M)  = \begin{cases}
		    C^n(M;\Z)\times C^{n-1}(M;\R)\times \Omega^n(M)
		    \quad &n\ge q \\ 
		    C^n(M;\Z)\times C^{n-1}(M;\R)
		    \quad &n< q,
		    \end{cases}
\]
with differential
\begin{align*}
d(c,h,\omega) &= (\delta c, \omega-c-\delta h, d\omega) \\
d(c,h) &= \begin{cases}
	  (\delta c, -c-\delta h,0)&\qquad (c,h)\in\chcochain q^{q-1} \\
          (\delta c, -c-\delta h)  &\qquad \text{otherwise}.
	  \end{cases}
\end{align*}
It will be
convenient to write
\[
\chcochain 0^{*} =\chzero^{\ast}
\]
and regard 
\[
\chcochain q^{\ast} \subset\chzero^{\ast}
\]
as the subcomplex consisting of triples $(c,h,\omega)$ for which
$\omega=0$ if $\deg\omega<q$.

The Mayer-Vietoris sequence associated to~\eqref{eq:44}
\begin{multline*}
\to\dots \chh q^n(M) \to H^n(M;\Z)\times
H^n(\Omega^{\ast\ge q}) \\
\to H^n(M;\R) \to 
\chh q^{2k}(M) \to\dots
\end{multline*}
leads to natural isomorphisms
\[
\chh q^n(M) = \begin{cases}
			    H^n(M;\Z) &\quad n > q \\
			    H^{n-1}(M;\rz) &\quad n < q,
			    \end{cases}
\]
and a short exact sequence
\[
0\mapsto H^{q-1}(M)\otimes\rz \to \chh q^q(M)
\to A^q(M) \to 0.
\]
Here $A^q(M)$ is defined by the pullback square
\[
\begin{CD}
A^k(M) @>>> \Omega^k_{\text{cl}} \\
@VVV @VVV \\
H^k(M;\Z) @>>> H^k(M;\R),
\end{CD}
\]
and is thus the subgroup of $H^q(M;\Z)\times \Omega^q_{\text{cl}}$
consisting of pairs for which the (closed) form is a representative of the
image of the cohomology class in deRham cohomology.  This sequence can
be arranged in three ways
\begin{equation}\label{eq-sequences}
\begin{aligned}
0\to H^{q-1}(M;\R/\Z) \to&\chh q^{q}(M)\to \Omega^q_0 \to 0 \\
0\to H^{q-1}\otimes \R/\Z \to & \chh q^{q}(M)\to A^q(M)\to 0\\
0\to \Omega^{q-1}/\Omega^{q-1}_0 \to&\chh q^{q}(M)\to H^q(M;\Z)\to 0,
\end{aligned}
\end{equation}
In this, $\Omega^j_0$ denotes the space of closed $j$-forms
with integral periods.

A cocycle for the group $\chh q^q$ is a triple
\[
(c^q, h^{q-1},\omega),
\]
for which
\[
\delta h^{q-1} = \omega - c^q.
\]
The equivalence class of 
\[
(c^q, h^{q-1},\omega)
\]
in $\chh q^q(M)$
determines, and is determined by the pair
\[
(\chi,\omega),
\]
where
\[
\chi: Z_{q-1}\to\rz
\]
is given by
\[
\chi z = h^{q-1}(z) \mod \Z.
\]

\begin{defin}
Let $M$ be a smooth manifold.  A {\em differential character} of $M$
of degree ($k-1$) is a pair $(\chi, \omega)$ consisting of a character
\[
\chi:\Z_{k-1}\to\rz
\]
of the group of smooth $k$-cycles, and a $k$-form $\omega$
with the property that for every smooth $k$-chain $B$, 
\begin{equation}\label{eq-chs-condition}
\chi(\partial B) = \int_{B}\omega.
\end{equation}
\end{defin}

The group of differential characters was introduced by
Cheeger-Simons~\cite{cheeger85:_differ} and is denoted
\[
\chs^{k-1}(M).
\]
We will refer to $\chs^{k-1}(M)$ as the {\em $(k-1)^{\text{st}}$
Cheeger-Simons cohomology group.}  As indicated, the map
\[
(c^q, h^{q-1},\omega)\mapsto (\chi,\omega)
\]
gives an isomorphism
\[
\chh q^q(M)\approx \chs^{q-1}(M).
\]

The cup product in cohomology and the wedge product of forms 
lead to pairings 
\begin{gather}
\chcochain{k}^{\ast}(M)\otimes\chcochain{l}^{\ast}(M)\to\chcochain{k+l}^{\ast}\\
\chh{k}^{\ast}(M)\otimes\chh{l}^{\ast}(M)\to\chh{k+l}^{\ast}
\end{gather}
making $\chh{\ast}^{\ast}$ into a graded commutative ring.  As Cheeger
and Simons point out, the formula for these pairings is complicated
by the fact that the map from forms to cochains does not map the wedge
product to the cup product.  For
\[
\omega\in\Omega^r,\qquad \eta\in\Omega^s,
\]
let 
\[
\omega\cup\eta \in C^{r+s}
\]
the cup product of the cochains represented by $\omega$ and
$\eta$, (using, for example the Alexander--Whitney chain approximation to
the diagonal).  Let 
\begin{equation}\label{eq:24}
B(\omega,\eta)\in C^{r+s-1}
\end{equation}
be any natural chain homotopy between $\wedge$ and $\cup$:
\[
\delta B(\omega,\eta) + B(d\, \omega,\eta) + (-1)^{|\omega|}
B(\omega,d\,\eta)
= \omega\wedge\eta - \omega\cup\eta
\]
(any two choices of $B$ are naturally chain homotopic).
The product of cochains
$(c_1,h_1,\omega_1)$ and $(c_2,h_2,\omega_2)$
is given by the formula
\[
(c_1\cup c_2, (-1)^{|c_1|}c_1\cup h_2 + h_1\cup \omega_2 +
B(\omega_1,\omega_2), \omega_1\wedge\omega_2).
\]

As described in section~\ref{subsec-differential-cocycles} the group
\[
\chh 2^2(M)
\]
can be identified with the group of isomorphism classes of $U(1)$
bundles with connection on $M$, and the complex
$\chcochain{2}^\ast(M)$ can be used to give a more refined description
of the whole category of $U(1)$-bundles with connection on $M$.  We
refer the reader back to \S\ref{subsec-differential-cocycles} for an
interpretation of the groups $\chh k^k(M)$ and of the complex
$\chcochain{k}^\ast(M)$ for $k<2$.

\subsection{Characteristic classes} \label{sec:char-class} Let 
$G$ be a compact Lie group with Lie algebra $\mathfrak g$, and
let $W(\mathfrak g)$ be the Weyl algebra of polynomial functions on
$\mathfrak g$ invariant under the adjoint action.  We will
regard $W(\mathfrak g)$ as a differential graded algebra, with
differential $0$, and graded in such a way that $W(\mathfrak g)^{2n}$
consists of the polynomials of degree $n$, and $W(\mathfrak
g)^{\text{odd}}=0$.  Using~\cite{narasimhan61:_exist,narasimhan63:_exist}
choose a system of smooth $n$-classifying spaces $B^{(n)}(G)$ and
compatible connections on the universal bundles.  Write
$BG=\varinjlim B^{(n)}G$, 
\begin{align*}
C^{\ast}(BG;\Z) &=\varprojlim C^{\ast}(B^{(n)}G;\Z) \\
C^{\ast}(BG;\R) &=\varprojlim C^{\ast}(B^{(n)}G;\R),
\end{align*}
$\nabla_{\text{univ}}$ for the universal connection and
$\Omega_{\text{univ}}$ for its curvature.  The Chern--Weil homomorphism
(with respect to the universal connection) is a co-chain homotopy
equivalences
\[
W(\mathfrak g)\xrightarrow{\approx}{} C^{\ast}(BG;\R).
\]
Define $\check C^{\ast}(BG)$ by the homotopy Cartesian square
\[
\begin{CD}
\check C^{\ast}(BG) @>>> W(\mathfrak g) \\
@VVV @VV\text{Chern--Weil}V \\ 
C^{\ast}(BG;\Z) @>>> C^{\ast}(BG;\R).
\end{CD}
\]
The map $\check C^{\ast}(BG) \to C^{\ast}(BG)$ is then a cochain
homotopy equivalence, and we can take $\check C^{\ast}(BG)$ as a
complex of integer cochains on $BG$.  Specifically, an element $\check
C^{n}(BG)$ is a triple
\[
(c,h,w) \in C^{n}(BG;\Z)\times C^{n-1}(BG;\R)\times W(\mathfrak g)^{n}
\]
and the differential is given by
\[
d(c,h,w) = (\delta c,  w-c-\delta h,d w)
= (\delta c,  w-c-\delta h, 0).
\]

Suppose that $M$ is a manifold, and $f:M\to BG$ is a map classifying a
principal $G$-bundle $P\to M$.  The map $f$ is determined up to
homotopy by $P$.  A {\em characteristic class for principal
$G$-bundles} is a cohomology class $\gamma\in H^{k}(BG)$.  Because
of the homotopy invariance of cohomology, the class $f^{\ast}\gamma\in
H^{k}M$ depends only on $P\to M$ and not on the map $f$ which
classifies $P$.

Now suppose we have, in addition to the above data, a connection
$\nabla$ on $P$ (with curvature $\Omega$), and that the map $f$ is
smooth.  Choose a cocycle $(c,h,w)\in \check C^{k}(BG)$
representing $\gamma$.  
As in~\cite{chern74:_charac} consider the connection
\[
\nabla_{t}=(1-t)f^{\ast}\nabla_{\text{univ}} + t\nabla
\]
on $P\times I$, with curvature $\Omega_{t}$.  Then the form
\[
\eta=\int_{0}^{1}w(\Omega_{t})
\]
satisfies
\begin{align*}
d \eta & = w(\Omega_{1})-w(\Omega_{0}) \\
      & = w(\Omega)-f^{\ast}w(\Omega_{\text{univ}}).
\end{align*}
We can associate to this data the {\em characteristic differential cocycle}
\begin{equation}\label{eq:62}
(f^{\ast}c, f^{\ast}h +\eta, w(\Omega))\in \chcocycle{k}^{k}(M).
\end{equation}
The characteristic differential cocycle depends on a connection on
$P$, a (universal) choice of cocycle representing $\gamma$, and a
smooth map $f$ classifying $P$ (but not the connection).  Since the
differential in $W(\mathfrak g)$ is zero, varying $(c,h,w)$ by a
coboundary changes~\eqref{eq:62} by the coboundary of a class in
$\chcochain{k}^{k-1}(M)$.  Varying $f$ by a smooth homotopy also
results in a change of~\eqref{eq:62} by the coboundary of a class in
$\chcochain{k}^{k-1}(M)$.  It follows that the underlying cohomology
class $\check\gamma$ of~\eqref{eq:62} depends only on the principal
bundle $P$ and the connection $\nabla$.  In this way we recover the
result of Cheeger-Simons which, in the presence of a connection,
refines an integer characteristic class to a differential cocycle.

To summarize, given a principal $G$-bundle $P\to M$ and a cohomology
class $\gamma\in H^{k}(BG;Z)$ one has a characteristic class 
\[
\gamma(P)\in H^{k}(M;\Z).
\]
A choice of connection $\nabla$ on $P$ gives a refinement of $\gamma(P)$
to a cohomology class 
\[
\check\gamma(P,\nabla)\in \chhs{k}(M;\Z).
\]
A choice of cocycle (in $\check C^{\ast}(BG)$) representing $\gamma$,
and a map $M\to BG$ classifying $P$ gives a cocycle representative of 
$\check\gamma(P,\nabla)$.

Now suppose that $V$ is an oriented orthogonal representation of $G$
of dimension $n$, and $P\to M$ is a principal $G$-bundle.  Associated
to the orientation is a Thom class $U\in H^{n}(\bar V;\Z)$.  By the
results of Mathai-Quillen~\cite[\S6]{mathai86:_super_thom}, and the
methods above, a choice of connection $\nabla$ on $P$ gives
a differential Thom cohomology class $[\check U]\in \chhs{n}(\bar V)$,
and the additional choice of a map $M\to BG$ classifying $P$ gives a
differential Thom cocycle $\check U\in \chcocycle{n}^{n}(\bar V)$.

\subsection{Integration}\label{subsec-integration} 

We begin with the {\em integration map}
\[
\int_{S\times\rnbar N/S}:\chcochain{p+N}_{c}^{q+N}(S\times \R^N)\to
\chcochain{p}^{q}(S),
\]
Choose a fundamental cycle
\[
Z_N\in C_N(\R^N;\Z)
\]
(for example by choosing $Z_1$ and then taking $Z_N$ to be
the $N$-fold product), and map
\begin{equation}\label{eq-integration-rn}
(c,h,\omega)\mapsto
\left(c\slant Z_N,h\slant Z_N, \int_{S\times\rnbar N/S}\omega\right),
\end{equation}
where $a\slant b$ denotes the slant product, as described
in~\cite{spanier19:_algeb}.  One checks that the cochain $\omega\slant
Z_N$ coincides with the form $\int_{S\times\rnbar N/S}\omega$
(regarded as a cochain), so the above expression could be written
\[
(c\slant Z_N,h\slant Z_N,\omega\slant Z_N).
\]
It will be convenient to write this expression as
\[
x\slant Z_N
\]
with $x=(c,h,\omega)$.  It follows immediately
that~\eqref{eq-integration-rn} is a map of complexes, since ``slant
product'' with a closed chain is a map of complexes.

\begin{defin}
Suppose that $p:E\to S$ is an \oriented{} map of smooth manifolds with
boundary of relative dimension $k$.  The ``integration'' map
\[
\int_{E/S}:\chcochain{p+k}^{q+k}(E)\to
\chcochain{p}^{q}(S)
\]
is defined to be the composite
\[
\chcochain{p+k}^{q+k}(E)\xrightarrow{\cup U}{}
\chcochain{p+N}_c^{q+N}(S\times \R^N)\xrightarrow{\int_{\R^N}}{}
\chcochain{p}_c^{q}(S).
\]
\end{defin}

\begin{rem}
\label{rem:22}
In the terminology of \S\ref{subsec-differential-cocycles},
when $E/S$ is a fibration over an open dense subspace of $S$,
the map $\int_{E/S}$ commutes with the formation of the ``connection
form''; ie
\[
\text{connection form}\left(\int_{E/S}(c,h,\omega)\right)
= \int_{E/S}\omega
\]
in which the right hand integral indicates ordinary integration over
the fibers.  In particular, the connection form of
$\int_{E/S}(c,h,\omega)$ depends only on the orientation of the
relative normal bundle, and not the other choices that go into the
\orientation{} of $E/S$.  As we will see at the end of this section,
up to natural isomorphism the {\em integration functor} depends only
on the orientation of the relative normal bundle.
\end{rem}

Suppose that
\[
\begin{CD}
E_1 @>f_E>> E_2 \\
@VVV @VVV \\
S_1 @>>f_S> S_2
\end{CD}
\]
is a transverse pullback square in which $f_{S}$ is a closed embedding.  An
\orientation{} of $E_2/S_2$ induces an \orientation{} of $E_1/S_1$, and
the resulting integration functors are compatible with base change in
the sense that
\[
f_S^\ast\left(\int_{E_2/S_2}x\right) = \int_{E_1/S_1}f_E^\ast x.
\]

We now turn to our version of Stokes Theorem as described
in~\S\ref{subsec-differential-cocycles}.  Let $p:E\to S$ be an
orientable map in which $E$ is a manifold with boundary, and $S$ is
closed.  Choose a defining function
\[
f:E\to [0,1]
\]
for the boundary of $E$.  Then
\[
f\times p:E\to [0,1]\times S
\]
is a neat map of manifolds with boundary, and
\[
\begin{CD}
\partial E @>>> E \\
@VVV @VVV \\
\partial[0,1]\times S @>> \iota > [0,1]\times S
\end{CD}
\]  
is a transverse pullback square.  Choose an \orientation{} of $f\times
p$, and let $Z_I$ be the fundamental chain of $I=[0,1]$.  The
expression
\[
\int_0^1\int_{E/[0,1]\times S}x :=
\left.\left(\int_{E/[0,1]\times S}x\right)\right\slant Z_I
\]
satisfies Stokes theorem:
\begin{equation}\label{eq-stokes}
\delta\int_0^1\int_{E/[0,1]\times S}x 
=\int_0^1\int_{E/[0,1]\times S}\delta x
-(-1)^{|x|}
\left(\int_{E_1/S}x-
\int_{E_0/S}x\right)
\end{equation}
where $E_i=f^{-1}(i)$.  This follows easily from the formula
\[
\delta(a\slant b) = (\delta a)\slant b+(-1)^{|a|+|b|} a\slant\partial b
\]
and naturality.  For a discussion of this sign, and the slant product
in
general see \S\ref{subsec-slant-products}

When $x$ is closed, \eqref{eq-stokes} can be re-written as
\[
\delta\int_0^1\int_{E/[0,1]\times S}x=
-\int_{\partial E/S}x.
\]
Put more prosaically, this says that up to sign
\[
\int_0^1\int_{E/[0,1]\times S}x
\]
is a trivialization of
\[
\int_{\partial E/S}x.
\]

Returning to the situation of Remark~\ref{rem:22}, suppose 
given two \orientation{s} $\mu_{0}$ and $\mu_{1}$ of $E/S$
refining the same orientation of the relative normal bundle.  To
distinguish the two integration functors we will write them as 
\[
\int_{E/S}(\slot)\,d\mu_{0}\quad\text{and}\quad
\int_{E/S}(\slot)\,d\mu_{1}.
\]
Choose an \orientation{} $\mu$ of $E\times \Delta^{1}\to S\times
\Delta^{1}$ restricting to $\mu_{i}$ at $E\times \{i \}$.  For a
differential cocycle $\alpha=(c,h,\omega)\in \chcocycle{p}^{q}(E)$ set
\[
\beta=\int_{0}^{1}\int_{E\times \Delta^{1}/S\times
\Delta^{1}}p_{2}^{\ast}\alpha \,d\mu,
\]
with $p_{2}:E\times \Delta^{1}\to E$ the projection.  By the above
\[
d\beta=\int_{E/S}\alpha\,d\mu_{1}-\int_{E/S}\alpha\,d\mu_{0}.
\]
Now with no assumption, $\beta$ is a cochain in
$\chcochain{p-n-1}^{q-n-1}$.  But since $E/S$ is a fibration over an
open dense set, and integration commutes with the formation of the
curvature form, the curvature form of $\beta$ is zero.  It follows
that 
\[
\beta\in\chcochain{p-n-1}^{q-n}.
\]
This construction can then be regarded as giving a natural
isomorphism between
\[
\int_{E/S}(\slot)\,d\mu_{0}\quad\text{and}\quad
\int_{E/S}(\slot)\,d\mu_{1}.
\]

\subsection{Slant products}
\label{subsec-slant-products}
In this section we summarize what is needed to extend the definition
of the slant product from singular cochains to differential cochains.
The main thing to check is that the slant product of a form along a
smooth chain is again a form (Lemma~\ref{lem-slant-form} below).

Suppose that $M$ and $N$ are smooth manifolds, and $R$ a ring.  The
{\em slant product} is the map of complexes
\begin{equation}\label{eq:81}
C^{p+q}(M\times N;R)\otimes C_{q}(N)\to C^p(M;R)
\end{equation}
adjoint to the contraction
\begin{multline*}
C^{p+q}(M\times N;R)\otimes C_{p}(M)\otimes C_{q}(N) \\
\to
C^{p+q}(M\times N;R)\otimes C_{p+q}(M\times N)\to R.
\end{multline*}
There is a sign here, which is made troublesome by the fact that the
usual convention for the differential in the cochain complex does not
make it the dual of the chain complex (i.e.\ the evaluation map is not a
map of complexes).  To clarify, the map~\eqref{eq:81} is a map of complexes,
provided the differential in $C_{\ast}(M)$ is modified to be
\[
b\mapsto (-1)^{|b|}\partial b.
\]
Thus the relationship between the (co-)boundary and the slant product
is given by
\[
\delta(a\slant b) = (\delta a)\slant b+(-1)^{|a|+|b|} a\slant\partial b.
\]

To extend the slant product to the complex of differential cochains
\[
\chzero{\ast}(M\times N)\otimes C_{\ast}(N)\to \chzero^{\ast}(M)
\]
amounts to producing maps
\begin{align*}
\Omega^{\ast}(M\times N)\otimes C_{\ast}(N) &\xrightarrow{\slant}
\Omega^{\ast}(M) \\
C^{\ast}(M\times N;\R)\otimes C_{\ast}(N) &\xrightarrow{\slant}
C^{\ast}(M;\R) \\
C^{\ast}(M\times N;\Z)\otimes C_{\ast}(N) &\xrightarrow{\slant}
C^{\ast}(M;\R) \\
\end{align*}
compatible with the inclusions
\[
C^{\ast}(\slot;\Z)\hookrightarrow
C^{\ast}(\slot;\R)\hookleftarrow
\Omega^{\ast}(\slot;\R).
\]
This reduces to checking that the slant product of a form along a smooth
chain is again a form.
\begin{lem}\label{lem-slant-form}
Suppose that $\omega$ is a $(p+q)$-form on $M\times N$,
regarded as a cochain, and that $Z_p$ and $Z_q$ are $p$ and $q$-chains
on $M$ and $N$ respectively.  Then the value of $\omega\slant Z_q$
on $Z_p$ is
\[
\int_{Z_p\times Z_q}\omega.
\]
In other words, the cochain $\omega\slant Z_q$ is represented by the
form
\[
\int_{M\times Z_q/M}\omega.
\]
\end{lem}

\begin{pf}
By naturality we are reduced to the case in which
\begin{align*}
M&=\Delta^p \\
N&=\Delta^q
\end{align*}
and $Z_p$ and $Z_q$ are the identity maps respectively.  The value of
$\omega\slant Z_q$ on $Z_p$ is
\[
\int_{Z_p\otimes Z_q}\omega.
\]
But $Z_p\otimes Z_q$ is by definition the sum of all of the
non-degenerate $(p+q)$-simplices of $\Delta^p\times\Delta^q$, with
orientation derived from that of $\Delta^p\times\Delta^q$.  In other
words, $Z_p\otimes Z_q$ is the fundamental chain of of
$\Delta^p\times\Delta^q$.  It follows that
\[
\int_{Z_p\otimes Z_q}\omega=
\int_{Z_p\times Z_q}\omega.
\]
\end{pf}

\section{Generalized differential cohomology}
\label{sec:gener-diff-cohom}

\subsection{Differential function spaces} \label{sec:diff-funct-spac}
In \S\ref{sec-db-cs} we introduced the cohomology groups $\chh{q}^{n}(S)$
by combining differential forms and ordinary cohomology.  These groups
are formed from triples $(c,h,\omega)$ with $c$ a cocycle, $h$ a
cochain and $\omega$ a form.  For practical purposes an $n$-cocycle
can be regarded%
\footnote{It is well-known~\cite{Spanier} that for a CW complex $S$,
the set $[S,K(\Z,n)]$ can be identified with the cohomology group
$H^{n}(S;\Z)$.  In particular, any cocycle is cohomologous to one
which is pulled back from a map to an Eilenberg-MacLane space.  A more
refined statement is that the space of maps from $S$ to $K(\Z,n)$ has
the homotopy type of the ``space of $n$-cocycles'' on $S$.  For a more
precise discussion see Appendix~\ref{sec:simplicial-methods}}
as a map to the Eilenberg-MacLane space $K(\Z,n)$.  
We take this as our point of departure, and in this section shift the
emphasis from cocycles to maps.  Consider a topological space $X$
(with no particular smooth structure) a cocycle $\iota\in
Z^{n}(X;\R)$, and a smooth manifold $S$.

\begin{defin}
A {\em differential function} $t:S\to (X;\iota)$ is a triple
$(c,h,\omega)$
\[
c:S\to X, \quad h\in C^{n-1}(S;\R),\quad \omega\in \Omega^{n}(S)
\]
satisfying 
\[
\delta\,h = \omega-c^{\ast}\iota.
\]
\end{defin}

\begin{eg}\label{eg:5}
Suppose $X$ is the Eilenberg-MacLane space $K(\Z,2)$, which we can
take to be $\cp$, though we don't make use of its smooth
structure.  Choose a two-cocycle $\iota_{Z}$ representing the first
Chern class of the universal line bundle $L$, and let $\iota$ be its
image in $\Z^{2}(X;\R)$.  As described in Example~\ref{eg:6}, to
refine a map $c$ to a differential function amounts to putting a
connection on the line bundle $c^{\ast}L$.  The $U(1)$-bundle with
connection is the one associated to
\[
(c^{\ast}\iota_{Z},h,\omega)\in Z(2)^{2}(S).
\]
\end{eg} 

To form a {\em space} of differential functions we use the {\em
singular complex} of space $X^{S}$---a combinatorial object that
retains the homotopy type of the function space.  This requires the
language of simplicial sets.  We have outlined what is needed in
Appendix~\ref{sec:simplicial-methods}.  For further details
see~\cite{goerss99:_simpl,may92:_simpl,curtis71:_simpl}.  Rest assured
that we recover the complex for differential cohomology.  (See
Example~\ref{eg:8}, and Appendix~\ref{sec:comp-checkhznks-chsn}).

Let 
\[
\Delta^{n}=\{(t_{0},\dots,t_{n})\mid 0\le t_{i}\le 1, \sum t_{i}=1 \}
\]
be the standard $n$-simplex.  The {\em singular complex} of a space
$M$ is the simplicial set $\sing M$ with $n$-simplices
\[
{\sing M}_{n} = \left\{z:\Delta^{n}\to M\right\}.
\]
The evaluation map~\eqref{eq:75} 
\[
|\sing M|\to M
\]
induces an isomorphism of homotopy groups and of singular homology
groups.  The simplicial set $\sing M$ always satisfies the Kan
extension condition (\ref{def:6}), so the maps
\[
\spi_{n}\sing M\to \pi_{n}|\sing M|\to \pi_{n}(M)
\]
are all isomorphisms.

The singular complex of the function space $X^{S}$
is the simplicial set whose $k$-simplices are maps
\[
S\times \Delta^{k}\to X. 
\]

\begin{defin}\label{def:4}
Suppose $X$ is a space, $\iota\in
Z^{n}(X;\R)$ a cocycle, and $S$ a smooth manifold.
The {\em differential function complex} 
\[
\check X^{S}=(X;\iota)^{S}
\]
is the simplicial set whose $k$-simplices are differential functions
\[
S\times \Delta^{k}\to (X;\iota).
\]
\end{defin}

\begin{rem}\label{rem:16}
Phrased slightly differently, a $k$-simplex of $(X;\iota)^S$ is a
$k$-simplex 
\[
c:\Delta^{k}\times S\to X
\]
of $X^{S}$, together with a refinement of $c^{\ast}\iota$ to a
differential cocycle.  We will refer to the differential cocycle
refining $c^{\ast}\iota$ as the {\em underlying differential cocycle}.
\end{rem}

The differential function complex has an important filtration.  Define
a filtration of the deRham complex of $\Delta^{k}$ by
\[
\fil_{s}\Omega^{\ast}(\Delta^{k})= \left(\Omega^{0}(\Delta^{k})\to\dots
\to \Omega^{s}(\Delta^{k}) \right),
\]
and let
$\fil_{s}\Omega^{\ast}(S\times\Delta^{k})$ be the subspace
generated by
\[
\Omega^{\ast}(S)\otimes\fil_{s}\Omega^{\ast}(\Delta^{k}).
\]
In other words, a form $\omega\in\Omega^{\ast}(S\times\Delta^{k})$
lies in filtration $\le s$ if it vanishes on all sequences of vectors
containing $(s+1)$ vectors tangent to $\Delta^{k}$.  
Note that exterior differentiation shifts this filtration:
\[
d:\fil_{s}\to \fil_{s+1}.
\]

\begin{defin}\label{def:7}
With the notation of Definition~\ref{def:4}, a $k$-simplex $(c,h,\omega)$
of the differential function complex $(X;\iota)^S$
has {\em weight filtration $\le s$} if the form $\omega$ satisfies
\[
\omega\in\fil_{s}\Omega^{n}(S\times\Delta^{k})_{\text{cl}}. 
\]
\end{defin}

We will use the notation
\[
\fil_{s}(X;\iota)^S
\]
to denote the sub-simplicial set of $(X;\iota)^S$ consisting of
elements of weight filtration $\le s$.

\begin{eg}
Continuing with Example~\ref{eg:5}, a $1$-simplex of $(\cp,\iota)^{S}$
gives rise to a $U(1)$-bundle with connection $L$ over $S\times[0,1]$.  
which, in turn, leads to an isomorphism of $U(1)$-bundles
\begin{equation}\label{eq:22}
L\vert_{S\times \{0\}}\to
L\vert_{S\times \{1\}}
\end{equation}
by parallel transport along the paths
\[
t\mapsto (x,t).
\]
Note that this isomorphism need {\em not} preserve the connections.
\end{eg}

\begin{eg}
A $1$-simplex of $\fil_{0}(\cp;\iota)^{S}$ also gives rise to a
$U(1)$-bundle with connection $L$ over $S\times[0,1]$.  But this time
the curvature form must be pulled back%
\footnote{By definition it must be of the form $g(t)\omega$, where
$\omega$ is independent of $t$.  But it is also {\em closed}, so $g$
must be constant.}
from a form on $S$.  In this case the isomorphism~\eqref{eq:22} is
{\em horizontal} in the sense that it {\em does} preserve the
connections.
\end{eg}

\begin{rem}\label{rem:10}
We will show later (Lemma~\ref{thm:29} and Proposition~\ref{thm:31})
that this construction leads to a simplicial homotopy equivalence
between $\fil_{s}(\cp;\iota)^{S}$ and the simplicial abelian group
associated to the chain complex
\[
Z(2)^{2}(S)\leftarrow
C(2)^{1}(S)\leftarrow
C(0)^{0}(S).
\]
Thus the fundamental groupoid%
\footnote{A more precise description of the relationship between the
category of line bundles and the simplicial set $(\cp;\iota)^{S}$ can
be formulated in terms of the {\em fundamental groupoid.}  Recall that
the fundamental groupoid~\cite{reidemeister72:_einfueh_topol} of a
space $E$ is the groupoid $\pi_{\le 1}E$ whose objects are the points
of $E$, and in which a morphism from $x$ to $y$ is an equivalence of
paths starting at $x$ and ending at $y$.  The equivalence relation is
that of homotopy relative to the endpoints, and composition of maps is
formed by concatenation of paths.  The fundamental groupoid of a
simplicial set is defined analogously.} of $\fil_{0}(\cp;\iota)^{S}$ is
equivalent to the groupoid of $U(1)$-bundles with connection over $S$.
The fundamental groupoid of $\fil_{1}(\cp;\iota)^{S}$ is equivalent to
the groupoid of $U(1)$-bundles over $S$ and isomorphisms, and the
fundamental groupoid of $\fil_{2}(\cp;\iota)^{S}$ is equivalent to the
fundamental groupoid of $(\cp)^{S}$ which is equivalent to the
category of principal $U(1)$-bundles over $S$, and homotopy classes of
isomorphisms.
\end{rem}

\begin{eg}\label{eg:8}
Take $X$ to be the Eilenberg-MacLane space $K(\Z,n)$, and let
$\iota_{\R}\in Z^{n}(X;\R)$ be a fundamental cocycle: a cocycle whose
underlying cohomology class corresponds to the inclusion $\Z\subset\R$
under the isomorphism
\[
H^{n}(K(Z,n);\R)\approx \hom(\Z,\R).
\]
We show in Appendix~\ref{sec:comp-checkhznks-chsn} that the simplicial
set
\[
\fil_{s}(X;\iota)^S
\]
is homotopy equivalent to the simplicial
abelian group associated with the chain complex
(\S\ref{subsec-differential-character}) 
\[
\chcocycle{n-s}^{n}(S)
\xleftarrow{}{}
\chcochain{n-s}^{n-1}(S)\dots 
\xleftarrow{}{}
\chcochain{n-s}^{0}(S).
\]
The equivalence is given by slant product of the underlying
differential cocycle with the fundamental class of the variable
simplex.  It follows (Proposition~\ref{thm:5}) that the homotopy
groups of $\fil_{s}\check{X}^{S}$ are given by
\[
\pi_{i}\fil_{s}\check{X}^{S}=\chh{n-s}^{n-i}(S),
\]
and we recover the differential cohomology of $S$.  In this way the
homotopy groups of differential function spaces generalize
the Cheeger-Simons cohomology groups.
\end{eg}

\begin{rem}
We will also be interested in the situation in which we have several
cocycles $\iota$ of varying degrees.  These can be regarded as a
single cocycle with values in a graded vector space.  We will use the
convention 
\[
\grv_{j}=\grv^{-j},
\]
and so grade cochains and forms with values in a graded vector space
$\grv$ in such a way that the $C^{i}(X;\grv_{j})$ and
$\Omega^{i}(X;\grv_{j})$ have total degree $(i-j)$.  We will write
\begin{align*}
C^{\ast}(X;\grv)^{n} &=\bigoplus_{i+j=n} C^{i}(X; \grv^{j})\\
\Omega^{\ast}(S;\grv)^{n} &=\bigoplus_{i+j=n} \Omega^{i}(X; \grv^{j}) \\
Z^{\ast}(X;\grv)^{n} &=\bigoplus_{i+j=n} Z^{i}(X; \grv^{j}),
\end{align*}
and
\[
H^{n}(X;\grv)= \bigoplus_{i+j=n} H^{i}(X; \grv^{j}).
\]
We define $\fil_{s}\Omega^{\ast}(S\times\Delta^{k};\grv)$ to be the
subspace generated by
\[
\Omega^{\ast}(S;\grv)\otimes\fil_{s}\Omega^{\ast}(\Delta^{k}).
\]
As before, a form $\omega\in\Omega^{\ast}(S\times\Delta^{k};\grv)$ lies
in filtration $\le s$ if it vanishes on all sequences of vectors
containing $(s+1)$ vectors tangent to $\Delta^{k}$.
\end{rem}

We now turn to the analogue of the square~\eqref{eq:44} and the second
of the fundamental exact sequences~\eqref{eq-sequences}.  Using the
equivalence between simplicial abelian groups and chain complexes
(see~\S\ref{sec:simpl-abel-groups}), we can fit the differential
function complex $(X;\iota)^S$ into a homotopy Cartesian square
\begin{equation}\label{eq:23}
\begin{CD}
\fil_{s}(X;\iota)^S @>>> \fil_{s}
\Omega^{\ast}(S\times\Delta^{\bullet};\grv)^{n}_{\text{cl}} \\
@VVV @VVV \\
\sing X^{S} @>>> Z^{\ast}(S\times\Delta^{\bullet};\grv)^{n}.
\end{CD}
\end{equation}
By Corollary~\ref{thm:7}
\begin{align*}
\pi_{m}Z(S\times\Delta^{\bullet};\grv)^{n}&=H^{n-m}(S;\grv)\\
\pi_{m}\fil_{s}
\Omega^{\ast}(S\times\Delta^{\bullet};\grv)^{n}_{\text{cl}} &=
\begin{cases}
H^{n-m}_{\text{DR}}(S;\grv) &\qquad m<s\\
\Omega^{\ast}(S;\grv)^{n-s}_{\text{cl}} &\qquad m=s\\
0&\qquad m>s.
\end{cases}
\end{align*}
This gives isomorphisms
\begin{align*}
\pi_{k}\fil_{s}(X;\iota)^S &\xrightarrow{\approx}{}\pi_{k}X^{S} &\qquad
k < s \\
\pi_{k}\fil_{s}(X;\iota)^S & \xrightarrow{\approx}{}
H^{n-k-1}(S;\grv)/\pi_{k+1}X^{S} &\qquad k > s, 
\end{align*}
and a short exact sequence
\[
H^{n-s-1}(S;\grv)/\pi_{s+1}X^{S} \rightarrowtail
\pi_{s}\fil_{s}(X;\iota)^S \twoheadrightarrow
A^{n-s}(S;X,\iota),
\]
where $A^{n-s}(S;X,\iota)$ is defined by the pullback square
\[
\begin{CD}
A^{n-s}(S;X,\iota) @>>> \Omega(S;\grv)^{n-s}_{\text{cl}}  \\
@VVV @VVV \\
\pi_{s}X^{S} @>>> H^{n-s}(S;\grv).
\end{CD}
\]

\subsection{Naturality and homotopy}\label{sec:naturality-homotopy} We
now describe how the the differential function complex $(X;\iota)^{S}$
depends on $X$, $S$, and $\iota$.  First of all, a smooth map
\[
g:S\to T
\]
gives a map
\[
\fil_{s}(X;\iota)^{T} \to \fil_{s}(X;\iota)^{S}
\]
sending $(c,h,\omega)$ to $(c,h,\omega)\circ g = (c\circ g, c^{\ast}h,
c^{\ast}\omega)$.  We will refer to this map as {\em composition with
$g$}.

Given a map $f:X\to Y$ and a cocycle 
\[
\iota\in Z^{n}(Y;\grv),
\]
composition with $f$ gives a map
\[
\check f:\fil_{s}(X;f^{\ast}\iota)^{S}\to
\fil_{s}(Y;\iota)^{S}
\]
sending $(c,h,\omega)$ to $(c\circ f,h,\omega)$.  

\begin{prop}\label{thm:20}
Suppose that $f:X\to Y$ is a (weak) homotopy equivalence, and
$\iota\in Z^{n}(Y;\grv)$ is a cocycle.   Then for  each manifold $S$, the
map 
\[
\check f:\fil_{s}(X;f^{\ast}\iota)^{S}\to
\fil_{s}(Y;\iota)^{S}
\]
is a (weak) homotopy equivalence.
\end{prop}

\begin{pf}
When $f$ is a homotopy equivalence, the vertical maps in the following
diagram are homotopy equivalences (two of them are the identity map).  
\[
\begin{CD}
\sing X^{S} @>>> Z^{\ast}(S\times\Delta^{\bullet};\grv)^{n} @<<<
\fil_{s}
\Omega^{\ast}(S\times\Delta^{\bullet};\grv)^{n}_{\text{cl}} \\
@VVV @VVV @VVV \\
\sing Y^{S} @>>> Z^{\ast}(S\times\Delta^{\bullet};\grv)^{n} @<<< 
\fil_{s}\Omega^{\ast}(S\times\Delta^{\bullet};\grv)^{n}_{\text{cl}} \\
\end{CD}
\]
It follows that the map of homotopy pullbacks (see~\eqref{eq:23})
\begin{equation}\label{eq:29}
\check f:\fil_{s}(X;f^{\ast}\iota)^{S}\to
\fil_{s}(Y;\iota)^{S}
\end{equation}
is a homotopy equivalence.  If the map $f$ is merely a weak
equivalence, one needs to use the fact that a manifold with corners
has the homotopy type of a CW complex to conclude that the
left vertical map is a weak equivalence.  Since the formation of
homotopy pullbacks preserves weak equivalences the claim again
follows.
\end{pf}

\begin{rem}\label{rem:14}
Suppose that we are given a homotopy
\[
H:X\times\left[0,1 \right]\to Y,
\]
with $H(x,0)=f(x)$ and $H(x,1)=g(x)$.  We then have a diagram of
differential function spaces
\[
\xymatrix{
\fil_{s}\left(X,f^{\ast}\iota \right)^{S} \ar[r]\ar[dr]_{\check
f} &
\fil_{s}\left(X\times I ;
H^{\ast}\iota \right)^{S} \ar[d] &
\fil_{s}\left(X;g^{\ast}\iota \right)^{S} \ar[l]\ar[dl]^{\check
g} \\
& \fil_{s}\left(Y,\iota \right)^{S}. & 
}
\]
By Proposition~\ref{thm:20}, the horizontal maps are homotopy equivalences,
and so the construction can be regarded as giving a homotopy between
$\fil_{s}\check f$ and $\fil_{s}\check g$.
\end{rem}

Given two cocycles $\iota_{1},\iota_{2}\in Z^{\ast}(X;\grv)^{n}$, and a cochain
$b\in C^{\ast}(X;\grv)^{n-1}$ with $\delta b=\iota_{1}-\iota_{2}$, we get a map
\begin{equation}\label{eq:36}
\begin{aligned}
(X;{\iota_{1}})^{S} &\to (
X;{\iota_{2}})^{S} \\
(c,h,\omega) &\mapsto (c, h+c^{\ast}b,\omega).
\end{aligned}
\end{equation}
This map is an isomorphism, with inverse $-b$.  In particular, the
group $Z^{\ast}(X;\grv)^{n-1}$ acts on the differential function complex
$(X;\iota)^{S}$.

Finally, suppose given a map $t:\grv\to \grw$ of graded vector spaces,
and a cocycle $\iota\in Z^{n}(X;\grv)$.  Composition with $t$
defines a cocycle 
\[
t\circ\iota\in Z^{n}(X;W),
\]
and a map of differential function complexes
\[
\left(X;\iota \right)^{S}
\to\left(X;t\circ\iota \right)^{S}.
\]

Combining these, we see that given maps 
\[
f:X\to Y\qquad t:\grv\to \grw,
\]
cocycles 
\[
\iota_{X}\in Z^{n}(X;\grv)\qquad\iota_{Y}\in Z^{n}(Y;W),
\]
and a cochain $b\in C^{n-1}(X;W)$ with $\delta
b=t\circ \iota_{X} - f^{\ast}\iota_{Y}$ we get a map of differential
function complexes
\begin{align*}
\left(X;\iota_{X}\right)^{S} &\to \left(Y;\iota_{Y}\right)^{S} \\
(c,h,\omega) &\mapsto (f\circ c, t\circ h + c^{\ast}b,t\circ\omega)
\end{align*}
which is a weak equivalence when $f$ is a weak equivalence and $t$ is
an isomorphism.

\begin{rem}\label{rem:15}
All of this means that the homotopy type of 
\[
\fil_{s}(X,\iota)^{S}
\]
depends only on the cohomology class of $\iota$
and and the homotopy type of $X$.  For example, suppose $f$ and $g$
are homotopic maps $X\to Y$, and that $\alpha, \beta\in Z^{n}(X)$ are
cocycles in the cohomology class of $f^{\ast}\iota$ and
$g^{\ast}\iota$ respectively.  A choice of cochains 
\begin{align*}
b_{1}, b_{2} &\in C^{n-1}(X,\R) \\
\delta b_{1} &=\alpha-f^{\ast}\iota \\
\delta b_{2} &=\beta-g^{\ast}\iota \\
\end{align*}
gives isomorphisms
\begin{align*}
\fil_{s}(X;\alpha)^{S} &\approx
\fil_{s}(X;f^{\ast}\iota)^{S} \\
\fil_{s}(X;\beta)^{S} &\approx
\fil_{s}(X;g^{\ast}\iota)^{S}, 
\end{align*}
and defines maps
\begin{align*}
\fil_{s}\check f:\fil_{s}(X;\alpha)^{S} &\to
\fil_{s}(Y;\iota)^{S} \\
\fil_{s}\check g:\fil_{s}(X;\beta)^{S} &\to
\fil_{s}(Y;\iota)^{S}.
\end{align*}
A homotopy $H:X\times\left[0,1 \right]\to Y$ from $f$ to
$g$, then leads to the diagram of Remark~\ref{rem:14}, and hence a
homotopy between $\fil_{s}\check f$ and $\fil_{s}\check g$.
\end{rem}

\subsection{Thom complexes}\label{sec:thom-complexes} In this section
we describe Thom complexes and the Pontryagin-Thom construction in 
the context of differential function complexes.

Recall that the {\em Thom complex} of a vector bundle $V$ over a
compact space $M$ is the $1$-point compactification of the total
space.  We will use the notation $\thom(M;V)$, or $\bar V$ to denote
the Thom complex.  (We are avoiding the more traditional $M^{V}$
because of the conflict with the notation $X^{S}$ for function
spaces).  When $M$ is not compact, one sets
\[
\thom(M;V)= \bar V=
\bigcup_{\substack{K_{\alpha}\subset M \\ \text{compact}}}
\thom(K_{\alpha};V_{\alpha})\qquad V_{\alpha}=V\vert_{K_{\alpha}}.
\]

Suppose $V$ is a vector bundle over a manifold $M$.  We will call a
map $g:S\to \bar V$ {\em smooth} if
its restriction to
\[
g^{-1}(V) \to  V
\]
is smooth.  We define the deRham complex $\Omega^{\ast}(\bar V)$ to be
the sub-complex of $\Omega^{\ast}(V)$ consisting of forms which are
fiber-wise compactly supported.  With these definitions, the
differential function complex 
\[
\fil_{s}(X;\iota)^{\bar V}
\]
is defined, and a smooth map $S\to \bar V$ gives a map of differential
function complexes
\[
\fil_{s}(X;\iota)^{\bar V}
\to
\fil_{s}(X;\iota)^{S}
\]

\begin{defin}
Let
\[
W_{1}\to X\quad\text{and}\quad W_{2}\to Y
\]
be vector bundles, and $\iota\in
Z^{k}\left(\bar{W}_{2},\{\infty\};\grw\right)$ a cocycle.  A {\em
(vector) bundle map} $W_{1}\to W_{2}$ is a pullback square
\[
\begin{CD}
W_{1} @>t>> W_{2} \\
@VVV @VVV \\
X @>>f> Y
\end{CD}
\]
for which the induced isomorphism $W_{1}\to f^{\ast}W_{2}$ is an
isomorphism of vector bundles.  A {\em differential bundle map} is a
differential function
\[
\check t=(c,h,\omega): \bar W_{1} \to \left(\bar W_{2};\iota \right).
\]
which is a bundle map.  The {\em complex of differential bundle maps}
\[
\fil_{s}\left(W_{2};\iota \right)^{W_{1}}\subset 
\fil_{s}\left(\bar{W}_{2};\iota \right)^{\bar{W}_{1}}
\]
is the subcomplex consisting of vector bundle maps.  
\end{defin}

In case $W_{2}$ is the universal bundle over some kind of classifying
space, then a vector bundle map $W_{1}\to W_{2}$ is a classifying map
for $W_{1}$.  

The set of bundle maps is topologized as a closed subspace of the
space of all maps from $\thom(X;W_{1})\to \thom(Y;W_{2})$ (which, in
turn is a closed subspace of the space of all maps $W_{1}\to W_{2}$).

Suppose that $B$ is a topological space, equipped with a vector bundle
$V$, and $\iota\in Z^{k}\left(\bar V,\{\infty \}\right)$ is a cocycle.

\begin{defin}
A {\em $B$-oriented embedding} is a neat embedding $p:E\to S$, a
tubular neighborhood\footnote{Recall from footnote~\ref{fn:1} that a {\em tubular neighborhood} of
$p:E\hookrightarrow S$ is a vector bundle $W$ over $E$, and an
extension of $p$ to a diffeomorphism of $W$ with a neighborhood of
$p(E)$.}%
\ $W\hookrightarrow S$ of $p:E\to S$, and a vector
bundle map $W\to V$ classifying $W$.  A {\em differential $B$-oriented
embedding} is a neat embedding $p:E\to S$, a tubular neighborhood
$W\hookrightarrow S$ of $p:E\to S$, and a differential vector bundle
map
\[
\check{t}: W\to \left(V;\iota \right)
\]
classifying $W$.
\end{defin}

Let $p:E\hookrightarrow S$ be a differential oriented embedding.  The
construction of Pontryagin-Thom gives a smooth map
\begin{equation}\label{eq:83}
S\to \thom(E,W).
\end{equation}
Composition with~\eqref{eq:83} defines the {\em push-forward}
\[
p_{!}:\fil_{s}(V;\iota)^{W} \subset 
\fil_{s}\left(X;\iota \right)^{\thom(E,W)}
\to \fil_{s}(X;\iota)^{S}.
\]

\subsection{Interlude: differential
$K$-theory}\label{sec:interl-diff-k}

Before turning to the case of a general cohomology theory, we apply
the ideas of the previous section to the case of $K$-theory.  The
resulting {\em differential $K$-theory} originally came up in anomaly
cancellation problems for $D$-branes in
$M$-theory~\cite{freed:_ramon_ramon_k,freed00:_dirac_charg_quant_gener_differ_cohom}.
The actual anomaly cancellation requires a refinement of the families
index theorem, ongoing joint work of the authors and Dan
Freed. 

Let $\fred$ be the space of Fredholm operators.  We remind the reader
that the space $\fred$ is a classifying space for $K$-theory and in
particular that any vector bundle can be obtained as the index bundle
of a map into $\fred$.  Let
\[
\iota=(\iota_{n})\in \prod Z^{2n}(\fred;\R)=Z^{0}(\fred;\grv)
\]
be a choice of cocycle representatives for the universal Chern
character, so that if $f:S\to\fred$ classifies a vector bundle $V$,
then the characteristic class $\ch_{n}(V)\in H^{2n}(S;\R)$ is
represented by $f^{\ast}\iota_{n}$.

\begin{defin}
The differential $K$-group $\check K^{0}(S)$ is the group 
\[
\pi_{0}\fil_{0}\left(\fred;\iota \right)^{S}
\]
\end{defin}

\begin{rem}
In other words, an element of $\check K^{0}(S)$  is represented by a
triple $(c,h,\omega)$ where $c:S\to\fred$ is a map,
$\omega=(\omega_{n})$ is a sequence of $2n$-forms, and $h=(h_{n})$ is
a sequence of $(2n-1)$ cochains satisfying 
\[
\delta h = \omega-c^{\ast}\iota.
\]
Two triples $(c^{0},h^{0},\omega^{0})$ and $(c^{1},h^{1},\omega^{1})$
are equivalent if there is a $(c,h,\omega)$ on $S\times I$, with
$\omega$ constant in the $I$ direction, and with 
\begin{align*}
(c,h,\omega)\vert_{\{0 \}} &= (c^{0},h^{0},\omega^{0}) \\
(c,h,\omega)\vert_{\{1 \}} &= (c^{1},h^{1},\omega^{1}).
\end{align*}
\end{rem}

We will use the symbol $\check{ch}(\check V)$ to denote the
differential cocycle underlying a differential function $\check
V:S\to(\fred;\iota)$.

The space $\Omega^{i}\fred$ is a classifying space for $K^{-i}$.  Let
\[
\iota^{-i}=(\iota^{-i}_{2n-i})\in\prod Z^{2n-i}(\Omega^{i}\fred;\R)
\]
be the cocycle obtained by pulling $\iota$ back along the evaluation
map
\[
S^{i}\times\Omega^{i}\fred\to\fred
\]
and integrating along $S^{i}$ (taking the slant product with the
fundamental cycle of $S^{i}$).  The cohomology classes of the
$\iota^{-i}_{2n-i}$ are the universal even (odd) Chern character
classes, when $i$ is even (odd).

\begin{defin}
The differential $K$-group $\check K^{-i}(S)$ is the group 
\[
\pi_{0}\fil_{0}\left(\Omega^{i}\fred;\iota^{-i} \right)^{S}
\]
\end{defin}

\begin{rem}
As above, an element of $\check K^{-i}(S)$  is represented by a
triple $(c,h,\omega)$ where $c:S\to\Omega^{i}\fred$ is a map,
$\omega=(\omega_{n})$ is a sequence of $2n-i$-forms, and $h=(h_{n})$ is
a sequence of $(2n-i-1)$ cochains satisfying 
\[
\delta h = \omega-c^{\ast}\iota^{-i}.
\]
Two triples $(c^{0},h^{0},\omega^{0})$ and $(c^{1},h^{1},\omega^{1})$
are equivalent if there is a $(c,h,\omega)$ on $S\times I$, with
$\omega$ constant in the $I$ direction, and with 
\begin{align*}
(c,h,\omega)\vert_{\{0 \}} &= (c^{0},h^{0},\omega^{0}) \\
(c,h,\omega)\vert_{\{1 \}} &= (c^{1},h^{1},\omega^{1}).
\end{align*}
\end{rem}

These differential $K$-groups lie in short exact sequences
\begin{gather*}
0\to K^{-1}(S;\R/\Z) \to \check K^{0}(S)\to
\prod\Omega^{2n}(S) \to 0 \\ 
0\to K^{-1}(S)\otimes \R/\Z \to   \check K^{0}(S)\to A_{K}^0(S)\to 0
\end{gather*}
and 
\begin{gather*}
0\to K^{-i-1}(S;\R/\Z) \to \check K^{-i}(S)\to
\prod\Omega^{2n-i}(S) \to 0 \\ 
0\to K^{-i-1}(S)\otimes \R/\Z \to   \check K^{-i}(S)\to A_{K}^{-i}(S)\to 0,
\end{gather*}
where $A_{K}^{-i}$ is defined by the pullback square
\[
\begin{CD}
A_{K}^{-i}(S) @>>> \prod \Omega^{2n-i}_{\text{cl}}(S) \\
@VVV @VVV \\
K^{2n-i}(S) @>>\text{ch}> \prod H^{2n-i}(S;\R).
\end{CD}
\]
So an element of $A_{K}^{-i}(S)$ consists of an element
$x\in K^{2n-i}(S)$ and a sequence of closed $2n-i$-forms representing
the Chern character of $x$.

Bott periodicity provides a homotopy equivalence 
\[
\Omega^{n}\fred\approx \Omega^{n+2}\fred
\]
under which $\iota_{n}$ corresponds to $\iota_{n+2}$.  This gives an
equivalence of differential function spaces
\[
\fil_{t}\left(\Omega^{n}\fred;\iota_{n} \right)^{S}\approx
\fil_{t}\left(\Omega^{n+2}\fred;\iota_{n+2}\right)^{S} 
\]
and in particular, of differential $K$-groups
\[
\check{K}^{-n}(S)\approx
\check{K}^{-n-2}(S).
\]
This allows one to define differential $K$-groups $\check K^{n}$ for
$n>0$, by
\[
\check{K}^{n}(S)=
\check{K}^{n-2N}(S)\qquad n-2N <0.
\]

In~\cite{lott94:_r_z} Lott defines the group $K^{-1}(S;\rz)$ in
geometric terms, and proves an index theorem, generalizing the index
theorem for flat bundles in~\cite{aps3}.  Lott takes as
generators of $K^{-1}(S;\rz)$ pairs $(V,h)$ consisting of a (graded)
vector bundle $V$ with a connection, and sequence $h$ of odd forms
satisfying
\[
dh=\text{Chern character forms of $V$}.
\]
Our definition is close in spirit to Lott's, with $c$ corresponding to
$V$, $c^{\ast}\iota$ to the Chern character forms of the connection,
$h$ to $h$ and $\omega=0$.  In a recent preprint~\cite{lott:_higher},
Lott constructs the abelian gerbe with connection whose curvature is
the $3$-form part of the Chern character of the index of a family of
self-adjoint Dirac-like operators.  In our terminology he constructs
the degree $3$ part of the differential Chern character going from
differential $K$-theory to differential cohomology.

\subsection{Differential cohomology theories}  In this section we
build on
our theory of differential functions and define differential
cohomology theories.  We explained in the introduction why we wish to
do so.

Let $E$ be a cohomology theory.  The spaces representing $E$
cohomology groups ($\Omega^{n}\fred$ in the case of $K$-theory,
Eilenberg-MacLane spaces for ordinary cohomology) fit together to form
a {\em spectrum}.

\begin{defin}
(see~\cite{LMayS,HopGoerss:Mult,elmendorf95:_moder,elmendorf97:_rings,Ad:SHGH})
A {\em spectrum} $E$ consists of a sequence of pointed spaces $E_{n}$,
$n=0,1,2,\dots$ together with maps
\begin{equation}
s_{n}^{E}:\Sigma E_{n}\to E_{n+1}
\end{equation}
whose adjoints
\begin{equation}\label{eq:4}
t_{n}^{E}:E_{n}\to \Omega E_{n+1}
\end{equation}
are homeomorphisms.
\end{defin}

\begin{rem}
By setting, for $n>0$
\[
E_{-n}=\Omega^{n}E_{0}=\Omega^{n+k}E_{k},
\]
a spectrum determines a sequence of spaces $E_{n}$, $n\in\Z$ together
with homeomorphisms
\[
t_{n}^{E}:E_{n}\to
\Omega E_{n+1}
\]
\end{rem}

If $X$ is a pointed space and $E$ is a spectrum, then $E$-cohomology
groups of $X$ are given by
\[
E^{k}(X) = [X,E_{k}]=[\Sigma^{N}X, E_{N+k}]
\]
and the $E$-homology groups are given by
\[
E_{k}(X) = \lim_{N\to\infty} \pi_{N+k}E_{N}\wedge X.
\]

\begin{eg}\label{eg:10}
Let $A$ be an abelian group, and $K(A,n)$ an Eilenberg-MacLane space.
Then $[S,K(A,n)]=H^{n}(S;A)$.  To assemble these into a {\em spectrum}
we need to construct homeomorphisms
\begin{equation}\label{eq:53}
K(A,n)\to \Omega K(A,n+1)
\end{equation}
By standard algebraic topology methods, one can choose
the spaces $K(A,n)$ so that~\eqref{eq:53} is a closed
inclusion.  Replacing $K(A,n)$ with  
\[
\varinjlim \Omega^{N}K(A,N+n)
\]
then makes~\eqref{eq:53} a homeomorphism.  This is the {\em
Eilenberg-MacLane spectrum}, and is denoted $HA$.  By construction
\begin{align*}
\HA^{n}(S) &=H^{n}(S;A)\\
\HA_{n}(S) &=H_{n}(S;A).
\end{align*}
\end{eg}

\begin{eg}
Now take $E_{2n}=\fred$, and $E_{2n-1}=\Omega\fred$.  We have a
homeomorphism $E_{2n-1}=\Omega E_{2n}$ by definition, and a homotopy
equivalence $E_{2n}\to\Omega E_{2n+1}$ by Bott periodicity:
\[
E_{2n}\xrightarrow{}{}\Omega^{2}E_{4k}=\Omega E_{2n+1}.
\]
As in Example~\ref{eg:10}, the spaces $E_{n}$ can be modified so that
the maps $E_{n}\to \Omega E_{n+1}$ form a spectrum, the $K$-theory
spectrum.
\end{eg}

Multiplication with the fundamental cycle $Z_{S^{1}}$ of $S^{1}$ gives
a map of singular chain complexes
\[
C_{\ast}E_{n}\to
C_{\ast+1}\Sigma E_{n}\to
C_{\ast+1} E_{n+1}.
\]

\begin{defin}
Let $E$ be a spectrum.  The {\em singular chain complex} of $E$ is the
complex
\[
C_{\ast}(E) =\varinjlim C_{\ast+n}E_{n},
\]
and the {\em singular cochain complex} of $E$ (with coefficients in an
abelian group $A$) is the cochain complex
\[
C^{\ast}(E;A)=\hom(C_{\ast}E,A)=\varprojlim(C^{\ast+{n}}E_{n})
\subset \prod C^{\ast+{n}}E_{n}.
\]
\end{defin}

The {\em homology} and {\em cohomology} groups of $E$ are the homology
groups of the complexes $C_{\ast}(E)$ and $C^{\ast}E$.  Note that these
groups can be non-zero even when $k$ is negative.

Now fix a cocycle $\iota\in Z^{p}(E;\grv)$.  By definition, this means
that there are cocycles $\iota_{n}\in Z^{p+n}(E_{n};\grv)$, $n\in\Z$
which are compatible in the sense that
\[
\iota_{n}= \left(s_{n}^{\ast}\,\iota_{n+1} \right)\slant Z_{S^{1}}.
\]
Once in a while it will be convenient to denote the pair
$(E_{n};\iota_{n})$ as
\begin{equation}\label{eq:86}
(E_{n};\iota_{n})=(E;\iota)_{n}.
\end{equation}

\begin{defin}
Let $S$ be a manifold.  The {\em differential $E$-cohomology group}
\[
E(n-s)^{n}(S;\iota)
\]
is the homotopy group
\[
\pi_{0}\fil_{s}\left(E_{n};\iota_{n}\right)^{S}
=
\pi_{0}\fil_{s}\left(E;\iota\right)_{n}^{S}.
\]
\end{defin}

\begin{rem}
We will see in the next section that the
spaces 
\[
\left(E_{n};\iota_{n} \right)^{S}
\]
come equipped with homotopy equivalences
\begin{equation}\label{eq:40}
\fil_{s+n}(E_{n};\iota_{n})^{S} \to
\Omega\fil_{s+n+1}(
E_{n+1};\iota_{n+1})^{S}.
\end{equation}
Thus there exists a spectrum 
\[
\fil_{s}(E;\iota)^{S}
\]
with 
\[
\left(\fil_{s}(E;\iota)^{S} \right)_{n}=
\fil_{s+n}(E;\iota)^{S}_{n}=
\fil_{s+n}(E_{n};\iota_{n})^{S},
\]
and the higher homotopy group 
\[
\pi_{t}\fil_{s}\left(E_{n};\iota_{n}\right)^{S}
\]
is isomorphic to 
\[
\pi_{0}\fil_{s-t}\left(E_{n-t};\iota_{n-t}\right)^{S} = 
E(n-s)^{n-t}(S;\iota).
\]
\end{rem}

\begin{eg}\label{eg:14}
We will show in Appendix~\ref{sec:comp-checkhznks-chsn} that when $E$
is the Eilenberg-MacLane spectrum $\HZ$, one has
\[
\HZ(n)^{k}(S)=\chh{n}^{k}(S).
\]
\end{eg}

\subsection{Differential function spectra}
\label{sec:diff-funct-spectr}

In this section we will establish the homotopy
equivalence~\eqref{eq:40}.  As a consequence there exists a {\em
differential function spectrum} $\fil_{s}\left(E;\iota \right)^{S}$
whose $n^{\text{th}}$ space has the homotopy type of
$\fil_{s+n}\left(E_{n};\iota_{n} \right)^{S}$.

There is one point here on which we have been deliberately vague, and
which needs to be clarified.  Differential function complexes are not
spaces, but simplicial sets.  Of course they can be made into spaces
by forming their geometric realizations, and, since they are Kan
complexes, no homotopy theoretic information is gained or lost by
doing so.  But this means that the notation
$\Omega\fil_{s+n+1}\left(E_{n+1};\iota_{n+1} \right)^{S}$ is
misleading, and that the object we need to work with is the space of
{\em simplicial loops}.

For any simplicial set $X$ let $\Omegasimp X$ be the {\em simplicial
loop space of $X$}; the simplicial set whose $k$-simplices are the
maps of simplicial sets $h:\Delta^{k}_{\bullet}\times
\Delta^{1}_{\bullet}\to X_{\bullet}$ for which $h(x,0)=h(x,1)=\ast$.
Using the standard simplicial decomposition of
$\Delta^{k}_{\bullet}\times \Delta^{1_{\bullet}}$, a $k$-simplex of
$\Omegasimp X$ can be described as a sequence
\[
h_{0},\dotsc ,h_{k}\in X_{k+1}
\]
of $(k+1)$-simplices of $X$ satisfying
\begin{equation}\label{eq:34}
\begin{gathered}
\partial_{i}h_{i} = \partial_{i}h_{i-1} \\
\partial_{0}h_{0} = \partial_{k+1}h_{k} = \ast.
\end{gathered}
\end{equation}
There is a canonical map
\[
\left|\Omegasimp X \right|\to \Omega\left|X \right|
\]
which is a homotopy equivalence if $X$ satisfies the Kan extension
condition (\ref{def:6}).  The simplicial set $\fil_{s+n+1}\left(\check
E_{n+1}^{S};\iota_{n+1} \right)$ satisfies the Kan extension
condition, and in this section we will actually produce a simplicial
homotopy equivalence
\begin{equation}\label{eq:76}
\fil_{s+n}(E_{n};\iota_{n})^{S} \to
\Omegasimp\fil_{s+n+1}(E_{n+1};\iota_{n+1})^{S}.
\end{equation}

Let
\[
E_{n+1,c}^{S\times\R}
\]
be the space of ``compactly supported'' functions, i.e., the fiber of
the map
\[
E_{n+1}^{S\times \bar\R}\to E_{n+1}^{S\times \left\{\infty \right\}},
\]
in which $\bar\R=\R\cup\{\infty \}$ is the one point compactification
of $\R$.  Of course, $E_{n+1,c}^{S\times\R}$ is simply the loop space
of $E_{n+1}^{S}$, and is homeomorphic to $E_{n}^{S}$.  Define a
simplicial set $\fil_{s}(E_{n+1};\iota_{n+1})^{S\times\R}_{c}$ by the homotopy Cartesian square
\begin{equation}\label{eq:31}
\begin{CD}
\fil_{n+s+1}(E_{n+1};\iota_{n+1})^{S\times\R}_{c} @>>> 
\fil_{n+s+1}\Omega_{c}^{\ast}\left(S\times\R\times\Delta^{\bullet};\grv
\right)^{n+p+1}_{\text{cl}} \\
@VVV @VVV \\
\sing E_{n+1,c}^{S\times\R} @>>> Z_{c}^{\ast}\left(S\times\R\times\Delta^{\bullet};\grv
\right)^{n+p+1}
\end{CD}
\end{equation}
The vector space $\Omega_{c}^{\ast}\left(S\times\R\times\Delta^{k};\grv
\right)$ is the space of forms which are compactly supported
(along $\R$), and the subspace
\[
\fil_{t}\Omega_{c}^{\ast}\left(S\times\R\times\Delta^{k};\grv
\right)
\]
consists of those whose Kunneth components along
$\Delta^{k}\times\R$ of degree greater than $t$ vanish.

We will construct a diagram of simplicial homotopy equivalences
\begin{equation}\label{eq:78}
\begin{CD}
\fil_{n+s+1}(E_{n+1};\iota_{n+1})^{S\times\R}_{c}
@>>> \Omegasimp\fil_{n+s+1}(E_{n+1};\iota_{n+1})^{S} \\
@VVV @. \\
\fil_{s+n}(E_{n};\iota_{n})^{S} @. 
\end{CD}
\end{equation}
Choosing a functorial section of the leftmost map gives~\eqref{eq:76}.

For the left map of~\eqref{eq:78}, note that the homeomorphism
$E_{n+1,c}^{S\times\R}\approx E_{n}^{S}$ prolongs to an isomorphism of
simplicial sets $\sing E_{n+1,c}^{S\times\R} \approx \sing E_{n}^{S}$
which is compatible with the maps
\begin{gather*}
\fil_{n+s+1}\Omega_{c}^{\ast}\left(S\times\R\times\Delta^{\bullet};\grv
\right)^{n+p+1}_{\text{cl}} \to
\fil_{n+s}\Omega^{\ast}\left(S\times\Delta^{\bullet};\grv
\right)^{n+p}_{\text{cl}}  \\
Z_{c}^{\ast}\left(S\times\R\times\Delta^{\bullet};\grv
\right)^{n+p+1} \to
Z^{\ast}\left(S\times\Delta^{\bullet};\grv
\right)^{n+p}
\end{gather*}
given by ``integration over $\R$'' and ``slant product with the
fundamental class of $\bar\R$.''   By Corollary~\ref{thm:24} these maps
are simplicial homotopy equivalences.  Passing to homotopy pullbacks then
gives a homotopy equivalence
\begin{equation}\label{eq:32}
\fil_{n+s+1}(E_{n+1,c};\iota_{n+1})^{S\times\R} \xrightarrow{\sim}{}
\fil_{s+n}({E}_{n};\iota_{n})^{S}.
\end{equation}

For the weak equivalence
\begin{equation}\label{eq:39}
\fil_{n+s+1}( E_{n+1,c};\iota_{n+1})^{S\times\R}\to
\Omegasimp\fil_{n+s+1}(E_{n+1};\iota_{n+1})^{S}
\end{equation}
first apply the ``simplicial loops'' to the diagram defining
$\fil_{n+s+1}(E_{n+1};\iota_{n+1})^{S}$ to get a
homotopy Cartesian square
\[
\begin{CD}
\Omegasimp\fil_{n+s+1}( E_{n+1};\iota_{n+1})^{S}
@>>>\Omegasimp\fil_{n+s+1} 
  \Omega^{\ast}\left(S\times\Delta^{\bullet};\grv\right)^{n+p+1}_{\text{cl}} \\
@VVV @VVV \\
\Omegasimp\sing E_{n+1}^{S} @>>>
\Omegasimp Z^{\ast}\left(S\times\Delta^{\bullet};\grv
\right)^{n+p+1}
\end{CD}
\]
The adjunction between $\sing$ and ``geometric realization'' gives
an isomorphism 
\begin{equation}\label{eq:33}
\sing E_{n+1,c}^{S\times\R}\to \Omegasimp\sing E_{n+1}^{S}.
\end{equation}
A $k$-simplex of $\Omegasimp Z^{\ast}\left(S\times\Delta^{\bullet};\grv
\right)^{n+p+1}$ consists of a sequence of cocycles
\[
c_{0},\dotsc c_{k+1}\in Z^{\ast}\left(S\times\Delta^{k};\grv
\right)^{n+p+1}
\]
satisfying the relations corresponding to~\eqref{eq:34}.  On the other
hand, a $k$-simplex of $Z_{c}^{\ast}\left(S\times\R\times\Delta^{\bullet};\grv
\right)^{n+p+1}$ can be identified with a cocycle 
\[
c\in Z^{\ast}\left(S\times\Delta^{1}\times\Delta^{k};\grv
\right)^{n+p+1}
\]
which vanishes on $S\times\partial\Delta^{1}\times\Delta^{k}$.
Restricting $c$ to the standard triangulation of $\Delta^{1}\times
\Delta^{k}$ leads to a map of simplicial sets
\begin{equation}\label{eq:38}
Z_{c}^{\ast}\left(S\times\R\times\Delta^{\bullet};\grv
\right)^{n+p+1}\to
\Omegasimp Z^{\ast}\left(S\times\Delta^{\bullet};\grv
\right)^{n+p+1},
\end{equation}
which by explicit computation is easily checked to be a weak
equivalence.  Similarly, a $k$-simplex of 
\[
\Omegasimp\fil_{n+s+1} 
  \Omega^{\ast}\left(S\times\Delta^{\bullet};\grv\right)^{n+p+1}_{\text{cl}}
\]
consists of a
sequences of forms
\begin{equation}\label{eq:35}
\omega_{0},\dotsc \omega_{k+1}\in \fil_{n+s+1}\Omega^{\ast}
\left(S\times\Delta^{k};\grv\right)^{n+p+1}_{\text{cl}}
\end{equation}
satisfying the analogue of~\eqref{eq:34}.  A $k$-simplex of 
\[
\fil_{n+s+1}\Omega_{c}^{\ast}\left(S\times\R\times\Delta^{\bullet};\grv
\right)^{n+p+1}_{\text{cl}}
\]
can be identified with a form $\omega$ on $S\times\Delta^{1}\times
\Delta^{k}$ whose restriction to $S\times \partial \Delta^{1}\times
\Delta^{k}$ vanishes, and whose Kunneth components on
$\Delta^{1}\times \Delta^{k} $of degrees $<n+s+1$ vanish.  Restricting
$\omega$ to the simplices in the standard triangulation of
$\Delta^{1}\times \Delta^{k}$ leads to a sequence~\eqref{eq:35}, and
hence to a map of simplicial abelian groups
\begin{equation}\label{eq:37}
\fil_{n+s+1}\Omega_{c}^{\ast}\left(S\times\R\times\Delta^{\bullet};\grv
\right)^{n+p+1}_{\text{cl}} \to
\Omegasimp\fil_{n+s+1} 
  \Omega^{\ast}\left(S\times\Delta^{\bullet};\grv\right),
\end{equation}
which, also by explicit computation, is easily checked to be a weak
equivalence.  

The equivalences~\eqref{eq:33},~\eqref{eq:38}, and~\eqref{eq:37} are
compatible with pullback of cochains, and the inclusion of forms into
cochains, and so patch together, via homotopy pullback, to give the
desired weak equivalence~\eqref{eq:39}.

\subsection{Naturality and Homotopy for
Spectra}\label{sec:naturality-homotopy-1} 
We now turn to some further
constructions on spectra, and the analogues of the results of
\S\ref{sec:naturality-homotopy}.

\begin{defin}
A {\em map} between spectra 
\[
E=\{E_{n},t_{n}^{E}\}\quad\text{and}\quad
F=\{F_{n},t_{n}^{F}\}
\]
consists of a collection of maps $f_{n}:E_{n}\to F_{n}$ , which are
compatible with the structure maps~\eqref{eq:4}, in the sense that the
following diagram commutes
\[
\begin{CD}
E_{n} @>f_{n}>> F_{n} \\
@V t_{n}^{E}VV @VV t_{n}^{F}V \\
\Omega E_{n+1}@>> \Omega f_{n+1} > \Omega F_{n+1}.
\end{CD}
\]
\end{defin}

\begin{rem}
The set of maps is topologized as a subspace of 
\[
\prod F_{n}^{E_{n}}.
\]
\end{rem}

Given $f:E\to F$ and a cocycle 
\[
\iota\in Z^{p}(F;\R),
\]
composition with $f$ gives a map
\[
\check f:\fil_{s}(E;f^{\ast}\iota)^{S}\to
\fil_{s}(F;\iota)^{S}.
\]

\begin{prop}\label{thm:1}
Suppose that $f:E\to F$ is a (weak) homotopy equivalence, and
$\iota\in Z^{p}(F;\R)$ is a cocycle.   Then for  each manifold $S$, the
map 
\[
\check f:\fil_{s}(E;f^{\ast}\iota)^{S}\to
\fil_{s}(F;\iota)^{S}
\]
is a (weak) homotopy equivalence.
\end{prop}

\begin{pf}
This is a consequence of Proposition~\ref{thm:20}.
\end{pf}

As in \S\ref{sec:naturality-homotopy}, Proposition~\ref{thm:1} implies
that a homotopy between maps of spectra $f,g:E\to F$ leads to a
filtration preserving homotopy between maps of differential function spectra.

Given two cocycles $\iota^{1},\iota^{2}\in Z^{p}(E;\R)$, and a cochain
$b\in C^{p-1}(E;\R)$ with $\delta b=\iota^{2}-\iota^{1}$, the
isomorphisms~\eqref{eq:36} 
\begin{align*}
\fil_{s+n}(E_{n};{\iota^{1}_{n}})^{S} &\to 
\fil_{s+n}(E_{n};{\iota^{2}_{n}})^{S}
\end{align*}
fit together to give an isomorphism
\[
\fil_{s}(E ;\iota^{1})^{S}\to
\fil_{s}(E ;\iota^{2})^{S}.
\]
It follows that given a map $f:E\to F$,  and cocycles
$\iota_{E}\in Z^{p}(E;\R)$, $\iota_{F}\in Z^{p}(F;\R)$, and a cochain
$b\in C^{p-1}(E;\R)$ with $\delta b=f^{\ast}\iota_{F}-\iota_{E}$ we
get a map of differential function spectra
\[
\fil_{s} (E;\iota_{E})^{S} \to \fil_{s}(F;\iota_{F})^{S}
\]
which is a weak equivalence if $f$ is.  It also follows 
the (weak) homotopy type of 
\[
\fil_{s}(E,\iota)^{S}
\]
depends only on the cohomology class of $\iota$
and and the homotopy type of $E$ in the sense described in
Remark~\ref{rem:15}.

\begin{rem}
Associated to a space $X$ is its suspension spectrum,
$\Sigma^{\infty}X$, with 
\[
(\Sigma^{\infty}X)_{n}=\varinjlim \Omega^{k}\Sigma^{n+x}X.
\]
One can easily check that the space of maps from $\Sigma^{\infty}X$ to
a spectrum $E$ is simply the space $E_{0}^{X}$, with components, the
cohomology group $E^{0}(X)$.
\end{rem}

\begin{rem}
If $E=\{E_{n} \}$ is a spectrum, one can construct a new spectrum by
simply shifting the indices.  These are known as the {\em shift
suspensions} of $E$, and we will indicate them with the notation
$\Sigma^{k}E$.   To be specific
\[
\left( \Sigma^{k}E\right)_{n}=E_{n+k}.
\]
Note that $k$ may be any integer, and that $\Sigma^{k}E$ is the
spectrum representing the cohomology theory 
\[
X\mapsto E^{\ast+k}(X).
\]
\end{rem}

\begin{rem}
Sometimes the symbol $\Sigma^{k}E$ is used to denote the
spectrum with
\begin{equation}\label{eq:43}
\left(\Sigma^{k}E \right)_{n}=\varinjlim \Omega^{n}\Sigma^{k}E_{n}.
\end{equation}
The spectrum described by~\eqref{eq:43} is canonically homotopy
equivalent, but not equal to $\Sigma^{k}E$.
\end{rem}

\subsection{The fundamental cocycle}\label{sec:fundamental-cocycle}

Recall that for any compact $S$, and any cohomology theory $E$
there is a canonical isomorphism
\begin{equation}\label{eq:48}
E^{\ast}(S)\otimes\R = H^{\ast}(S;\pi_{\ast}E\otimes \R).
\end{equation}
When $E$ is $K$-theory, this isomorphism is given by the Chern character.
The isomorphism~\eqref{eq:48} arises from a universal cohomology class
\[
i_{E}\in H^{0}(E;\pi_{\ast}E\otimes\R) = \varprojlim
H^{n}(E_{n};\pi_{\ast+n}E\otimes\R),
\]
and associates to a map $f:S\to E_{n}$ representing an element of
$E^{n}(S)$ the cohomology class $f^{\ast}i_{n}$, where $i_{n}$
is the projection of $i_{E}$ to $H^{n}(E_{n};\pi_{\ast+n}E)$.  We
will call the class $i_{E}$ a {\em fundamental cohomology
class}--the term used in the case of the Eilenberg-MacLane spectrum,
with $E_{n}=K(\Z,n)$.  A cocycle representative of $i_{E}$ will be
called a {\em fundamental cocycle}.   More precisely,

\begin{defin}
Let $E$ be a spectrum.  A {\em fundamental cocycle} on $E$ is a
cocycle
\[
\iota^{E}\in Z^{0}(E;\pi_{\ast}E\otimes \R)
\]
representing the cohomology class corresponding to the map
\begin{align*}
\pi_{\ast}E &\to\pi_{\ast}E\otimes \R \\
a &\mapsto a\otimes 1
\end{align*}
under the Hurewicz isomorphism
\[
H^{0}(E;\grv)\xrightarrow{\approx}{} \hom(\pi_{\ast}E,\grv).
\]
\end{defin}

Any two choices of fundamental cocycle are cohomologous, and so lead
to isomorphic differential cohomology theories.  Unless otherwise
specified, we will use $\check E(n)$ to denote the differential
cohomology theory associated to a choice of fundamental cocycle.

When $\iota$ is a choice of fundamental cocycle, the
square~\eqref{eq:23} leads to short exact sequences
\begin{equation}\label{eq:45}
\begin{gathered}
0\to E^{q-1}(S;\R/\Z) \to \check E(q)^{q}(S)\to
\Omega^{\ast}(S;\pi_{\ast}E)^{q}_{\text{cl}} \to 0 \\ 
0\to E^{q-1}(S)\otimes \R/\Z \to   \check E(q)^{q}(S)\to A_{E}^q(S)\to 0\\
0\to
\Omega^\ast(S;\pi_{\ast}E)^{q-1}/\Omega^{\ast}(S;\pi_{\ast}E)^{q-1}_0
\to \check E(q)^{q}(S)\to E^q(M)\to 0.
\end{gathered}
\end{equation}
The group $A_{E}^{q}(S)$ is defined by the pullback square
\[
\begin{CD}
A_{E}^{q}(S) @>>> \Omega(S;\pi_{\ast}E\otimes\R)^{q}_{\text{cl}}  \\
@VVV @VVV \\
E^{q}(S) @>>> H^{q}(S;\pi_{\ast}E\otimes\R).
\end{CD}
\]

\begin{rem}
Any cocycle in $E$ is a specialization of the fundamental cocycle in
the sense that given any cocycle $\iota_{1}\in Z^{\ast}(E;\grv)^{0}$,
there is a map (the Hurewicz homomorphism)
\[
t:\pi_{\ast}E\otimes\R\to V
\]
and a cocycle $b\in C^{-1}(E;\grv)$ with 
\[
\delta b=\iota_{1}- t\circ\iota.
\]
\end{rem}

\subsection{Differential bordism}\label{sec:differential-bordism}

\subsubsection{Thom spectra}
\label{sec:thom-spectra}
Let $\R^{\infty}$ be the vector space with basis $\{e_{1},e_{2},\dots
\}$, and 
\[
\mok{k}=\thom\left(\grass{k}{\R^{\infty}},V_{k} \right)
\]
the Thom complex of the natural $k$-plane bundle ($V_{k}$) over
the Grassmannian $\grass{k}{\R^{\infty}}$.  The ``translation map''
\begin{align*}
T:\R^{\infty} &\to\R^{\infty} \\
T(e_{i}) &=e_{i+1}
\end{align*}
induces a map
\begin{align*}
\grass{k}{T}:\grass{k}{\R^{\infty}} &\to\grass{k+1}{\R^{\infty}} \\
V &\mapsto \R e_{1}\oplus V.
\end{align*}
The space $\varinjlim_{k}\grass{k}{\R^{\infty}}$ is a classifying
space for the stable orthogonal group, and will be denoted $BO$.

The pullback of $V_{k+1}$ along $\grass{k}{T}$ is $\R e_{1}\oplus
V_{k}$, giving closed inclusions 
\begin{align*}
S_{k} :\Sigma\mok k & \to \mok{k+1} \\
T_{k} :\mok{k} &\to
\Omega\mok{k+1}.
\end{align*}
The spaces
\[
\mo_{k}=\varinjlim_{n} \Omega^{n}\mok{k+n}
\]
form a spectrum--the {\em unoriented bordism} spectrum, $\mo$.

More generally, given a sequence of closed inclusions
\[
t_{k}:B(k)\to B(k+1),
\]
fitting into a diagram
\begin{equation}\label{eq:46}
\begin{CD}
\cdots @>>> B(k) @>>> B(k+1) @>>> \cdots  \\
@. @VV \xi_{k} V @VV \xi_{k+1} V @. \\
\cdots @>>> \grass{k}{\R^{\infty}} @>>\grass{k}{T}>
\grass{k+1}{\R^{\infty}}@ >>> \cdots
\end{CD}
\end{equation}
the Thom complexes
\[
B(k)^{\xi_{k}^{\ast}V_{k}}
\]
come equipped with closed inclusions
\[
\Sigma B(k)^{\xi_{k}}\to B(k+1)^{\xi_{k+1}}
\]
leading to a spectrum $X=\thom(B;\xi)$ with
\begin{equation}\label{eq:49}
X_{k}=\varinjlim \Omega^{n}\thom(B(n+k);\xi_{n+k}).
\end{equation}

For instance, $B(k)$ might be the space $\sgrass{k}{\R^{\infty}}$ of
oriented $k$-planes in $\R^{\infty}$, in which case the resulting
spectrum is the oriented cobordism spectrum $\mso$.

Spaces $B(k)$ can be constructed from a single map $\xi:B\to BO$, by forming
the homotopy pullback square
\[
\begin{CD}
B(k)@>>> B\\
@V\xi_{k}VV @VV\xi V \\
\grass{k}{\R^{\infty}} @>>> BO.
\end{CD}
\]
The resulting spectrum $\thom\left(B;\xi \right)$ is the {\em Thom
spectrum} of $\xi$.

\subsubsection{Differential $B$-oriented maps}\label{sec:b-oriented-maps}

Let $X=\thom(B;\xi)$ be the Thom spectrum of a map $\xi:B\to BO$, and $S$
a compact manifold.  
Suppose that we are given a cocycle $\iota\in Z^{n}(X;\grv)$ for some
real vector space $\grv$.

\begin{defin}
Let $p:E\to S$ be a (neat) map of manifolds of relative dimension $n$.
A {\em $B$-orientation} of $p$ consists of a neat embedding
$p_{N}=(p,p_{N}'):E\hookrightarrow S\times \R^{N}$, for some $N$, a
tubular neighborhood $W_{N}\hookrightarrow S\times\R^{N}$ and a vector
bundle map
\[
t_{N}:W_{N}\to\xi_{N-n}
\]
classifying $W_{N}$.
\end{defin}

\begin{rem}  Two $B$-orientations are equivalent if they are in the
equivalence relation generated by identifying 
\[
\xi_{N-n}\leftarrow W_{N}\subset S\times\R^{N}
\]
with 
\[
\xi_{N+1-n}\approx\xi_{N+1}\times\R^{1}\leftarrow (W_{N}\times
\R^{1})\approx W_{N+1}\subset S\times\R^{N+1}.
\]
\end{rem}

\begin{defin}
A {\em $B$-oriented map}  is a neat map $p:E\to S$ equipped with an
equivalence class of $B$-orientations.
\end{defin}

\begin{defin}
Let $p:E\to S$ be a (neat) map of manifolds of relative dimension $n$.
A {\em differential $B$-orientation} of $p$ consists of a neat
embedding $p_{N}=(p,p_{N}'):E\hookrightarrow S\times \R^{N}$, for some
$N$, a tubular neighborhood $W_{N}\hookrightarrow S\times\R^{N}$ and a
differential vector bundle map
\[
t_{N}:W_{N}\to(\xi_{N-n};\iota_{N-n})
\]
classifying $W_{N}$.
\end{defin}

\begin{rem}
In short, a differential $B$-orientation of $p:E\to S$ is a lift of
$p$ to a differential $B(N-n)$-oriented embedding $E\subset
\R^{N}\times S$.
\end{rem}

\begin{rem} Two differential $B$-orientations are equivalent if they
are in the equivalence relation generated by identifying
\[
(\xi_{N-n};\iota_{N-n})\xleftarrow{t_{N}}{} W_{N}\subset S\times\R^{N}
\]
with 
\[
(\xi_{N+1-n};\iota_{N+1-n})\xleftarrow{t_{N+1}}{} (W_{N}\times
\R^{1})\approx W_{N+1}\subset S\times\R^{N+1}
\]
where $t_{N}=(c_{N},h_{N},\omega_{N})$, $t_{N+1}=(c_{N}\times
\text{Id},h_{N+1}, \omega_{N+1})$, and
\begin{gather*}
h_{N}=h_{N+1}\slant Z_{\R^{1}} \\
\omega_{N}=\int_{\R^{1}}\omega_{N}.
\end{gather*}
\end{rem}

\begin{defin}
A {\em differential $B$-oriented map} is a neat map 
\[
E\to S
\]
together with an equivalence class of differential $B$-orientations.
\end{defin}

\begin{rem}
\label{rem:2} Let $p:E\to S$ be a differential $B$-oriented map of
relative dimension $n$.  Using \S\ref{sec:thom-complexes},
the construction of Pontryagin-Thom gives a differential function 
\begin{equation}
\label{eq:69}
S\to \left(\thom(B,\xi)_{n}, \iota_{-n} \right).
\end{equation}
In fact the homotopy type of 
\[
\fil_{s}\left(\thom(B,\xi)_{n}, \iota_{-n} \right)^{S}
\]
can be described entirely in terms of differential $B$-oriented maps 
\[
E\to S\times \Delta^{k}.
\]
The proof is basically an elaboration of the ideas of
Thom~\cite{Thom}, and we omit the details.  We will refer to this
correspondence by saying that the differential $B$-oriented map
$p:E\to S$ is {\em classified} by the differential
function~\eqref{eq:69}.
\end{rem}

\begin{rem}
Our definition of differential $B$-orientation depends on many choices
(a tubular neighborhood, a differential function, etc.) which
ultimately affect the push-forward or integration maps to be defined in
\S\ref{sec:push-forward-1}.  These choices are all homotopic, and so
aren't made explicit in the purely topological approaches.  A homotopy
between two differential $B$-orientations can be thought of as a
differential $B$-orientation of $E\times \Delta^{1}\to
S\times\Delta^{1}$, and the effect of a homotopy between choices can be
described in terms of integration along this map.
\end{rem}

\subsubsection{$\bso$-orientations and \orientation{s}}
\label{sec:bso-orient-orient}

Let $\bso=\varinjlim\sgrass{k}{\R^{\infty}}$ be the stable oriented
Grassmannian, as described in \S\ref{sec:thom-spectra}, and choose a
Thom cocycle $U\in Z^{0}(\mso)$.  The resulting notion of
$\bso$-oriented map is a refinement to differential algebraic topology
of the topological notion of an oriented map.  We there have
formulated two slightly different notions of a differential
orientation for a map $E\to S$: that of a differential
$\bso$-orientation that of an \orientation.  For practical purposes
these two notions are equivalent, and we now turn to making precise
the relationship between them.

Fix a manifold $S$.

\begin{defin}
The space of {\em \oriented-maps} of relative dimension $n$ is the
simplicial set whose $k$-simplices are \oriented{} maps $E\to S\times
\Delta^{k}$ of relative dimension $n$.  We denote this simplicial
set $A(S)$, and define $\fil_{t}A(S)$ by restricting the
Kunneth component of the Thom form along the simplices.
\end{defin}

\begin{defin}
The space of {\em differential $\bso$-oriented maps} of relative
dimension $n$, is the simplicial set with $k$-simplices the
differential $\bso$-oriented maps $E\to S\times \Delta^{k}$ of
relative dimension $n$.  We denote this space $B(S)$ and define
$\fil_{s}B(S)$ similarly.
\end{defin}

We also let $\tilde B(S)$ be the space whose $k$-simplices are the
$k$-simplices of $B(S)$, together with a differential cochain
\[
(b,k,0)\in C(N-n)^{N-n-1}(E),
\]
and define $\fil_{t}\tilde B(S)$ by restricting the Kunneth components
as above.  Forgetting about $(b,k,0)$ defines a filtration preserving
function
\[
\tilde B(S)\to B(S)
\]
and taking $(b,k,0)=(0,0,0)$ defines a filtration preserving function 
\[
B(S)\to \tilde B(S).
\]

Suppose that $E/S\times \Delta^{k}$ is a $k$-simplex of $B(S)$.  That
is, we are given
\[
(\xi_{N-n};U_{N-n})\xleftarrow{(c,h,\omega)}{} W\hookrightarrow
S\times \Delta^{k}\times\R^{N}.
\]
Then $W\hookrightarrow S\times \Delta^{k}\times\R^{N}$ together with
$(c^{\ast}U_{N-n},h,\omega)$ is a $k$-simplex of $A(S)$.  This defines
a forgetful map
\[
\fil_{s}B(S)\to \fil_{s}A(S).
\]
We define a map
\[
\fil_{s}\tilde B(S)\to \fil_{s}A(S)
\]
in a similar way, but with the Thom cocycle
\[
(c^{\ast}U_{N-n},h,\omega)+d (b,k,0).
\]
Note that there is a factorization 
\[
\fil_{s} B(S)\to \fil_{s} \tilde B(S)\to \fil_{s} A(S),
\]
and that 
\[
\fil_{s} B(S)\to \fil_{s} \tilde B(S)\to \fil_{s} B(S)
\]
is the identity map.  

\begin{lem}\label{thm:27}
For each $t$, the maps
\begin{gather*}
\fil_{s} \tilde B(S)\to \fil_{s}A(S) \\
\fil_{s} \tilde B(S)\to \fil_{s}B(S) 
\end{gather*}
are acyclic fibrations of simplicial sets.  In particular they are 
homotopy equivalences.
\end{lem}

\begin{cor}\label{thm:26}
For each $s$, the map $\fil_{s}B(S)\to
\fil_{s}A(S)$ is a simplicial homotopy equivalence.\qed
\end{cor}

\begin{pf*}{Proof of Lemma~\ref{thm:27}}
We'll give the proof for
\[
\fil_{s} \tilde B(S)\to \fil_{s}A(S).
\]
The situation with the other map is similar.  By definition, we must
show that a lift exists in every diagram of the form
\begin{equation}
\label{eq:84}
\xymatrix{
\partial \Delta^{k}_{\bullet} \ar[r]\ar[d] & \fil_{s}\tilde B(S)\ar[d] \\
\Delta^{k}_{\bullet}\ar[r] \ar@{-->}[ur]^{\exists} & \fil_{s} A(S).
}
\end{equation}
The bottom $k$-simplex classifies an \oriented{} map 
\begin{gather*}
p:E\to S\times \Delta^{k}\\
W\hookrightarrow S\times \Delta^{k}\times \R^{N}\\
(c,h,\omega)\in Z(N-n)^{N-n}(W), 
\end{gather*}
and the top map gives compatible 
$\bso$-orientations to the boundary faces
\begin{gather*}
\partial_{i}E\to S\times\partial_{i}\Delta^{n} \\
(\xi_{N-n};U_{N-n})\xleftarrow{(f_{i},h_{i},\omega_{i})}
\partial_{i}W\hookrightarrow S\times \partial_{i}\Delta^{k}\times\R^{N-n} \\
\end{gather*}
and compatible differential cocycles
\[
(b_{i},k_{i},0)\in C(N-n)^{N-n}(f_{i}^{\ast}\xi),
\]
satisfying 
\[
(c,h,\omega)\vert_{\partial_{i}W}
=(f_{i}^{\ast}U_{N-n},h_{i},\omega_{i})+d(b_{i},k_{i},0).
\]
Write 
\[
\partial E=\bigcup_{i} \partial_{i}E\qquad
\partial W=\bigcup_{i} \partial_{i}W=\bigcup_{i} W\vert_{\partial_{i}E}.
\]
The compatibility conditions imply that the functions $f_{i}$ together
form a vector bundle map
\[
\partial f:\partial W\to \xi_{N-n}.
\]
By the universal property of $\bso(N-n)$, $\partial f$ extends to a
vector bundle map
\[
f:W\to \xi_{N-n},
\]
inducing the same (topological) orientation on $W$ as the one given by
the Thom cocycle $(c,h,\omega)$.  In particular, this implies that
$f^{\ast}U(N-n)$ and $c$ represent the same cohomology class.

To construct a lift in~\eqref{eq:84}, it suffices to find a cocycle
$h'\in C^{N-n-1}(\bar W;\R)$ with 
\begin{equation}
\begin{gathered}
\delta h'=\omega-f^{\ast}U_{N-n} \\
h'\vert_{\partial_{i}W}=h_{i},
\end{gathered}
\end{equation}
and a differential cochain $(b,k,0)\in C(N-n)^{N-n-1}(\bar W)$
satisfying
\begin{equation}\label{eq:90}
\begin{gathered}
(b,k,0)\vert_{\partial_{i}W}=(b_{i},k_{i},0)\\
(c,h,\omega)=(f^{\ast}U_{N-n},h',\omega)+d(b,k,0).
\end{gathered}
\end{equation}
To construct $h'$, first choose any $h''\in C^{N-n-1}(\bar W;\R)$
with $h''\vert_{\partial_{i}W}=h_{i}$ for all $i$.  Then the cocycle 
\begin{equation}\label{eq:89}
\delta h''-\omega+f^{\ast}U_{N-n} 
\end{equation}
represents a relative cohomology class in 
\[
H^{N-n}(\bar
W,\bar{\partial W};\R)\xrightarrow{\approx}{}
H^{0}(E,\partial E;\R).
\]
Now this latter group is simply the ring of real-valued functions on
the set of path components of $E$ which do not meet $\partial E$.  But
on those components, the expression~\eqref{eq:89} represents $0$,
since the image of $f^{\ast}U_{N-n}$ in $Z^{N-n}(\bar W;\R)$ and
$\omega$ represent the same cohomology class.  It follows
that~\eqref{eq:89} is the coboundary of a relative cochain
\[
h'''\in C^{N-n-1}(\bar W,\bar{\partial W};\R).
\]
We can then take 
\[
h'=h''-h'''.
\]
A similar argument, using the fact that $H^{N-n-1}(\bar
W,\bar{\partial W};\rz)=0$ leads to the existence of $(b,k,0)$
satisfying~\eqref{eq:90}.
\end{pf*}

The fact that acyclic fibrations are preserved under change of base
makes the result of Proposition~\ref{thm:27} particularly convenient.
For instance, a single \oriented{} map $E/S$ defines a $0$-simplex in
$A(S)$, and Proposition~\ref{thm:27} asserts that the inverse image of
this $0$-simplex in $\fil_{s}\tilde B(S)$ is contractible.  This means that
for all practical purposes the notions of \oriented{} map and
$\bso$-oriented are equivalent.  This same discussion applies to many
combinations of geometric data.

\subsection{Integration}\label{sec:push-forward-1} 

In this section we define {\em integration} or {\em push-forward} in
differential cohomology theories.  We will see that because of the
results of \S\ref{sec:bso-orient-orient} and
Appendix~\ref{sec:comp-checkhznks-chsn} this recovers our theory of
integration of differential cocycles (\S\ref{subsec-integration}) in
the case of ordinary differential cohomology.  At the end of this
section we discuss the integration in differential $K$-theory.

The topological theory of integration is simply the interpretation in
terms of manifolds, of a map from a Thom spectrum to another spectrum.
Thus let $X=\thom(B;\xi)$ be a Thom spectrum, and $R$ a spectrum.
Consider 
\[
\thom(B\times R_{m};\xi\oplus 0)=X\wedge (R_{m})_{+}
\]
and suppose there is given a map
\[
\mu:X\wedge (R_{m})_{+}\to \Sigma^{m}R.
\]
In geometric terms, a map 
\[
S\to \left(X\wedge (R_{m})_{+} \right)_{-n}
\]
arises from a $B\times R_{m}$-oriented map $E\to S$ of relative
dimension $n$, or what amounts to the same thing, a $B$-oriented map,
together with a map $x:E\to R_{m}$, thought of as a ``cocycle''
representing an element of $R^{m}(E)$.  The composition
\[
S\to \left(X\wedge (R_{m})_{+} \right)_{-n}\xrightarrow{\mu}{}
\left(\Sigma^{m}R \right)_{-n}\approx
\left(R \right)_{m-n}
\]
is a map $y$ representing an element of $R^{m-n}(S)$.  Thus the
geometric interpretation of the map $\mu$ is an operation which
associates to every $B$-oriented map $E/S$ of relative dimension $n$
and every $x:E\to R_{m}$, a map $y:S\to R_{m-n}$.  We think
of $y$ as the {\em integral} of $x$.

\begin{rem}
The notation for the spectra involved in discussing integration tends
to become compounded.  Because of this, given a spectrum $R$ and a
cocycle $\iota\in C^{\ast}(R;\grv)$ we will denote pair
$(R_{n};\iota_{n})$ as
\[
(R_{n};\iota_{n})=(R;\iota)_{n}.
\]
\end{rem}

For the {\em differential} theory of integration, suppose that
cocycles $\iota_{1}$, $\iota_{2}$ have been chosen so as to refine
$\mu$ to a map of differential cohomology theories
\[
\check{\mu}:\left(X\wedge (R_{m})_{+} ;\iota_{1}\right)\to
\left(R;\iota_{2}\right).
\]  
The map $\check{\mu}$ associates
to every differential $B$-oriented map $E\to S$ of
relative dimension $n$ together with a differential function $x:E\to
R_{m}$, a differential function $y:S\to R_{m-n}.$ We will refer to
$y$ as the {\em push-forward} of $x$ and write $y=p_{!}(x)$, or
\[
y=\int^{\mu}_{E/S}x
\]
to emphasize the analogy with integration of differential cocycles.
When the map $\mu$ is understood, we will simply write
\[
y=\int_{E/S}x.
\]

\begin{rem}
The construction $p_{!}$ is a map
of differential function spaces 
\[
\int^{\mu}:\left(X\wedge (R_{m})_{+};\iota_{1} \right)_{-n}^{S}\to
\left(R ;\iota_{2}\right)_{m-n}^{S}.
\]
The weight filtration of $\int^{\mu}$ can be controlled geometrically,
but in the cases that come up in this paper it is something that can
be computed after the fact.
\end{rem}

The push-forward construction is compatible with changes in $R$.
Suppose $f:P\to R$ is a map of spectra, and that there are maps 
\begin{align*}
\mu_{P} :X\wedge (P_{m})_{+} \to \Sigma^{m}P \\
\mu_{R} :X\wedge (R_{m})_{+} \to \Sigma^{m}R
\end{align*}
and a homotopy 
\[
H:X\wedge (P_{m})_{+}\wedge \Delta^{1}_{+} \to \Sigma^{m}R
\]
between the two ways of going around
\begin{equation}\label{eq:52}
\begin{CD}
X\wedge (P_{m})_{+} @> \mu_{P}>>  \Sigma^{m}P \\
@VVV @VVV \\
X \wedge (R_{m})_{+} @>\mu_{R} >> \to \Sigma^{m}R.
\end{CD}
\end{equation}
Suppose also that cocycles 
\[
\iota_{P}\in Z^{k}(P;\grv)\qquad
\iota_{R}\in Z^{k}(R;\grw)
\]
have been chosen so that $f$ refines to a map of differential
cohomology theories
\[
\check f:\left(P;\iota_{P} \right)\to
\left(R;\iota_{R} \right).
\]
Then differential pushforward maps
\begin{align*}
\int^{\mu_{P}}(\slot) : (X\wedge (P_{m})_{+};\mu_{P}^{\ast}\iota_{P})_{-n}^{S}
\to (P;\iota_{P} )_{m-n}^{S}\\
\int^{\mu_{R}}(\slot) : (X\wedge (R_{m})_{+};\mu_{R}^{\ast}\iota_{R})_{-n}^{S}
\to (R;\iota_{R} )_{m-n}^{S}
\end{align*}
are defined, and the results of \S\ref{sec:naturality-homotopy-1} give a
canonical (weight filtration preserving) homotopy between 
\[
\check f\circ \int^{\mu_{P}}(\slot)\quad\text{and}\quad
\int^{\mu_{P}}\check f_{\ast}(\slot).
\]

\begin{rem}\label{rem:20}
In case $B=\text{pt}$, the Thom spectrum $X=S^{0}$ represents
framed cobordism, and every spectrum $R$ comes equipped with a
natural map
\[
X\wedge (R_{m})_{+}\to \Sigma^{m}R,
\]
namely the structure map of the spectrum.  This map is compatible with
every map of spectra $P\to R$ in the sense described above.  {\em Thus
pushforward maps exist in every cohomology theory for a {\em differential
framed} map $E\to S$, and these are compatible with all maps between
the cohomology theories.}  
\end{rem}

\begin{eg} {\em Ordinary differential cohomology.}\label{eg:12}
Let $X$ be the oriented bordism spectrum $\mso$ and $R$ the
Eilenberg-MacLane spectrum $H\Z$, with $R_{m}=K(\Z,m)$.  Choose a Thom
cocycle $U\in Z^{0}(\mso;\Z)\subset Z^{0}(\mso;\R)$ and a fundamental
cocycle $x\in Z^{0}(H\Z;\Z)\subset Z^{0}(H\Z;\R)$.  We take
\[
\mu: \mso\wedge K(\Z,m)_{+}\to \Sigma^{m}H\Z
\]
to be a map representing
\[
\iota = u\cup x_{m}\in Z^{m}(MSO\wedge K(\Z,m)_{+};\Z)\subset
Z^{m}(MSO\wedge K(\Z,m)_{+};\R).
\]

Let 
\begin{gather*}
E\to S \\
(\xi_{N-n};U_{N-n}) \xleftarrow{\check U}{}W\subset S\times\R^{N}
\end{gather*}
be a differential $\bso$-oriented map.  Given a differential function
\[
\check x:E\to \left(K(\Z,m);\iota_{m} \right)
\]
form the differential $\bso\times K(\Z,m)$-oriented map
\begin{equation}
\label{eq:85}
(\xi_{N-n}\oplus 0;U_{N-n}\cup \iota_{m}) \xleftarrow{\check U\cup
\check x}{}W\subset S\times\R^{N}.
\end{equation}
The map~\eqref{eq:85} is classified by a differential function
\[
S\to \left(\mso\wedge K(\Z,m)_{+}; U\cup \iota_{m}\right)_{-n},
\]
and composing this with $\int^{\mu}$ gives a differential function
\[
\int^{\mu}_{E/S}\check x:S \to \left(K(\Z,m-n);\iota_{m-n} \right).
\]
Replacing $E\to S$ with $E\times \Delta^{k}\to S\times \Delta^{k}$
leads to a map of differential function complexes
\[
\int^{\mu}:\fil_{s}\left(K(\Z,m);\iota_{m} \right)^{E}
\to \fil_{s}\left(K(\Z,m-n);\iota_{m-n} \right)^{S}.
\]

On the other hand, associated to the differential $\bso$-orientation
of $E/S$ is an \orientation{} of $E/S$, and, as described in
Appendix~\ref{sec:comp-checkhznks-chsn} the differential function complexes
\[
\fil_{s}\left(K(\Z,m);\iota_{m} \right)^{E}\quad\text{and}\quad
\fil_{s}\left(K(\Z,m-n);\iota_{m-n} \right)^{S}
\]
are homotopy equivalent to the simplicial abelian groups underlying
the chain complexes
\begin{gather*}
\chcocycle{m-s}^{m}(E)
\xleftarrow{}{}
\chcochain{n}^{m-1}(E)\dots 
\xleftarrow{}{}
\chcochain{m-s}^{0}(E) \quad\text{and}\\
\chcocycle{m-n-s}^{m-n}(S)
\xleftarrow{}{}
\chcochain{m-n-s}^{m-n-1}(S)\dots 
\xleftarrow{}{}
\chcochain{m-n-s}^{0}(S),
\end{gather*}
respectively.  As the reader can check, our conventions have been
chosen so that integration of differential cocycles and integration of
differential functions to an Eilenberg-MacLane space agree under
this correspondence.
\end{eg}

\begin{eg} {\em The case of differential $K$-theory.}\label{eg:13}
We now turn to the pushforward map in differential $K$-theory.  Let
$p:E\to S$ be a map of relative dimension $2n$ with a ``differential
$\spinc$-structure'' on the relative normal bundle.  The theory
described above gives a map
\[
p_{!}:\check K^{0}(E)\to \check K^{-2n}(S)\approx \check K^{0}(S).
\]
Actually $p_{!}$ sends a {\em differential
vector bundle} $\bar V$ on $E$, i.e., a differential function
\[
\bar V=(c,h,\omega):E\to \left(\fred;\iota \right),
\]
to a differential vector bundle on $S$.  In this section we will
describe this map in some detail.

If we drop the ``differential'' apparatus, we get the topological
push-forward map
\[
p_{!}^{\text{top}}: K^{0}(E)\to  K^{0}(S).
\]
We remind the reader how that goes.

To define the topological pushforward, embed $E$ in $\R^{2N}\times S$,
and let 
\[
\pi:\nu\to E
\]
be the normal bundle.   Clifford multiplication
\[
\nu\times S^{+}(\nu)\to \nu\times S^{-}(\nu)
\]
defines a $K$-theory class 
\[
\Delta\in K^{0}(\nu,\nu-E).
\]

Let 
\[
\nu\approx D\subset \R^{2N}\times S
\]
be a tubular neighborhood of $E$.  We have canonical isomorphisms
\[
K^{0}(\nu,\nu- E)\approx 
K^{0}(D,D- E)\xleftarrow{\approx}{}
K^{0}(\R^{2N}\times S,R^{2N}\times S- E).
\]
If $V$ is a vector bundle representing an element of $K^{0}(E)$,
then 
\[
\Delta\cdot\pi^{\ast}V
\]
is an element of $K^{0}(\nu,\nu-E)$, and the
topological pushforward of $V$ is the image of $\Delta\cdot\pi^{\ast}V$
under
\[
K^{0}(\nu,\nu- E)\approx K^{0}(\R^{2N}\times S,R^{2N}\times
S- E) 
\to K^{0}(\bar{\R}^{2N}\times S)\xrightarrow{}{}
K^{0}(S).
\]

We can now turn to differential $K$-theory.  To get a 
``differential $p_{!}$'' we need a differential
$\spinc$-structure--that is, a refinement of $\Delta$ to 
\[
\hat{\Delta}=(c_{\Delta}, h_{\Delta}, \omega_{\Delta})
\]
where
\[
c_{\Delta}:\bar\nu\to \fred
\]
is a map classifying $\Delta$, and 
\[
\quad h_{\Delta}\in C^{\text{odd}}(\nu;\R),\quad
\omega_{\Delta}\in \Omega^{\text{ev}}_{\text{cl}}(\bar\nu)
\]
satisfy
\[
\delta h_{\Delta}=\omega_{\Delta}-c_{\Delta}^{\ast}\iota,
\]
where $\iota$ is a choice of Chern character cocycles as in
\S\ref{sec:interl-diff-k}.  We now imitate the topological
construction.  Let $\check{V}=(c,h,\omega)$ represent an element of
$\check K^{0}(E)$.  Then $\pi^{\ast}\check V\in \check{K}^{0}(\nu)$,
and we define $p_{!}\check V$ to be the image of $\check \Delta\cdot
\pi^{\ast}\check V$ under
\[
\check K^{0}(\bar{\nu})\approx \check K^{0}(\R^{2N}\times S,R^{2N}\times
S- D)
\to \check K^{0}(\bar{\R}^{2N}\times S)\rightarrow
\check K^{0}(S).
\]
\end{eg}

\section{The topological theory}
\label{sec-topology}

\renewcommand{\bsobwu}{\bso\langle\beta\nu_{2k} \rangle}
\renewcommand{\msobwu}{\mso\langle\beta\nu_{2k} \rangle}
\renewcommand{\msobwuh}{\msobwu_{h\Z/2}}

\subsection{Proof of Theorem~\ref{thm-main}}\label{sec:proof-theorem-main}

We now turn to the proof of Theorem~\ref{thm-main}.  Our approach will
be to reformulate the result in terms of a transformation of
differential cohomology theories.  In fact we need only construct a
transformation of {\em cohomology theories}--the refinement to the
differential versions is more or less automatic.  It is rather easy to
show that many transformations exist which will do the job (see
Proposition~\ref{thm:t-22}).  But it turns out that results of
Milgram~\cite{milgram74:_surger}, and Morgan--Sullivan~\cite{morgan74}
make it possible to single one out.  In this first section we prove
Theorem~\ref{thm-main}, but do not single out a particular $\kappa$.
The remaining sections are devoted to making a particular choice.

The input of the functor $\kappa$ is an \oriented{} map $E\to S$
equipped with a ``twisted'' differential cocycle giving an integral
Wu-structure, and a differential $L$-cocycle.  We first turn to
interpreting this data in terms of cobordism.

Choose a map
\[
\bso\to K(\Z/2,2k)
\]
representing the Wu class, $\nu_{2k}$, of the
universal bundle, and $\bsobwu$ by the homotopy pullback square
\[
\begin{CD}
\bsobwu @> \lambda >> K(\Z,2k) \\
@VVV @VVV \\
\bso @>>> K(\Z/2,2k).
\end{CD}
\]
Let $\msobwu$ be the associated Thom spectrum.  The group
\[
\msobwu^{i-4k}(S)
\]
is the cobordism group of maps $E\to S$ of relative dimension $(4k-i)$
equipped with an orientation of the relative stable normal bundle, and
an integral Wu-structure.  

The Thom spectrum of the vector bundle classified by the projection map
\[
\bsobwu\times K(\Z,2k)\to BSO
\]
is the smash product
\[
\msobwu\wedge K(\Z,2k)_{+},
\]
The group
\[
\left(\msobwu\wedge K(\Z,2k)_{+} \right)^{i-4k}(S)
\]
is the cobordism group of maps $E\to S$ of relative dimension $(4k-i)$
equipped with an orientation of the relative stable normal bundle,
an integral Wu-structure, and a cocycle $x\in Z^{2k}(E)$.

The functors $\kappa$ and $q^{\wustruct}$ in Theorem~\ref{thm-main}
are constructed out of maps from $\msobwu$ and $\msobwu\wedge
K(\Z,2k)_{+}$ to some other spectrum.  There are many possible choices
of this other spectrum, and we will work with one which is
``universal'' in the sense that it {\em receives} a map from any other
choice.  This universal spectrum theory is the {\em Anderson dual} of
the sphere~\cite{anderson69:_univer_k,yosimura75:_univer_cw}, which we
denote $\tilde I$.  Roughly speaking, Anderson duality is like
Pontryagin duality, and for any spectrum $X$, any homomorphism
$\pi_{n}X\to \Z$ is represented by a map $X\to\Sigma^{n}\tilde I$.
More precisely, there is a (splittable) short exact sequence
\begin{equation}\label{eq:1}
\ext\left(\pi_{n-1}X,\Z \right)\rightarrowtail
[X,\Sigma^{n}\tilde I]=\tilde I^{n}(X)\twoheadrightarrow
\hom\left(\pi_{n}X,\Z \right).
\end{equation}
The spectrum $\tilde I$
is defined in Appendix~\ref{sec:picard-categ-anders}.  By~\eqref{eq:1}, one has
\[
\tilde I^{k}(\text{pt})= 
\begin{cases}
0 &\qquad \text{if $k<0$} \\
\Z &\qquad \text{if $k=0$} \\
0 &\qquad  \text{if $k=1$} \\
\hom(\pi_{k-1}^{\text{st}}S^{0},\qz) &\qquad k>1
\end{cases}
\]
where 
\[
\pi_{k-1}^{\text{st}}S^{0}= \varinjlim \pi_{k-1+N}S^{N}
\]
denotes the $(k-1)^{\text{st}}$ stable homotopy group of the sphere.
It follows from the Atiyah-Hirzebruch spectral sequence that for any
space $M$, one has
\begin{align*}
\tilde I^{0}(M)&= H^{0}(M;\Z) \\
\tilde I^{1}(M)&= H^{1}(M;\Z) \\
\tilde I^{2}(M)&= H^{2}(M;\Z)\times H^{0}(M;\Z/2).
\end{align*}
In particular, the group $\tilde I^{2}(M)$ contains the group of
complex line bundles.  In fact, thinking of $\Z/2$ as the group
$\pm1$, the whole group $\tilde I^{2}(M)$ can be identified with the
group of {\em graded} line bundles, the element of $H^{0}(M;\Z/2)$
corresponding to the degree.  Since the map $H\Z\to \tilde I$ is a
rational equivalence, we can choose a cocycle
\[
\iotai\in Z^{0}\left(\tilde I;\R \right)
\]
restricting to a fundamental cocycle in $Z^{0}\left(H\Z;\R \right)$.
This makes $\tilde I$ into a differential cohomology theory.  The
differential cohomology group $\tilde I^{2}(2)(M)$ can be interpreted
as the group of horizontal isomorphism classes of graded
$U(1)$-bundles with connection.

To a first approximation, the functor $\kappa$ is derived from a map
of spectra 
\begin{equation}\label{eq:t-64}
\msobwu\to \Sigma^{4k}\tilde I.
\end{equation}
But there is one more complication.  In order to establish the symmetry
\begin{equation}\label{eq:t-63}
\kappa(-\wustruct)\approx
\kappa(\wustruct)
\end{equation}
we need to put a $\Z/2$-action on $\msobwu$ in such a way that it
corresponds to the symmetry 
\[
\wustruct\mapsto -\wustruct.
\]
In \S\ref{sec:quadratic-functions} we will describe a
$\Z/2$-equivariant Eilenberg-MacLane space, $K(\Z(1),2k)$ with the
following properties:
\begin{textList}
\item The involution 
\[
\tau:K\left(\Z(1),2k \right)\to K\left(\Z(1),2k \right)
\]
has degree $-1$;

\item there is an equivariant map 
\[
K\left(\Z(1),2k \right)\to 
K(\Z/2,2k),
\]
with $\Z/2$ acting trivially on $K(\Z/2,2k)$, corresponding to
reduction modulo $2$.
\end{textList}
Using this, we define a $\Z/2$-equivariant $\bsobwu$ by the
homotopy pullback square
\[
\begin{CD}
\bsobwu @> \lambda >> K\left(\Z(1),2k \right) \\
@VVV @VVV \\
BSO @>>\nu_{2k}> K\left(\Z/2,2k \right).
\end{CD}
\]
The associated Thom spectrum $\msobwu$ then acquires a $\Z/2$-action,
and the existence of the symmetry isomorphism~\eqref{eq:t-63} can be
guaranteed by factoring~\eqref{eq:t-64} through a map
\begin{equation}\label{eq:t-65}
\msobwuh \to \Sigma^{4k}\tilde I.
\end{equation} 
We have used the notation
\[
X_{h\Z/2}=X\wedge_{\Z/2}E\Z/2_{+}
\]
to denote the {\em homotopy
orbit spectrum of} of a spectrum $X$ with a $\Z/2$-action.

The spectrum $\msobwuh$ is also a Thom spectrum
\[
\msobwuh=\thom\left(\bsobwu\times_{\Z/2}E\Z/2;V \right),
\]
with $V$ the stable vector bundle classified by
\begin{equation}\label{eq:88}
\bsobwu\times_{\Z/2}E\Z/2\to
\bso\times B\Z/2\to\bso.
\end{equation}
The group
\[
\pi_{4k}\msobwuh
\]
is the cobordism group of $4k$-dimensional oriented manifolds $M$,
equipped with a map $t:\pi_{1}M\to\Z/2$ classifying a local
system $\Z(1)$, and a cocycle $\lambda\in Z^{2k}(M;\Z(1))$ whose mod
$2$-reduction represents the Wu-class $\nu_{2k}$.  

Fix a cocycle $L_{4k}\in Z^{4k}(\bso;\R)$ representing the component
of the Hirzebruch $L$-polynomial of degree $4k$.  By abuse of notation
we will write
\[
U\in \Z^{0}(\msobwuh;\R)\qquad
L_{4k}\in \Z^{0}(\bsobwu\times_{\Z/2}E\Z/2;\R)
\]
for the pullback along the maps derived from~\eqref{eq:88} of the Thom
cocycle and $L_{4k}$, respectively.

The group of maps~\eqref{eq:t-65} sits in a (splittable) short exact
sequence
\begin{multline}
\label{eq:87}
\ext\left(\pi_{4k-1}\msobwuh,\Z \right)\rightarrowtail
\tilde I^{4k}\left(\msobwuh \right)   \\
\twoheadrightarrow
\hom\left(\pi_{4k}\msobwuh,\Z \right),
\end{multline}
whose rightmost terms is the group of integer-valued cobordism
invariants of manifolds $M$  of the type just described.  One example
is the signature $\sigma=\sigma(M)$ of the non-degenerate bilinear
form
\begin{equation}\label{eq:t-66}
B(x,y)=\int_{M}x\cup y:H^{2k}\left(M;\Q(1) \right)\times
H^{2k}\left(M;\Q(1) \right)\to \Q.
\end{equation}
Another is 
\[
\int_{M}\lambda^{2}.
\]
By definition, $\lambda$ is a characteristic element for the bilinear
form~\eqref{eq:t-66}, and so
\[
\int_{M}\lambda^{2}\equiv \sigma(M) \mod 8.
\]

\begin{prop}\label{thm:t-22}
Any map of spectra
\[
\kappa:\msobwu_{h\Z/2}\to 
\Sigma^{4k}\tilde I
\]
whose underlying homomorphism
\[
\pi_{4k}\msobwu_{h\Z/2}\to\Z
\]
associates to $M^{4k}$, $\lambda\in Z^{2k}(M;\Z(1))$, the integer
\begin{equation}\label{eq:t-27}
\frac18\left(\int_{M}\lambda^{2}- \sigma \right),
\end{equation}
gives a family of functors
\[
\kappa_{E/S}:\chcat^{2k}_{\nu}(E) \to \chcat^{i}(S)
\]
having the properties listed in Theorem~\ref{thm-main} and
Corollary~\ref{thm:28}.
\end{prop}

Because of the sequence~\eqref{eq:87} the set of maps $\kappa$
satisfying the condition of Proposition~\ref{thm:25} is a non-empty
principal homogeneous space for 
\[
\ext\left(\pi_{4k-1}\msobwuh,\Z.
\right)
\]

Before going through the proof we need one more topological result.
The change of an integral Wu-structure by a $2k$-cocycle is
represented by a map
\[
\left(\bsobwu\times K(\Z,2k) \right)\times_{\Z/2}E\Z/2\to
\bsobwu\times_{\Z/2}E\Z/2.
\]
We'll write the induced map of Thom spectra as 
\[
(\lambda-(2)x):\msobwukhp\to \msobwuh.
\]
Of course, the map $(\lambda-(2)x)-(\lambda)$ factors through
\[
\msobwukh\to \msobwuh.
\]
We will be interested in 
\[
b(x,y)=(\lambda-(2)(x+y))- (\lambda-(2)x) -
(\lambda-(2)y) + (\lambda)
\]
which is a map 
\[
\msobwukkh\to \msobwuh.
\]
\begin{lem}\label{thm:25}
If $\kappa:\msobwu_{h\Z/2}\to \Sigma^{4k}\tilde I$
is any map satisfying the condition of Proposition~\ref{thm:t-22}, then
the following diagram commutes up to homotopy
\[
\begin{CD}
\msobwukkh @> U\cup x\cup y>> \Sigma^{4k}\hz \\
@V b(x,y) VV @VVV \\
\msobwu @>> \kappa > \Sigma^{4k}\tilde I
\end{CD}
\]
\end{lem}

\begin{pf}
Since 
\[
\pi_{4k-1}\msobwu\wedge K(\Z(1),2k)\wedge
K(\Z(1),2k)_{h\Z/2}=0,
\]
it suffices to check that the two ways of going around the diagram
agree after passing to $\pi_{4k}$.   
The clockwise direction associates to $(M,\lambda,x,y)$ the integer
\[
q^{\lambda}(x+y)-q^{\lambda}(x)-q^{\lambda}(y)+q^{\lambda}(0),
\]
where we have written
\[
q^{\lambda}(x)=\kappa(\lambda-(2)x)=\tfrac12\int_{M}(x^{2}-x\lambda).
\]
The counter-clockwise direction is 
\[
\int_{M}x\cup y.
\]
The result follows easily. 
\end{pf}

\begin{pf*}{Proof of Proposition~\ref{thm:t-22}}
By construction, $\tilde I$ comes equipped with a fundamental cocycle
$\iotai\in Z^{0}(\tilde I;\R)$ which necessarily restricts to a
fundamental cocycle in $\hz$.  Now the signature of the bilinear
form~\eqref{eq:t-66} coincides with the signature of $M$.  Indeed, let
$\tilde M\to M$ be the double cover classified by the homomorphism
$t:\pi_{1}M\to\Z/2$.  Then by the signature theorem, $\sigma(\tilde
M)=2\sigma(M)$, and the claim then follows by decomposing
$H^{2k}(\tilde M;\Q)$ into eigenspaces under the action of $\Z/2$.  It
follows that the cohomology class of $\kappa^{\ast}\iotai$ coincides
with that of
\begin{equation}\label{eq:70}
\alpha=U\cdot\frac{(\lambda^{2}-L_{4k})}{8},
\end{equation}
so after choosing a cochain whose coboundary is the difference between 
$\kappa^{\ast}\iotai$ and~$\alpha$,  we have a map of
differential cohomology theories
\[
\left(\msobwuh; \alpha \right)\to
\Sigma^{4k}\left(\tilde I;\tilde\iota \right).
\]
Since $H^{4k-1}\left(\msobwuh;\R \right)=0$, this choice of cochain has
no effect on the maps of fundamental groupoids we derive from it.

Choose a point $b\in E\Z/2$.  We use $b$ to define a map
\[
i_{b}:\msobwu=
\msobwu\wedge\left\{b \right\}_{+}\hookrightarrow
\msobwu\wedge E\Z/2_{+}.
\]
And to keep the notation simple, we will not distinguish in notation
between $\alpha$ and $\iota_{b}^{\ast}\alpha$.  Suppose that $E/S$ is
an \oriented{} map of manifolds of relative dimension $(4k-i)$.  The
functor $\kappa_{E/S}$ is constructed from the map of fundamental
groupoids of differential function complexes
\begin{equation}\label{eq:68}
\pi_{\le 1}\fil_{0}\left(\msobwu_{i-4k};\alpha \right)^{S}\to
\pi_{\le 1}\fil_{0}\left(\tilde I_{i};\iotai \right)^{S}.
\end{equation}
In principle, a differential integral Wu-structure on $E/S$ defines a
$0$-simplex of the differential function space
\[
\fil_{0}\left(\msobwu_{i-4k};\alpha\right)^{S},
\]
in such a way as to give a functor
\begin{equation}\label{eq:67}
\chcat^{2k}_{\nu}(E) \to \pi_{\le
1}\fil_{0}\left(\msobwu_{i-4k};\alpha\right)^{S}.
\end{equation}
The functor $\kappa_{E/S}$ is just the composition of~\eqref{eq:67}
with~\eqref{eq:68}.  We have said ``in principle'' because the
geometric data we specified in the statement of Theorem~\ref{thm-main}
is not quite the data which is classified by a differential function
from $S$ to $\left(\msobwu;\alpha\right)_{i-4k}$.  To remedy this, we
use the technique described in \S\ref{sec:bso-orient-orient} to
produce a diagram
\[
\chcat^{2k}_{\nu}(E) \xleftarrow{\sim}{}
\chcat^{2k}_{\nu}(E)_{\text{geom}} \to 
\pi_{\le 1}
\fil_{0}\left(\msobwu_{i-4k};\alpha\right)^{S},
\]
in which the left map is an equivalence of groupoids.  The functor
$\kappa_{E/S}$ is then constructed as described, after choosing an
inverse to this equivalence.  

The category $\chcat^{2k}_{\nu}(E)_{\text{geom}}$ is the fundamental
groupoid of a certain simplicial set $\mathcal S_{\nu}^{2k}(E)$.  To
describe it, first note that by the results of
\S\ref{sec:bso-orient-orient} we may assume that the \orientation{} of
$E/S$ comes from a differential $\bso$-orientation
\[
(\xi_{N-i+4k},U_{N-i+4k})\xleftarrow{(c,h,\omega)} W\subset S\times \R^{N},
\]
and that the differential cocycle\footnote{We are using the symbol
$L_{4k}$ to denote both the universal signature cocycle and the chosen
refinement to a differential cocycle on $E$.  We hope this causes no
confusion.}  $L_{4k}$ refines the map classifying $W$ to a
differential function
\[
E\xrightarrow{(c_{0},h_{0},\omega_{0})}{} (\bso;L_{4k}).
\]
A $k$-simplex of $\mathcal
S_{\nu}^{2k}(E)$ consists of a differential function
\begin{equation}\label{eq:91}
E\times \Delta^{k}\xrightarrow{(c_{1},h_{1},\omega_{1})}{}
\left(\bsobwu;\lambda \right)
\end{equation}
of weight filtration $0$, for which the map $c_{1}$ fits into the
commutative diagram
\[
\xymatrix{
E\times\Delta^{k}\ar[rr]^{c_{1}}\ar[d] & & \bsobwu  \ar[d] \\
E \ar[r] &\bso(N-i+4k) \ar[r]&\bso.
}
\]
The functor
\begin{equation}\label{eq:92}
\chcat^{2k}_{\nu}(E)_{\text{geom}}=\pi_{\le 1}\mathcal
S_{\nu}^{2k}(E)\to \chcat^{2k}_{\nu}(E)
\end{equation}
sends~\eqref{eq:91} to the slant product of
$(c_{1}^{\ast}\lambda,h_{1},\omega_{1})$ with the fundamental class of
$\Delta^{k}$.  Using the results of Appendix~\ref{sec:comp-checkhznks-chsn} 
one checks that~\eqref{eq:92} is an equivalence of groupoids.

Write
\begin{align*}
\Check U &= \left(c^{\ast}U_{N-i+4k},h,\omega \right)\\
\Check \lambda &= \left(c_{1}^{\ast}\lambda,h_{1},\omega_{1} \right).
\end{align*}
The lift $c_{1}$ of the map classifying $W$, and the differential
cocycle
\[
\Check U\cdot\frac{\left(\Check\lambda^{2}-L_{4k} \right)}{8}
\]
combine to make 
\[
E\times\Delta^{k}\to S\times\Delta^{k}
\]
into a differential $\bsobwu$-oriented map, with respect to the
cocycle $\alpha$.   This defines 
\[
\mathcal
S_{\nu}^{2k}(E) \to
\fil_{0}\left(\msobwu_{i-4k};\alpha\right)^{S},
\]
and hence $\kappa_{E/S}$ as described above.

It is now a fairly routine exercise to verify that this functor
$\kappa_{E/S}$, and the resulting $q$ have the properties claimed by
Theorem~\ref{thm-main}.  Property~\thmItemref{item:2} is immediate
by assumption.  For property~\thmItemref{item:3}, write the
action of $\Z/2$ on $E\Z/2_{+}$ as $b\mapsto -b$.  Consider the
diagram of differential function complexes, in which the vertical maps
are obtained by smashing the identity map with the inclusions
$\{b\}\to E\Z/2$ and $\{-b\}\to E\Z/2$
\[
\xymatrix@=10pt{
&\fil_{0}\left(\msobwu_{i-4k};\alpha \right)^{S} \ar[d]^{b}\ar[dr] & \\
{\mathcal S}^{2k}_{\nu}(E) \ar[ur]\ar[r]\ar[dr] &
 \fil_{0}\left(\left(\msobwu\wedge E\Z/2_{+};\alpha \right)_{i-4k}
\right)^{S} \ar[r]^(.7){\kappa} & \fil_{0}\left(\tilde I_{i};\iotai
\right)^{S} \\ 
& \fil_{0}\left(\msobwu_{i-4k} \right)^{S} \ar[u]_{-b}\ar[ur] &.
}
\]

The map of fundamental groupoids induced by the upper composition is
$\kappa_{E/S}(\wustruct)$, by definition.  Since $\kappa$ factors
through the quotient by the diagonal $\Z/2$-action, the bottom
composition is $\kappa_{E/S}(-\wustruct)$.  The fact that the vertical
arrows are homotopy equivalences means that there is a homotopy
between the two ways of going around.  Choose one.  By definition this
homotopy is an isomorphism
\[
\tau(\wustruct):\kappa(-\wustruct)\xrightarrow{\approx}{}\kappa(\wustruct)
\]
in the fundamental groupoid 
\[
\pi_{\le 1}\left(\tilde I_{i},\iotai \right)^{S}.
\]
Any two homotopies extend over a disk, and so define the same isomorphism
in the fundamental groupoid.  The compositions of the
homotopies $\tau(\wustruct)$ and $\tau(-\wustruct)$ also extends over
the disk, and so is the identity map
\[
\tau(\wustruct)\circ\tau(-\wustruct)=\text{identity map of
$\kappa(\wustruct)$}.
\]

For the base change property~\thmItemref{item:4} note that if $E/S$ is
classified by a differential function 
\[
S\to \left(\msobwu_{i-4k};\alpha\right),
\]
then the map
\[
S'\to S\to \left(\msobwu_{i-4k};\alpha\right),
\]
is the map classifying a differential $\bsobwu$-orientation of
$E'/S'$.  The result then follows easily.

Now for the transitivity property~\thmItemref{item:5}.  Suppose that
$E/B$ has relative dimension $m$ and is classified by a differential
function
\begin{equation}\label{eq:t-73}
B\to \left(\msobwu,\alpha\right)_{-m},
\end{equation}
and that $B/S$ has relative dimension $\ell$.  Let
\begin{equation}\label{eq:t-70}
S^{0}\wedge \left(\msobwu_{-m} \right)_{+}
\end{equation}
denote the unreduced suspension spectrum of $\msobwu_{-m}$.
Interpreting~\eqref{eq:t-70} as the Thom spectrum
$\thom(\msobwu_{-m};0)$, we can regard the differential framing on
$B/S$ together with the map~\eqref{eq:t-73} as classified by a
differential function
\begin{equation}\label{eq:t-62}
S\to \left(S^{0}\wedge (\msobwu_{-m})_{+}, \alpha_{-m}
\right)_{-\ell}. 
\end{equation}
We will use the 
structure map
\begin{equation}\label{eq:t-71}
\left(S^{0}\wedge \left(\msobwu_{-m} \right)_{+};\alpha_{-m} \right)
\to \Sigma^{-m}\left(\msobwu;\alpha\right)
\end{equation}
of the differential spectrum $\left(\msobwu;\alpha \right)$.
Composing~\eqref{eq:t-62} with~\eqref{eq:t-71} gives a differential
function
\[
S\to \left(\msobwu;\alpha \right)_{-m-\ell}
\]
classifying a differential $\bsobwu$-orientation on $E/S$.  In
this way $E/S$ acquires an \orientation{}.  One easily checks that the 
differential cocycle refining the pullback of $\alpha$ is
\[
\frac18\check U\cup(\check \lambda^{2}-\check{L}_{4k})
\]
The transitivity isomorphism is then derived from
a choice of homotopy between the two ways of going around
the diagram of differential function spectra 
\[
\begin{CD}
\Sigma^{-\ell}\left(S^{0}\wedge (\msobwu_{-m})_{+}, \alpha_{-m}
\right) @>>> 
\Sigma^{-m-\ell}\left(\msobwu; \alpha\right) \\
@VVV @VVV \\
\Sigma^{-\ell}\left(S^{0}\wedge (\tilde I_{4k-m})_{+}, (\iotai)_{4k-m}
\right) @>>> 
\Sigma^{4k-m-\ell}\left(I, \iotai \right).
\end{CD}
\]
The clockwise composition associates to the differential function
classifying $E/B$ the value
\[
\kappa_{E/S}(\wustruct),
\]
and the counter-clockwise composition gives the value 
\[
\int_{B/S}\kappa_{E/B}(\wustruct).
\]

It remains to establish the properties of $q$ stated in
Corollary~\ref{thm:28}.
Properties~\thmItemref{item:7}--~\thmItemref{item:9} are formal
consequences of the corresponding properties of $\kappa$.
Property~\thmItemref{item:6} follows from Lemma~\ref{thm:25} using
methods similar to those we've been using for the properties of
$\kappa$.
\end{pf*}

\begin{rem}
The proofs of the symmetry, base change, and transitivity properties
didn't make use of the condition on $\kappa$ stated in
Proposition~\ref{thm:t-22}.  These properties are built into the
formalism of differential bordism theories, and would have held for
any $\kappa$.
\end{rem}

\subsection{The topological theory of quadratic functions}
\label{sec:topol-theory-quadr}

We now turn to the construction of a particular topological $\kappa$
\[
\kappa:\msobwu_{h\Z/2}\to \Sigma^{4k}\tilde I
\]
satisfying~\eqref{eq:t-27}.  To describe a map to $\tilde I$ requires a
slightly more elaborate algebraic object than an abelian group.  There
are several approaches, but for our purposes, the most useful involves
the language of {\em Picard categories} (see
Appendix~\ref{sec:picard-categ-anders}.)  Here we state the main
result without this language.

Let $\bar\nu_{2k}$ be the composite
\[
\bso\xrightarrow{\nu_{2k}}{}K(\Z/2,2k) \to
K(\qz(1),2k),
\]
and $\bsowubar$ its (equivariant) homotopy fiber.  The space $\bsowubar$
fits into a homotopy Cartesian square
\[
\begin{CD}
\bsowubar @>>> K(\qz(1),2k-1) \\
@VVV @VV\beta V\\
\bso @>>\nu_{2k}> K(\Z/2,2k).
\end{CD}
\]
The $\Z/2$-action corresponds to sending $\eta$ to $-\eta$, and the
associated Thom spectrum $\msowubar_{h\Z/2}$ is the bordism theory of
manifolds $N$ equipped with a homomorphism $t:\pi_{1}N\to\Z/2$
classifying a local system $\Z(1)$ and a cocycle $\eta\in
Z^{2k-1}(N;\qz(1))$ for which $\beta\eta\in
Z^{2k}(N;\Z(1))$ is an integral Wu-structure.

To ease the notation,
set 
\begin{gather*}
B=\bsobwu\times_{\Z/2}E\Z/2 \\
\bar{B}=\bsowubar\times_{\Z/2}E\Z/2.
\end{gather*}
Let $\mathcal C$ be the groupoid whose objects are closed $\bar B$-oriented
manifolds $(M,\eta)$ of dimension $(4k-1)$, and whose morphisms are
equivalence classes of $B$-oriented maps 
\begin{equation}\label{eq:t-39}
p:M\to \Delta^{1}
\end{equation}
equipped with $\bar B$-orientations $\eta_{0}$ and $\eta_{1}$ of
$\partial_{0}M=p^{-1}(0)$ and $\partial_{1}M=p^{-1}(1)$ which are
compatible with the $B$-orientation in the sense that
\[
\beta\eta_{i}=\lambda\vert_{\partial_{i}M},\quad i=0,1.
\]
The equivalence relation and the composition law are described in
terms of $B$-oriented maps to $f:E\to\Delta^{2}$ equipped with
compatible $\bar B$ orientations of $f^{-1}(e_{i})$, $i=0,1,2$.  Write
\[
\partial_{i}E=f\vert_{p^{-1}\partial_{i}\Delta^{2}}.
\]
Then for each such $E/\Delta^{2}$ we set
\begin{equation}\label{eq:t-40}
\partial_{0}E \circ
\partial_{2}E\sim
\partial_{1}E.
\end{equation}

\begin{rem}
Strictly speaking we define $\mathcal C$ to be the quotient of the
category freely generated by the maps~\eqref{eq:t-39}, by the
relations~\eqref{eq:t-40}.  This category happens to be a groupoid;
every morphism is represented by some map~\eqref{eq:t-39}, and the
composition law can be thought of as derived from the operation of
gluing together manifolds along common boundary components.  It works
out that the identity morphism $\text{Id}_{N}$ is represented by the
projection
\[
N\times \Delta^{1}\to \Delta^{1}.
\]
\end{rem}

\begin{prop}\label{thm:t-23}
The set of homotopy classes of maps
\begin{equation}\label{eq:t-43}
\msobwu_{h\Z/2}\to\Sigma^{4k}\tilde I
\end{equation}
can be identified with the set of equivalence classes of pairs of
invariants
\begin{align*}
\kappa_{4k} &:\left\{\text{morphisms of $\mathcal C$} \right\} \to \Q\\
\kappa_{4k-1} &:\left\{\text{objects of $\mathcal C$} \right\} \to \qz 
\end{align*}
satisfying
\begin{equation}\label{eq:t-4}
\begin{aligned}
\kappa_{4k}(M_{1}/\Delta^{1}\amalg M_{2}/\Delta^{1}) &=
\kappa_{4k}(M_{1}/\Delta^{1})+\kappa_{4k}(M_{2}/\Delta^{1}) \\
\kappa_{4k-1}(N_{1}\amalg N_{2}) &=
\kappa_{4k-1}(N_{1})+\kappa_{4k-1}(N_{2}),
\end{aligned}
\end{equation}
\begin{equation}\label{eq:t-41}
\kappa_{4k}(M/\Delta^{1}) \equiv
\kappa_{4k-1}(\partial_{0}M)-\kappa_{4k-1}(\partial_{1}M)\mod \Z,
\end{equation}
and for each $B$-oriented $f:E\to \Delta^{2}$ equipped with
compatible $\bar B$ orientations of $f^{-1}(e_{i})$, $i=0,1,2$, 
\begin{equation}\label{eq:t-42}
\kappa_{4k}\left(\partial_{0}E \right) -
\kappa_{4k}\left(\partial_{1}E \right) +
\kappa_{4k}\left(\partial_{2}E \right)=0.
\end{equation}
Two pairs $(\kappa_{4k},\kappa_{4k-1})$ and
$(\kappa'_{4k},\kappa'_{4k-1})$ are equivalent if 
there is a map 
\[
h:\left\{\text{objects of $\mathcal C$} \right\} \to\Q
\]
with
\begin{align*}
h\left(M_{1}\amalg M_{2} \right)
&=h\left(M_{1} \right)+h\left(M_{2} \right), \\
h(N) &\equiv \kappa'_{4k-1}(N)-\kappa_{4k-1}(N)\mod\Z
\end{align*}
and for each $M/\Delta^{1}$, 
\[
h\left(\partial_{0}M \right)- h\left(\partial_{1}M
\right)=\kappa_{4k}(M/\Delta^{1}) - \kappa_{4k}(M/\Delta^{1}).
\]
\end{prop}

\begin{rem}\label{rem:1}
In order to keep the statement simple, we have been deliberately
imprecise on one point in our statement of Proposition~\ref{thm:t-23}.
With our definition, the disjoint union of two $B$-oriented manifolds
doesn't have a canonical $B$-orientation.  In order to construct one,
we need to choose a pair of disjoint cubes embedded in $\R^{\infty}$.
The more precise statement is that for {\em any $B$-orientation
arising from any such choice of pair of cubes}, one has
\begin{align*}
\kappa_{4k}(M_{1}\amalg M_{2}) &=
\kappa_{4k}(M_{1})+\kappa_{4k}(M_{2}) \\
\kappa_{4k-1}(M_{1}\amalg M_{2}) &=
\kappa_{4k-1}(M_{1})+\kappa_{4k}(M_{2}).
\end{align*}
In particular, the value of $\kappa_{i}$ on a disjoint union is
assumed to be independent of this choice.  For more on this, see
Example~\ref{eg:1}, and the remark preceding it.
\end{rem}

\begin{rem}
Proposition~\ref{thm:t-23} can be made much more succinct.  Let
$\left(\Q\to\qz \right)$ denote the category with objects $\qz$ and in
which a map from $a$ to $b$ is a rational number $r$ satisfying
\[
r\equiv b-a\mod\Z.
\]
Then the assertion of Proposition~\ref{thm:t-23} is that the set of
homotopy classes~\eqref{eq:t-43} can be identified with the set of
``additive'' functors
\[
\kappa:\mathcal C\to \left(\Q,\qz \right)
\]
modulo the relation of ``additive'' natural equivalence.  The language
of Picard categories is needed in order to make precise this notion of
``additivity.''  We have chosen to spell out the statement of
Proposition~\ref{thm:t-23} in the way we have in order to make clear
the exact combination of invariants needed to construct the
map~\eqref{eq:t-43}.
\end{rem}

Aside from the fact that the integral Wu-structure on the boundary has
a special form, Proposition~\ref{thm:t-23} is a fairly straightforward
consequence of Corollary~\ref{thm:t-27}.  We now turn to showing
that this boundary condition has no real effect.

\begin{lem}\label{thm:t-18}
The square
\begin{equation}\label{eq:t-25}
\begin{CD}
\bsowubar @>>> \bso \\
@VVV @VVV \\
\bsobwu @>>> \bso\times K(\Q(1),2k),
\end{CD}
\end{equation}
is homotopy co-Cartesian.
The components of the bottom map are the defining projection to
$\bso$ and the image of $\lambda$ in 
\[
Z^{2k}(\bsobwu;\Q(1)).
\]
\end{lem}

\begin{pf}
It suffices to prove the result after profinite completion and after
localization at $\Q$.  Since 
\begin{align*}
K(\qz,2k-1) &\to K(\Z,2k) \quad\text{and}\\
\text{pt} &\to K(\Q,2k)
\end{align*}
are equivalences after profinite completion, so are the vertical maps
in~\eqref{eq:t-25} and so the square is homotopy co-Cartesian after
profinite completion.  The horizontal maps become equivalences after
$\Q$-localization, and so the square is also homotopy co-Cartesian
after $\Q$-localization.  This completes the proof.
\end{pf}

Passing to Thom spectra gives
\begin{cor}\label{thm:t-19}
The square
\[
\begin{CD}
\msowubar @>>> \mso \\
@VVV @VVV \\
\msobwu @>>> \mso\wedge\left(K(\Q(1),2k)_{+} \right)
\end{CD}
\]
is homotopy co-Cartesian.
\qed
\end{cor}
 
Finally, passing to homotopy orbit spectra, and using the fact that 
\[
\iota^{2}\times 1:K(\Q(1),2k)\times_{\Z/2}E\Z/2 \to
K(\Q,4k)\times B\Z/2
\]
is a stable weak equivalence (its cofiber has no homology) gives
\begin{cor}\label{thm:t-29}
The square
\[
\begin{CD}
\msowubar_{h\Z/2} @>>> \mso\wedge B\Z/2_{+} \\
@VVV @VVV \\
\msobwu_{h\Z/2} @>>> \mso\wedge \left(K(\Q,4k)\times B\Z/2\right)_{+}
\end{CD}
\]
is homotopy co-Cartesian.
\qed
\end{cor}

\begin{cor}\label{thm:t-20}
The map
\[
\msowubar_{h\Z/2}\to\msobwu_{h\Z/2}
\]
is $(4k-1)$-connected.  In particular, there are cell decompositions
of both spectra, for which the map of $(4k-1)$-skeleta
\begin{equation}\label{eq:t-34}
\left(\msowubar_{h\Z/2} \right)^{(4k-1)}
\to\left(\msobwu_{h\Z/2} \right)^{(4k-1)}
\end{equation}
is a weak equivalence.\qed
\end{cor}

\begin{pf*}{Proof of Proposition~\ref{thm:t-23}}
We will freely use the language of Picard categories, and the results
of Appendix~\ref{sec:picard-categ-anders}.

Choose a cell decomposition of $\bar B$, and take 
\[
\msowubarh^{(t)}=\thom\left(\bar B^{(t)};V \right).
\]
Let $\msobwuh^{(t)}$ be any cell decomposition
satisfying~\eqref{eq:t-34}.  Finally, let 
\[
\pi'_{\le1}
(\msobwuh)_{(1-4k)}
\]
be the Picard category whose objects are transverse
maps
\[
S^{4k-1}\to (\msobwuh^{(4k-1)} =
\msowubarh^{(4k-1)},
\]
and whose morphisms are homotopy classes, relative to
\[
S^{4k-1}\wedge \partial\Delta^{1}_{+}
\]
of transverse maps of pairs
\[
\left(S^{4k-1}\wedge \Delta^{1}_{+},S^{4k-1}\wedge
\partial\Delta^{1}_{+} \right) \\
\to \left(\msobwuh,\msobwuh^{(4k-1)} \right).
\]
By the cellular approximation theorem, and the geometric
interpretation of the homotopy groups of Thom spectra, the
Pontryagin-Thom construction gives an equivalence of Picard categories.  
\[
\pi'_{\le1}
(\msobwuh)_{(1-4k)}  \to
\pi_{\le 1}\sing\left(\left(\msobwuh\right)_{(1-4k)}\right).
\]
Now, by definition, a pair of
invariants $(\kappa_{4k},\kappa_{4k-1})$ is a functor (of Picard
categories)
\begin{equation}\label{eq:t-35}
\pi'_{\le1}(\msobwuh)_{(1-4k)}\to \left(\Q\to\qz \right).
\end{equation}
It is easy to check, using the exact sequences
\[
\ext\left(\pi_{4k-1}\msobwuh,\Z \right) \\ \rightarrowtail A
\twoheadrightarrow \hom\left(\pi_{4k}\msobwuh,\Z \right),
\]
(where $A$ either the group of homotopy classes of
maps~\eqref{eq:t-43}, or the natural equivalences classes of
functors~\eqref{eq:t-35}) that this gives an isomorphism of the group
of equivalence classes of pairs $(\kappa_{4k},\kappa_{4k-1})$ and the
group of natural equivalence classes of functors of Picard
categories~\eqref{eq:t-35}.  The result now follows from
Corollary~\ref{thm:t-27}.
\end{pf*}

\subsection{The topological $\kappa$}
\label{sec:topological-kappa}
By Proposition~\ref{thm:t-23}, in order to define
the map
\[
\kappa:\msobwu_{h\Z/2}\to\Sigma^{4k}\tilde I.
\]
we need to construct invariants $(\kappa_{4k}, \kappa_{4k-1})$,
satisfying~\eqref{eq:t-4}--\eqref{eq:t-42}.  In the case $\lambda=0$,
these invariants appear in surgery theory in the work of
Milgram~\cite{milgram74:_surger}, and Morgan--Sullivan~\cite{morgan74}.
The methods described in~\cite{milgram74:_surger,morgan74} can easily
be adapted to deal with general $\lambda$.  We will adopt a more
homotopy theoretic formulation, which is more convenient for our
purposes.

Suppose that $M/\Delta^{1}$ is a morphism of $\mathcal C$, and let
$\sigma$ be the signature of the non-degenerate bilinear form
$\int_{M}x\cup y$ on the subgroup of $H^{2k}(M;\Q(1))$ consisting of
elements which vanish on the boundary.  Because of the boundary
condition the image of $\lambda$ in $H^{2k}(M;\Q(1))$ is in this
subgroup, and so
\[
\int_{M}\lambda^{2}
\]
is also defined.  We set
\begin{equation}\label{eq:t-36}
\kappa_{4k}(M)=\frac18\left(\int_{M}\lambda^{2}-\sigma \right).
\end{equation}

We will show in Proposition~\ref{thm:t-4} of
\S\ref{sec:quadratic-functions} that the ``Wu-structure'' 
\[
\eta\in
Z^{2k-1}(N;\qz(1))
\]
on an object $N\in \mathcal C$ gives rise to a
quadratic refinement
\[
\phi=\phi_{N}=\phi_{N,\eta}:H^{2k}(N;\Z(1))_{\text{tor}}\to\qz
\]
of the link pairing
\[
H^{2k}(N;\Z(1))_{\text{tor}}\times
H^{2k}(N;\Z(1))_{\text{tor}}\to \qz.
\]
The function $\phi$ is a cobordism invariant in the sense that if
$M/\Delta^{1}$ is a morphism in $\mathcal C$, and $x$ is an element of
$H^{2k}(M;\Z(1))$ whose restriction to $\partial M$ is torsion, then
\begin{equation}\label{eq:t-44}
\begin{gathered}
\phi_{\partial_{0}M}(\partial_{0}x)
-\phi_{\partial_{1}M}(\partial_{1}x)
\equiv \frac12\int_{M} (x^{2}-x\lambda)\mod \Z \\
\partial_{i}x=x\vert_{\partial_{i}M}.
\end{gathered}
\end{equation}
Note that the integral is well-defined since it depends only on the
image of $x$ in $H^{2k}(M;\Q(1))$, which vanishes on $\partial M$.
Note also that the integral vanishes when $x$ is itself torsion.

We define $\kappa_{4k-1}(N)$ by
\begin{equation}\label{eq:t-37}
\begin{gathered}
e^{\tpi\,\kappa_{4k-1}(N)} = \frac1{\sqrt{d}}\sum_{x\in
H^{2k}(N;\Z(1))_{\text{tor}}} e^{-\tpi \phi(x)} \\
d=\# H^{2k}(N;\Z(1))_{\text{tor}}. 
\end{gathered}
\end{equation}

With these definitions property~\eqref{eq:t-4} is immediate, and
property~\eqref{eq:t-42} is nearly so.  In the situation
of~\eqref{eq:t-42}, note that $\partial E$ is homeomorphic to a smooth
manifold whose signature is zero since it bounds an
oriented manifold.  Novikov's
additivity formula for the signature then gives
\[
\sigma\left(\partial_{0}E \right) -
\sigma\left(\partial_{1}E \right) +
\sigma\left(\partial_{2}E \right)=
\sigma\left(\partial E \right)=0
\]
which is property~\eqref{eq:t-42}.

For~\eqref{eq:t-41} we need an algebraic result of
Milgram~\cite{milgram74:_surger,morgan74,milnor73:_symmet}.  Suppose
that $V$ is a vector space over $\Q$ of finite dimension, and
$q_{V}:V\to\Q$ is a quadratic function (not necessarily even) whose
underlying bilinear form
\[
B(x,y)=q_{V}(x+y)-q_{V}(x)-q_{V}(y) 
\]
is non-degenerate.  One easily checks that $q_{V}(x)-q_{V}(-x)$ is
linear, and so there exists a unique $\lambda\in V$ with
\[
q_{V}(x)-q_{V}(-x)=-B(\lambda,x).
\]
By definition
\[
q_{V}(x)+q_{V}(-x) = B(x,x).
\]
Adding these we have
\begin{equation}\label{eq:t-6}
q_{V}(x)  = \frac{B(x,x)-B(x,\lambda)}{2}.
\end{equation}
The class $\lambda$ is a {\em characteristic element} of $B(x,y)$.

Suppose that $L$ is a lattice in $V$ on which $q_{V}$ takes integer
values, and let 
\[
L^{\ast}=\left\{x\in V\mid B(x,y)\in\Z\quad\forall y\in L \right\}
\]
be the dual lattice.  Then $L\subset L^{\ast}$, and $q_{V}$ descends
to a non-degenerate quadratic function
\begin{align*}
q : L^{\ast}/L &\to\qz \\
x &\mapsto q_{V}(x) \mod\Z.
\end{align*}

We associate to $q$ the Gauss sum
\begin{gather*}
g(q)=\frac1{\sqrt{d}}\sum_{x\in L^{\ast}/L}e^{-\tpi\,q(x)} \\
d=\left| L^{\ast}/L \right|
\end{gather*}
Finally, let $\sigma=\sigma(q_{V})$ denote the signature of $B$.

\begin{prop}\label{thm:t-5}
With the above notation, 
\begin{equation}\label{eq:t-5}
g(q)=e^{\tpi (B(\lambda,\lambda)-\sigma)/8}.
\end{equation}
\end{prop}

\begin{pf}
The result is proved in Appendix 4 of~\cite{milnor73:_symmet} in case
$\lambda=0$.  To reduce to this case, first note that both the left
and right hand sides of~\eqref{eq:t-5} are independent of the choice of
$L$.  Indeed, the right hand side doesn't involve the choice of $L$,
and the fact that the left hand side is independent of $L$ is Lemma~1,
Appendix~4 of~\cite{milnor73:_symmet} (the proof doesn't make use of
the assumption that $q_{V}$ is even).  Replacing $L$ by $4L$, if
necessary, we may assume that
\[
B(x,x)\equiv 0\mod 2\qquad x\in L,
\]
or, equivalently that
\[
\eta=\frac12\lambda
\]
is in $L^{\ast}$.  Set 
\[
q'_{V}(x)=\frac12 B(x,x).
\]
The function $q_{V}$ and $q'_{V}$ are quadratic refinements of the same
bilinear form, $q'_{V}$ takes integer values on $L$, and
\[
q_{V}(x)=q'_{V}(x-\eta)-q'_{V}(\eta).
\]
The case $\lambda=0$ of~\eqref{eq:t-5} applies to $q'_{V}$ giving
\[
g(q')=e^{\tpi \sigma/8}.
\]
But 
\begin{align*}
g(q')&=\frac1{\sqrt{d}}\sum_{x\in L^{\ast}/L}e^{-\tpi\,q'(x)} \\
&=\frac1{\sqrt{d}}\sum_{x\in L^{\ast}/L}e^{-\tpi\,q'(x-\eta)} \\
&=\frac1{\sqrt{d}}\sum_{x\in L^{\ast}/L}e^{-\tpi\,q(x)-q'(\eta)},
\end{align*}
and the result follows since 
\[
q'(\eta)=\frac{B(\lambda,\lambda)}{8}.
\]
\end{pf}

\begin{prop}\label{thm:t-16}
The invariants $(\kappa_{4k},\kappa_{4k-1})$ defined above
satisfy the conditions of Proposition~\ref{thm:t-23}:
\begin{equation}\label{eq:t-23}
\begin{aligned}
\kappa_{4k}(M_{1}/\Delta^{1}\amalg M_{2}/\Delta^{1}) &=
\kappa_{4k}(M_{1}/\Delta^{1})+\kappa_{4k}(M_{2}/\Delta^{1}) \\
\kappa_{4k-1}(N_{1}\amalg N_{2}) &=
\kappa_{4k-1}(N_{1})+\kappa_{4k-1}(N_{2}),
\end{aligned}
\end{equation}
\begin{equation}\label{eq:t-24}
\kappa_{4k}(M/\Delta^{1}) \equiv
\kappa_{4k-1}(\partial_{0}M)-\kappa_{4k-1}(\partial_{1}M)\mod \Z,
\end{equation}
and for each $B$-oriented $f:E\to \Delta^{2}$ equipped with
compatible $\bar B$ orientations of $f^{-1}(e_{i})$, $i=0,1,2$, 
\begin{equation}\label{eq:t-38}
\kappa_{4k}\left(\partial_{0}E \right) -
\kappa_{4k}\left(\partial_{1}E \right) +
\kappa_{4k}\left(\partial_{2}E \right)=0.
\end{equation}
\end{prop}

\begin{pf}
As we remarked just after defining $\kappa_{4k-1}$
(equation~\eqref{eq:t-37}), property~\eqref{eq:t-23} is immediate,
and~\eqref{eq:t-38} follows from Novikov's additivity formula for the
signature.  For property~\eqref{eq:t-24}, let $M'/\Delta^{1}$ be a
morphism in $\mathcal C$, and $M/\Delta^{1}$ any morphism for which
\[
M/\Delta^{1} + \text{Id}_{\partial_{1}M'}\sim M'/\Delta^{1}.
\]
Note that this identity forces $\partial_{1}M$ to be empty.
Applying~\eqref{eq:t-38} to $N\times \Delta^{2}\to \Delta^{2}$ one
easily checks that both the left and right sides of~\eqref{eq:t-24}
vanish for $\text{Id}_{\partial_{1}M'}$.  It then follows
from~\eqref{eq:t-23} that to check~\eqref{eq:t-24} for general
$M'/\Delta^{1}$, it suffices to check~\eqref{eq:t-24} for morphisms
$M/\Delta^{1}$ with the property that
\[
\partial_{1}M=\phi.
\]
We apply Milgram's result to the
situation
\begin{align*}
V&=\text{image}\left(H^{2k}(M,\partial M;\Q(1))\to H^{2k}(M;\Q(1)) \right) \\
&= \ker\left(H^{2k}(M;\Q(1))\to H^{2k}(\partial M;\Q(1)) \right),\\
L&=\text{image}\left(H^{2k}(M,\partial M;\Z(1))\to V\right),
\end{align*}
with $q_{V}$ the quadratic function 
\[
\phi(x)=\phi_{M}(x)\int_{M}\frac{x^{2}-x\lambda}{2}.
\]
Poincar\'e duality gives
\[
L^{\ast}=V\cap\text{image}\left(H^{2k}(M;\Z(1))\to H^{2k}(M;\Q(1))\right).
\]
The right hand side of~\eqref{eq:t-5} is by definition 
\[
e^{\tpi\,\kappa_{4k}(M)}.
\]
We need to identify the left hand side with 
\[
\kappa_{4k-1}\left(\partial_{0}M \right) -
\kappa_{4k-1}\left(\partial_{1}M \right)=
\kappa_{4k-1}\left(\partial_{0}M \right).
\]
Write 
\begin{align*}
A &=H^{2k}(\partial M;\Z(1))_{\text{tor}},\\
A_{0} &=\text{image}\left(H^{2k}(M)_{\text{tor}}\to 
H^{2k}(\partial M)_{\text{tor}} \right).
\end{align*}
and let $B$ be the torsion
subgroup of the image of 
\[
H^{2k}(M)\to H^{2k}(\partial M).
\]
Then
\[
L^{\ast}/L = B/A_{0}.
\]
By cobordism invariance~\eqref{eq:t-44} the restriction of 
$\phi_{\partial M}$ to $B$ is compatible with this isomorphism,
and so
\[
\frac1{\sqrt{\left|L^{\ast}/L \right|}}
\sum_{x\in L^{\ast}/L}e^{-\tpi\,q(x)}=
\frac1{\sqrt{\left|B/A_{0} \right|}}
\sum_{x\in B/A_{0}}e^{-\tpi\,\phi_{\partial M}(x)}
\]
By Lemma~\ref{thm:t-28} below, the subgroup $B$ coincides with the
annihilator $A_{0}^{\ast}$ of $A_{0}$, and so we can further re-write
this expression as
\[
\frac1{\sqrt{\left|A_{0}^{\ast}/A_{0} \right|}}
\sum_{x\in A_{0}^{\ast}/A_{0}}e^{-\tpi\,\phi_{\partial M}(x)}.
\]
The identification of this with 
\[
\kappa_{4k-1}(\partial M)=\kappa_{4k-1}(\partial_{0}M)
\]
is then given by Lemma~\ref{thm:t-7} below.
\end{pf}

\begin{lem}\label{thm:t-7}
Suppose that $A$ is a finite abelian group, and 
\[
q:A\to \qz
\]
a (non-degenerate) quadratic function with underlying bilinear form
$B$.  Given a subgroup $A_{0}\subset A$ on which $q$ vanishes, let
$A_{0}^{\ast}\subset A$ be the dual of $A_{0}$:
\[
A_{0}^{\ast}=\left\{x\in A\mid B(x,a)=0, a\in A_{0} \right\}.
\]
Then $q$ descends to a quadratic function on $A_{0}^{\ast}/A_{0}$, and
\[
\frac{1}{\sqrt{\left|A \right|}}\sum_{x\in A}e^{-\tpi\,q(a)}=
\frac{1}{\sqrt{\left|A'/A_{0} \right|}}\sum_{x\in
A_{0}^{\ast}/A_{0}}e^{-\tpi\,q(a)}.
\]
\end{lem}

\begin{pf}
Since $q(a)=0$ for $a\in A_{0}$, we have 
\[
q(a+x)=q(a)+B(x,a).  
\]
Choose coset representatives $S\subset A$ for $A/A_{0}$, and write
$S_{0}=S\cap A_{0}^{ast}$.  The set $S_{0}$ is a set of coset
representatives for $A_{0}^{\ast}/A_{0}$.  Now write
\begin{multline*}
\frac{1}{\sqrt{\left|A \right|}}\sum_{x\in A}e^{-\tpi\,q(a)}
=
\frac{1}{\sqrt{\left|A \right|}}
\sum_{x\in S}\sum_{a\in A_{0}}e^{-\tpi\,q(x+a)} \\
= \frac{1}{\sqrt{\left|A \right|}} \sum_{x\in
S}\left(e^{-\tpi\,q(x)}\sum_{a\in A_{0}}e^{-\tpi\,B(x,a)} \right).
\end{multline*}
By the linear independence of characters
\[
\sum_{a\in A_{0}}e^{-\tpi\,B(x,a)} =
\begin{cases}
0 &\quad B(x,a)\ne 0 \\
\left|A_{0}\right| &\quad B(x,a) = 0.
\end{cases}
\]
It follows that 
\[
\frac{1}{\sqrt{\left|A \right|}} \sum_{x\in
S}\left(e^{-\tpi\,q(x)}\sum_{a\in A_{0}}e^{-\tpi\,B(x,a)} \right)
\\ =
\frac{\left|A_{0}\right|}{\sqrt{\left|A \right|}} \sum_{x\in
S_{0}}\left(e^{-\tpi\,q(x)} \right).
\]
This proves the result, since the bilinear form $B$ identifies
$A/A_{0}^{\ast}$ with the character group of $A_{0}$, giving
\[
\left|A/A_{0}^{\ast} \right|
=
\left|A_{0}\right|
\]
and
\[
\left|A \right|=\left|A_{0} \right|\cdot \left|
A_{0}^{\ast}/A_{0}\right|\cdot
\left|A/A_{0}^{\ast} \right|=
\left|A_{0} \right|^{2}\left|A_{0}^{\ast}/A_{0} \right|.
\]
\end{pf}

We have also used
\begin{lem}\label{thm:t-28}
With the notation Lemma~\ref{thm:t-7}, the subgroup $A_{0}^{\ast}$
coincides with the torsion subgroup of the image of
\[
H^{2k}(M)\to H^{2k}(\partial M).
\]
\end{lem}

\begin{pf}
Poincar\'e duality identifies the $\qz$ dual of 
\begin{equation}\label{eq:t-45}
H^{2k}(\partial M)\leftarrow H^{2k}(M)
\end{equation}
with
\begin{equation}\label{eq:t-46}
H^{2k-1}(\partial M;\qz)\to 
H^{2k}(M,\partial M;\qz).
\end{equation}
Thus the orthogonal complement of the image of~\eqref{eq:t-45} 
is the image of~\eqref{eq:t-46}.  The claim follows easily.
\end{pf}

\subsection{The quadratic functions}\label{sec:quadratic-functions} 

We now turn to the relationship between integral Wu-structures and
quadratic functions.  That there is a relationship at all has a simple
algebraic explanation.  Suppose that $L$ is a finitely generated free
abelian group equipped with a non-degenerate symmetric bilinear form
\[
B:L\times L\to\Z.
\]
A {\em characteristic element} of $B$ is an element $\lambda\in L$
with the property
\[
B(x,x)\equiv B(x,\lambda)\mod 2.
\]
If $\lambda$ is a characteristic element, then
\begin{equation}\label{eq:t-7}
q(x)=\frac{B(x,x)-B(x,\lambda)}{2}
\end{equation}
is a quadratic refinement of $B$.  Conversely, if 
$q$ is a quadratic refinement of $B$,
\[
q(x+y)-q(x)-q(y)+q(0)=B(x,y),
\] 
then 
\[
q(x)-q(-x) :L\to \Z
\]
is linear, and so there exists $\lambda\in L$ with 
\[
q(x)-q(-x)= -B(x,\lambda).
\]
We also have
\[
q(x)+q(-x)=2\,q(0)+B(x,x).
\]
If we assume in addition that $q(0)=0$ then it follows that $q(x)$ is
given by~\eqref{eq:t-7}.   Thus the set of quadratic refinements $q$ of
$B$, with $q(0)=0$ are in one to one correspondence with the
characteristic elements $\lambda$ of $B$.  

We apply the above discussion to the situation in which
$M$ is an oriented manifold of dimension $4k$,
\[
L=H^{2k}(M;\Z)/\text{torsion},
\]
and 
\[
B(x,y)=\int_M x\cup y.
\]
By definition, the Wu-class 
\[
\nu_{2k}\in H^{2k}(M;\Z/2)
\]
satisfies
\[
\int_{M}x^{2} = \int_{M}x\,\nu_{2k}\in \Z/2, \quad x\in H^{2k}(M;\Z/2).
\]
Thus the integer lifts $\lambda$ are exactly the
characteristic elements of $B$, and correspond to quadratic
refinements $q^{\lambda}$ of the intersection pairing.  

\begin{lem}
The function 
\[
q^{\lambda}(x)=\tfrac12\int_{M} (x^{2}-x\lambda)
\]
defines a homomorphism
\[
\pi_{4k}\msobwu\wedge K(\Z,2k)\to\Z.
\]
\end{lem}

\begin{pf}
The group $\pi_{4k}\msobwu\wedge K(\Z,2k)_{+}$ is the cobordism group
of triples $(M,x,\lambda)$ with $M$ an oriented manifold of dimension
$4k$, $\lambda\in Z^{2k}(M;\Z)$ a lift of $\nu_{2k}$ and $x\in
Z^{2k}(M;\Z)$.  The group $\pi_{4k}\msobwu\wedge K(\Z,2k)$ is the
quotient by the subgroup in which $x=0$.  
Let $(M_{1},\lambda_{1},x_{1})$ and
$(M_{2},\lambda_{2},x_{2})$ be two such manifolds.  Since
\[
x_{1}\cup x_{2}=0\in H^{4k}(M_{1}\amalg M_{2})
\]
we have
\[
q^{\lambda_{1}+\lambda_{2}}(x_{1}+ x_{2})=
q^{\lambda_{1}}(x_{1})+
q^{\lambda_{2}}(x_{2}),
\]
and so $q^{\lambda}$ is additive.  If $M=\partial N$ and both
$\lambda$ and $x$ extend to $N$, then
\[
q^{\lambda}(x)=\frac12\int_{M}x^{2}-x\lambda
=\frac12\int_{N}d(x^{2}-x\lambda)=0,
\]
and $q$ is a cobordism invariant.   The result now follows since $q(0)=0$.
\end{pf}

We will now show that the homomorphism 
\begin{align*}
\pi_{4k}\msobwu\wedge K(\Z,2k) &\to \Z \\
(M,\lambda,x) &\mapsto \tfrac12\int_{M}x^{2}-x\lambda
\end{align*}
has a canonical lift to a map 
\begin{equation}\label{eq:t-47}
\msobwu\wedge K(\Z,2k)\to \Sigma^{4k}\tilde I.
\end{equation}
It suffices to produce the adjoint to~\eqref{eq:t-47} which is a map
\[
\msobwu \to \Sigma^{4k}\tilde I\left(K(\Z,2k) \right).
\]
Since $\msobwu$ is $(-1)$-connected this will have to factor 
through the $(-1)$-connected cover%
\footnote{%
The symbol $X\langle n,\dots,m \rangle$ indicates the Postnikov
section of $X$ having homotopy groups only in dimension $n\le i\le
m$.  The angled brackets are given lower precedence than suspensions, so
that the notation $\Sigma X\langle n,\dots,m \rangle$ means
$\left(\Sigma X \right)\langle n,\dots,m \rangle$ and coincides with 
$\Sigma\left(X \langle n+1,\dots,m+1 \rangle \right)$.
}%
\[
\msobwu \to \Sigma^{4k}\tilde I\left(K(\Z,2k)
\right)\langle0,\dots,\infty \rangle.
\]

To work out the homotopy type of $\Sigma^{4k}\tilde I\left(K(\Z,2k)
\right)\langle 0,\dots,\infty \rangle$ we will follow the approach of
Browder and Brown~\cite{browder69:_kervair,brown65:_note_kervair,
brown72:kervair,brown71:_kervair}.

For a space $X$ and a cocycle $x\in Z^{2k}(X;\Z/2)$ the theory of
Steenrod operations provides a universal $(2k-1)$-cochain $h$ with the
property
\[
\delta h=x\cup x + \Sq^{2k}(x).
\]
Taking the universal case $X=K(\Z,2k)$ $x=\iota$, this can be
interpreted as as an explicit homotopy making the diagram 
\[
\begin{CD}
\Sigma^{\infty} K(\Z,2k) @>\iota\cup \iota>>
\Sigma^{4k}\hz \\
@V\iota VV @VVV \\
\Sigma^{2k}\hz @>> \Sq^{2k}> \Sigma^{4k}\hz/2.
\end{CD}
\]
Passing to Postnikov sections gives an explicit homotopy making
\begin{equation}\label{eq:t-50}
\begin{CD}
\Sigma^{\infty} K(\Z,2k)\langle-\infty, 4k \rangle @>\iota\cup \iota>>
\Sigma^{4k}\hz \\
@V\iota VV @VVV \\
\Sigma^{2k}\hz @>> \Sq^{2k}> \Sigma^{4k}\hz/2.
\end{CD}
\end{equation}
commute.

\begin{prop}\label{thm:t-15}
The square~\eqref{eq:t-50} is homotopy Cartesian.
\end{prop}

\begin{pf}
This is an easy consequence of the Cartan--Serre computation of the
cohomology of $K(\Z,2k)$.
\end{pf}

It will be useful to re-write~\eqref{eq:t-50}.  Factor
the bottom map as 
\[
\Sigma^{2k}\hz\to \Sigma^{2k}\hz/2\to \Sigma^{4k}\hz/2,
\]
and define a spectrum $X$ by requiring that the squares in
\[
\begin{CD}
\Sigma^{\infty} K(\Z,2k)\langle-\infty, 4k \rangle @>>> X
@>>>\Sigma^{4k}\hz \\
@V\iota VV @VVV @VVV \\
\Sigma^{2k}\hz @>>> \Sigma^{2k}\hz/2@>> \Sq^{2k}> \Sigma^{4k}\hz/2.
\end{CD}
\]
be homotopy Cartesian.  Then~\eqref{eq:t-50} determines a homotopy
Cartesian square
\[
\begin{CD}
\Sigma^{2k-1}\hz/2 @>\beta >> \Sigma^{2k}\hz \\
@V\beta \Sq^{2k}VV @VVV \\
\Sigma^{4k}\hz @>>> \Sigma^{\infty}K(\Z,2k)\langle-\infty,\dotsc ,0 \rangle,
\end{CD}
\]
in which the horizontal map is the inclusion of the fiber of $\iota$,
and the vertical map is the inclusion of the fiber of the map to $X$.
Taking Anderson duals then gives
\begin{prop}\label{thm:t-3}
The diagram~\eqref{eq:t-50} determines a homotopy Cartesian square
\begin{equation}\label{eq:t-51}
\begin{CD}
\Sigma^{4k}\tilde IK(\Z,2k)\langle 0,\dots ,\infty \rangle @> b >> \Sigma^{2k}H\Z \\
@VaVV @VVV \\
H\Z @> \chi\Sq^{2k}>>\Sigma^{2k}H\Z/2.
\end{CD}
\end{equation}
\qed
\end{prop}

By Proposition~\ref{thm:t-3}, to give a map
\[
X\to \Sigma^{4k}\tilde IK(\Z,2k)\langle 0,\dots ,\infty \rangle
\]
is to give cocycles $a\in Z^{0}(X;\Z)$, $b\in Z^{2k}(X;\Z)$, and a
cochain $c\in C^{2k-1}(X;\Z/2)$ satisfying
\[
\chi Sq^{2k}(a)-b\equiv\delta c\mod 2.
\]
In case $X=\msobwu$ and $a=U$ is the Thom cocycle, the theory of Wu
classes gives a universal $c$ for which 
\[
\delta c=\chi\Sq^{2k}a+\nu_{2k}\cdot U.
\]
Taking $b=-\lambda\cdot U$ leads to a canonical map
\begin{equation}\label{eq:t-49}
\msobwu\to X\to \Sigma^{4k}\tilde IK(\Z,2k).
\end{equation}

\begin{prop}
The adjoint to~\eqref{eq:t-49}
\begin{equation}\label{eq:t-53}
\msobwu\wedge K(\Z,2k)\to \Sigma^{4k}\tilde I
\end{equation}
is a lift of the homomorphism 
\[
\tfrac12\int_{M}\left(x^{2}-x\lambda \right).
\]
\end{prop}

\begin{pf}
We'll use the notation of~\eqref{eq:t-50} and~\eqref{eq:t-51}.   The maps
\begin{gather*}
\Sigma^{\infty}K(\Z,2k)\langle-\infty,\dotsc ,0 \rangle\xrightarrow{(\iota,\iota^{2})}{} 
\Sigma^{2k}\hz\vee \Sigma^{4k}\hz \\
\Sigma^{4k}\tilde I\left(K(\Z,2k) \right)\langle0,\dots,\infty \rangle
\xrightarrow{(a,b)}{}\hz\vee \Sigma^{2k}\hz,
\end{gather*}
are rational equivalences, and by construction
\begin{align*}
\iota\circ \Sigma^{4k}\tilde I(b)&=2 &
\iota\circ \Sigma^{4k}\tilde I(a)&=0 \\
\iota\cup\iota\circ \Sigma^{4k}\tilde I(b)&=0 &
\iota\cup\iota\circ \Sigma^{4k}\tilde I(a)&=2.
\end{align*}
It follows that the rational evaluation map
\[
\Sigma^{4k}\tilde I\left(K(\Z,2k) \right)\langle0,\dots,\infty
\rangle 
\wedge
\Sigma^{\infty}K(\Z,2k)\langle-\infty,\dotsc ,0 \rangle \\\to
\Sigma^{4k}\tilde I\to \hq
\]
is $(a\iota^{2}+b\iota)/2$.
This means that
\[
\msobwu\wedge K(\Z,2k)\to \Sigma^{4k}\tilde IK(\Z,2k)\wedge
K(\Z,2k) \\ \to \Sigma^{4k}\tilde I\to \Sigma^{4k}\hq
\]
is given by
\[
(U\cdot x^{2} - U\cdot\lambda\cdot x)/2,
\]
and the claim follows.
\end{pf}

\begin{rem}

The topological analogue of a characteristic element is an integral
Wu-structure, and the object in topology corresponding to a integer
invariant is a map to $\tilde I$.  We've shown that an integral
Wu-structure gives a function to $\tilde I$.  That this function is
quadratic is expressed in topology by the following diagram which is
easily checked to be homotopy commutative (and in fact to come with an
explicit homotopy):
\begin{equation}\label{eq:t-52}
\begin{CD}
\msobwu\wedge K(\Z,2k)\wedge K(\Z,2k) @> \int x\cup y >>
\Sigma^{4k}\hz \\
@V(x+y)-(x)-(y) VV @VVV \\
\msobwu\wedge K(\Z,2k) @>>>
\Sigma^{4k}\tilde I.
\end{CD}
\end{equation}
\end{rem}

There is more information in the map~\eqref{eq:t-53}.
Consider the following diagram in which the bottom vertical arrows are
localization at $\hq$ and the vertical sequences are fibrations
\begin{equation}
\begin{CD}
\msobwu\wedge K(\qz,2k-1)  @>>> \Sigma^{4k-1}I \\
@VVV @VVV \\
\msobwu\wedge K(\Z,2k)  @>>> \Sigma^{4k}\tilde I \\
@VVV @VVV \\
\msobwu\wedge K(\Q,2k)  @>>> \Sigma^{4k}\hq.
\end{CD}
\end{equation}
The top map is classified by a homomorphism
\[
\pi_{4k-1}\msobwu\wedge K(\qz,2k-1)\to\qz.
\]
In geometric terms, this associates to each integral Wu-structure
on an oriented manifold of $N$ of dimension $(4k-1)$ a function
\[
\phi_{4k-1}=\phi_{4k-1}(N,\slot):H^{2k-1}(N;\qz)\to\qz
\]
One easily checks (using~\eqref{eq:t-52}) that this map is a quadratic
refinement of the link pairing
\[
\phi_{4k-1}(x+y)-\phi_{4k-1}(x)-\phi_{4k-1}(y)=\int_{N}x\cup\beta y.
\]
The bottom map is classified by a homomorphism 
\[
\phi_{4k}:\pi_{4k}\msobwu\wedge K(\Q,2k)\to \Q.
\]
Thinking of $\pi_{4k}\msobwu\wedge K(\Q,2k)$ as the relative homotopy
group
\[
\pi_{4k}\left(\msobwu\wedge K(\Z,2k),\msobwu\wedge
K(\qz,2k-1) \right)
\]
we interpret it in geometric terms as the bordism group of manifolds
with boundary $M$, together with cocycles
\begin{align*}
x &\in Z^{2k}(M;\Z)\\
y &\in Z^{2k-1}(M;\qz)
\end{align*}
satisfying:
\[
\beta(y)=x\vert_{\partial M}.
\]
Because of this identity, the cocycle $x$ defines a class
$x\in H^{2k}(M,\partial M;\Q)$ and in fact
\[
\phi_{4k}(M,x,y)=\phi_{4k}(M,x)=\int_{M}(x^{2}-x\lambda)/2
\]
The compatibility of the maps $\phi$ with the connecting homomorphism 
gives
\[
\phi_{4k-1}(\partial M,y)\equiv \phi_{4k}(M,x)\mod\Z.
\]

The functions $\phi_{4k}$ and $\phi_{4k-1}$ share the defect that they
are not quite defined on the correct groups.  For instance $\phi_{4k}$
is defined on the middle group in the exact sequence
\[
\cdots \to H^{2k-1}(\partial M;\Q)\xrightarrow{\delta}{}
H^{2k}(M,\partial M;\Q)\to H^{2k}(M;\Q)\to\cdots
\]
but does not necessarily factor through the image of this group in
\[
H^{2k}(M;\Q).
\]
The obstruction is
\[
\phi_{4k}(M,\delta x)=\frac12\int_{\partial M}x\cdot \lambda.
\]
Similarly, $\phi_{4k-1}$ is defined on $H^{2k-1}(N;\qz)$, but does not
necessarily factor through $H^{2k}(N;\Z)_{\text{torsion}}$, and the
obstruction is 
\[
\phi_{4k-1}(N,x)=\frac12\int_{N}x\cdot \lambda,\quad x\in H^{2k-1}(N;\Q).
\]
Both of these obstructions vanish if it happens that the integral
Wu-structures in dimension $(4k-1)$ are torsion.  This is the
situation considered in \S\ref{sec:topol-theory-quadr}.

We still haven't quite constructed functions $\phi_{N,\eta}$ needed
in \S\ref{sec:topological-kappa}.  For one thing, we've only described
the compatibility relation in the case of a manifold with boundary.
To get a formula like~\eqref{eq:t-44} one simply considers maps of pairs
\begin{multline*}
\left(S^{4k-1}\wedge \Delta^{1}_{+},S^{4k-1}\wedge
\partial\Delta^{1}_{+} \right) \\
\to \left(\msobwu\wedge K(\Z,2k),\msobwu\wedge
K(\qz,2k-1) \right)
\end{multline*}
instead of relative homotopy groups.  We also haven't incorporated the
symmetry.  We now indicate the necessary modifications.

We need to work with spaces and spectra%
\footnote{%
In the terminology of~\cite{LMayS} these are {\em naive}
$\Z/2$-spectrum.  }%
\ equipped with an action of the group $\Z/2$, with the convention
that an equivariant map is regarded as a a weak equivalence if the
underlying map of spaces or spectra is a weak
equivalence.\footnote{This is sometimes known as the ``coarse'' model
category structure on equivariant spaces or spectra.  The alternative
is to demand that the map of fixed points be a weak equivalence as
well.}  With this convention the map
\[
X\times E\Z/2\to E\Z/2
\]
is a weak equivalence, and a spectrum $E$ with the
trivial $\Z/2$-action represents the equivariant cohomology theory
\[
X\mapsto E_{\Z/2}^{\ast}(X)=E^{\ast}(E\Z/2\times_{\Z/2}X).
\]
Let $\Z(1)$ denote the local system which on each $\Z/2$-space
$E\Z/2\times X$ is locally $\Z$, and has monodromy given by the
homomorphism
\[
\pi_{1}E\Z/2\times X\to \pi_{1}B\Z/2=\Z/2.
\]
We'll write $\hz(1)$ for the equivariant spectrum representing the
cohomology theory
\[
X\mapsto H^{\ast}_{\Z/2}(X;\Z(1))=H^{\ast}(E\Z/2\times_{\Z/2}X;\Z(1)),
\]
and more generally $HA(1)$ for the equivariant spectrum representing
\[
H^{\ast}(E\Z/2\times_{\Z/2}X;A\otimes\Z(1)).
\]
The associated Eilenberg-MacLane spaces are denoted $K(\Z(1),n)$, etc.
Since 
\[
\Z/2(1)=\Z/2,
\]
the spectrum $\hz/2(1)$ is just $\hz/2$, and
reduction modulo $2$ is represented by a map
\[
\hz(1)\to\hz/2(1).
\]

Note that if $M$ is a space and $t:\pi_{1}M\to\Z/2$ classifies a double
cover $\tilde M\to M$, and local system $\Z(1)$ on $M$, then 
\[
E\Z/2\times_{\Z/2} \tilde M\to M
\]
is a homotopy equivalence (it is a fibration with contractible fibers) 
and we have an isomorphism
\[
H^{\ast}(M;\Z(1))=
H^{\ast}_{\Z/2}(\tilde M;\Z(1)).
\]

As described in \S\ref{sec:proof-theorem-main} we define an
equivariant $\bsobwu$ by the homotopy pullback square
\[
\begin{CD}
\bsobwu @> \lambda >> K\left(\Z(1),2k \right) \\
@VVV @VVV \\
BSO @>>> K\left(\Z/2,2k \right),
\end{CD}
\]
and in this way give the associated Thom spectra $\msobwu$ and
\begin{multline*}
\msobwu\wedge K(\Z(1),2k)_{+} \\=
\thom\left(\bsobwu\times K(\Z(1),2k),\xi\oplus 0 \right)
\end{multline*}
$\Z/2$-actions.  The group
\begin{multline*}
\pi_{4k}\msobwukhp \\
=\pi_{4k}\thom\left(\left(\bsobwu\times K(\Z(1),2k) \right)\times_{Z/2}E\Z/2,
\xi\oplus 0 \right)
\end{multline*}
is the cobordism group of $4k$-dimensional oriented manifolds $M$,
equipped with a map $t:\pi_{1}M\to\Z/2$ classifying a local system
$\Z(1)$, a cocycle $\lambda\in Z^{2k}(M;\Z(1))$ whose mod
$2$-reduction represents the Wu-class $\nu_{2k}$, and a cocycle $x\in
Z^{2k}(M;\Z(1))$.  As in the non-equivariant case, the integral
Wu-structure $\lambda\in H^{2k}(M;\Z(1))$ defines the quadratic
function
\[
\phi_{4k}=\phi_{M,\lambda}(x)=\tfrac12\int_{M}(x^{2}-x\lambda)\qquad x\in
H^{2k}(M;\Z(1)),
\]
and a homomorphism
\[
\pi_{4k}\msobwukh\to\Z.
\]
An analysis of the Anderson dual of $\Sigma^{\infty}K(\Z(1),2k)$ 
leads to a canonical refinement of this to a map 
\[
\msobwukh\to \tilde I,
\]
and checking the rationalization sequence leads to the desired
quadratic refinement of the link pairing in dimension $(4k-1)$.

The the analysis of $\tilde I(K(\Z(1),2k))$ works more or less the
same as in the non-equivariant case.  Here are the main main points:

\begin{textList}
\item
With our conventions, a homotopy class of maps 
\[
\msobwukh\to \Sigma^{4k}\tilde I
\]
is the same as an {\em equivariant} homotopy class of maps
\[
\msobwuk\to \Sigma^{4k}\tilde I.
\]
We are therefore looking for an equivariant map
\begin{equation}\label{eq:t-59}
\msobwu\to \Sigma^{4k}\tilde I\left(K(\Z(1),2k)\right).
\end{equation}
\item 
By attaching cells of the form $D^{n}\times \Z/2$, The Postnikov
section $X\langle n,\dots,m \rangle$ of an equivariant spectrum or
space can be formed.  Its non-equivariant homotopy groups vanish
outside of dimensions $n$ through $m$, and coincide with those of $X$
in that range.  If $E$ is a space or spectrum with a free
$\Z/2$-action, and having only cells in dimensions less than $n$, then
\[
\left[E,X\langle n,\dots ,m \rangle \right]=0.
\]
The same is true if $E$ is an equivariant space or spectrum with
$\pi_{k}E=0$ for $k>m$.  Since the equivariant spectrum $\msobwu$ is
$(-1)$-connected, this means, as before, that we can
replace~\eqref{eq:t-59} with
\begin{multline*}
\msobwu\to \Sigma^{4k}\tilde I\left(K(\Z(1),2k)\right) 
\langle0,\dots,\infty\rangle\\
= \Sigma^{4k}\tilde I\left(\Sigma^{\infty}K(\Z(1),2k)\langle -\infty,\dots,4k
\rangle\right).
\end{multline*}

\item The cup product goes from
\[
H^{\ast}_{\Z/2}(X;\Z(1))\times H^{\ast}_{\Z/2}(X;\Z(1))\to 
H^{\ast}_{\Z/2}(X;\Z)
\]
and is represented by a map of spectra
\[
\hz(1)\wedge \hz(1)\to \hz.
\]
\item The adjoint of 
\[
\hz(1)\wedge \hz(1)\to \hz\to \tilde I
\]
is an equivariant weak equivalence 
\[
\hz(1)\to\tilde I\hz(1).
\]
\item The square
\[
\begin{CD}
\Sigma^{\infty} K(\Z(1),2k)\langle-\infty, 4k \rangle @>\iota\cup \iota>>
\Sigma^{4k}\hz \\
@V\iota VV @VVV \\
\Sigma^{2k}\hz(1) @>> \Sq^{2k}> \Sigma^{4k}\hz/2.
\end{CD}
\]
is a homotopy pullback square.

\item  As a consequence, so is
\[
\begin{CD}
\Sigma^{4k}\tilde I\left(K(\Z(1),2k) \right)\langle0, \infty \rangle @>\iota\cup \iota>>
\Sigma^{2k}\hz(1) \\
@V\iota VV @VVV \\
\hz @>> \chi \Sq^{2k}> \Sigma^{2k}\hz/2.
\end{CD}
\]
\end{textList}

The following proposition summarizes the main result of this
discussion.  As in \S\ref{sec:topol-theory-quadr} we'll use the notation
\begin{gather*}
B=\bsobwu\times_{\Z/2}E\Z/2 \\
\bar{B}=\bsowubar\times_{\Z/2}E\Z/2.
\end{gather*}

\begin{prop}\label{thm:t-4}
Let $N$ be a $\bar B$-oriented manifold of dimension $(4k-1)$.
Associated to the ``Wu-cocycle'' $\eta\in Z^{2k-1}(N;\qz)$ is a
quadratic function
\[
\phi=\phi_{N}=\phi_{N,\eta}:H^{2k}(N;\Z(1))_{\text{tor}}\to\qz
\]
whose associated bilinear form is the link pairing
\[
H^{2k}(N;\Z(1))_{\text{tor}}\times
H^{2k}(N;\Z(1))_{\text{tor}}\to \qz.
\]
If $M/\Delta^{1}$ is a $B$-oriented map of relative dimension $(4k-1)$,
equipped with compatible $\bar B$ orientations on $\partial_{0}M$ and
$\partial_{1}M$, and $x$ is an element of
$H^{2k}(M;\Z(1))$ whose restriction to $\partial M$ is torsion, then
\begin{equation}\label{eq:t-60}
\begin{gathered}
\phi_{\partial_{0}M}(\partial_{0}x)
-\phi_{\partial_{1}M}(\partial_{1}x)
\equiv \frac12\int_{M} (x^{2}-x\lambda)\mod \Z \\
\partial_{i}x=x\vert_{\partial_{i}M}.
\end{gathered}
\end{equation}
\qed
\end{prop}

\appendix

\section{Simplicial methods}\label{sec:simplicial-methods}

\subsection{Simplicial set and simplicial
objects}\label{sec:simpl-set-simpl} 

A {\em simplicial set} $\simp X$ consists of a sequence of sets $X_n$,
$n\ge 0$ together with ``face'' and ``degeneracy'' maps
\begin{align*}
d_i: X_n\to X_{n-1} &\qquad i=0,\dots,n \\
s_i: X_{n-1}\to X_{n} &\qquad i=0,\dots,n-1
\end{align*}
satisfying
\begin{gather*}
d_{j}d_{i} = d_{i} d_{j+1} \qquad j\ge i \\
s_{j}s_{i} = s_{i} s_{j-1} \qquad j > i \\
d_{j}s_{i} = \begin{cases}
             s_{i-1}d_{j}\qquad & j<i \\
             \text{identity}\qquad & j=i,i+1 \\
             s_{i}d_{j-1}\qquad & j>i+1 \\
             \end{cases}
\end{gather*}
The set $X_n$ is called the set of $n$-simplices of $\simp X$.

\begin{rem}
Let $\Delta$ be the category whose objects are the finite ordered sets
\[
[n]=\{0\le 1\dots\le n\},\qquad n\ge 0
\]
and whose morphisms are the order preserving maps.  The data
describing a simplicial set $\simp X$ is equivalent to the data
describing a contravariant functor
\[
X:\Delta\to\sets
\]
with $X[n]$ corresponding to the set $X_{n}$ of $n$-simplices, and the
face and degeneracy maps $d_i$ and $s_i$ corresponding to the values of
$X$ on the maps
\begin{equation}\label{eq-face-degeneracy}
\begin{aligned}
{}[n-1]&\hookrightarrow [n] \\
{}[n]&\twoheadrightarrow [n-1]
\end{aligned}
\end{equation}
which, respectively ``skip $i$'' and ``repeat $i$.''  
\end{rem}

\begin{rem}
A simplicial object in a category $\cat{C}$ is a contravariant functor
\[
\simp X:\Delta\to \cat{C}.
\]
Thus one speaks of simplicial abelian groups, simplicial Lie algebras, etc.
\end{rem}

One basic example of a simplicial set is the {\em singular
complex}, $\sing S$ of a space $S$, defined by
\[
(\sing S)_n = \text{set of maps from $\Delta^n$ to $S$},
\]
where 
\[
\Delta^{n}=\{(t_{0},\dots,t_{n})\mid 0\le t_{i}\le 1, \sum t_{i}=1 \}
\]
is the standard $n$-simplex with vertices
\[
e_{i}=(0,\dots,\overset{i}{1},\dots,0)\qquad i=0,\dots,n.
\]
The face and degeneracy maps are derived from the linear extensions
of~\eqref{eq-face-degeneracy}.

\begin{defin}
The {\em geometric realization} $|\simp X|$ of a simplicial set $\simp X$
is the space
\[
\coprod X_n\times\Delta^n/\sim
\]
where $\sim$ is the equivalence relation generated by
\begin{gather*}
(d_i x,t) = (x,d^i t) \\
(s_i x,t) = (x,s^i t)
\end{gather*}
and 
\begin{gather*}
d^i:\Delta^{n-1}\to\Delta^n \\
s^i:\Delta^{n}\to\Delta^{n-1} 
\end{gather*}
are the linear extensions of the maps~\eqref{eq-face-degeneracy}.
\end{defin}

The evaluation maps assemble to a natural map
\begin{equation}\label{eq:75}
|\sing S|\to S
\end{equation}
which induces an isomorphism of both singular homology and homotopy
groups.

The {\em standard (simplicial) $n$-simplex} is the simplicial set
$\Delta^{n}_{\bullet}$, given as a contravariant functor by
\[
\Delta^{n}([m])=\Delta\left([m],[n] \right).
\]
One easily checks that 
\[
\left|\Delta^{n}_{\bullet} \right|=\Delta^{n}.
\]

A {\em simplicial homotopy} between maps $f:X_{\bullet}\to
Y_{\bullet}$ is a simplicial map 
\[
h:X_{\bullet}\times \Delta^{1}_{\bullet}\to Y_{\bullet}
\]
for which $\partial_{1}h=f$, and $\partial_{0}h=g$.

\subsection{Simplicial homotopy groups}\label{sec:simpl-homot-groups}

\begin{defin}
Let $\simp X$ be a pointed simplicial set.  The set
$\spi_{n}\simp X$ is the set of simplicial maps 
\[
(\Delta^{n}_{\bullet},\partial\Delta^{n}_{\bullet})\to (\simp X,\ast)
\]
modulo the relation of simplicial homotopy
\end{defin}

There is a map 
\begin{equation}\label{eq:7}
\spi_{n}\simp X\to \pi_{n}(|\simp X|).
\end{equation}

The {\em $k$-horn} of $\Delta^{n}_{\bullet}$ is the simplicial set 
\[
V^{n,k}_{\bullet}=\bigcup_{i\ne
k}\partial_{i}\Delta^{n}_{\bullet}\subset\partial \Delta^{n}_{\bullet}.
\]

\begin{defin}\label{def:6}
A simplicial set $X_{\bullet}$ satisfies the {\em Kan extension
condition} if for every $n,k$, every map $V^{n,k}_{\bullet}\to
X_{\bullet}$ extends to a map $\Delta^{n}_{\bullet}\to X_{\bullet}$.
\end{defin}

\begin{prop}\label{thm:6}
If $\simp X$ satisfies the Kan extension condition then~\eqref{eq:7}
is an isomorphism.
\end{prop}

\subsection{Simplicial abelian groups}\label{sec:simpl-abel-groups}

The category of simplicial abelian groups is equivalent to the
category of chain complexes, making the linear theory of simplicial
abelian groups particularly simple.  The equivalence associates to a
simplicial abelian group $A$ the chain complex $NA$ with
\begin{align*}
NA_n &= \left\{a\in A_n\mid d_i a =0,\quad i=1,\dots n\right\} \\
d &= d_0.
\end{align*}
The inverse correspondence associates to a chain complex $C_{\ast}$,
the simplicial abelian group with $n$-simplices
\[
\bigoplus_{\substack{f:[n]\to[k]\\\text{$f$ surjective}}} C_k.
\]

\begin{rem}\label{rem:11}
The normalized complex $NA$ of a simplicial abelian group $A$ is a
subcomplex of the simplicial abelian group $A$ regarded as a chain
complex with differential $\sum (-1)^{i}\partial_{i}$.  In fact $NA$
is a retract of this chain complex, and the complementary summand is
contractible.  It is customary not to distinguish in notation between a
simplicial abelian group $A$ and the chain complex just described.  
\end{rem}

It follows immediately from the definition that 
\[
\spi_{n} A= H_{n}(NA).
\]
Moreover, any simplicial abelian group satisfies the Kan extension condition.
Put\-ting this together gives

\begin{prop}\label{thm:5}
For a simplicial abelian group $\simp A$, the homotopy groups of
$|\simp A|$ are given by
\[
\pi_{n}|\simp A|= H_{n}(NA).
\]
\qed
\end{prop}

Fix a fundamental cocycle $\iota\in Z^{n}\left(K(\Z,n);\Z \right)$.  

\begin{prop}\label{thm:30}
The map 
\begin{align*}
\sing K(A,n)^{X} &\to Z^{n}(X\times \Delta^{\bullet}) \\
f &\mapsto f^{\ast}\iota
\end{align*}
is a simplicial homotopy equivalence.
\end{prop}

\begin{pf}
On path components, the map induces the isomorphism
\[
[X,K(\Z,n)]\approx H^{n}(X;\Z).
\]
One can deduce from this that the map is an isomorphism of higher
homotopy groups by replacing $X$ with $\Sigma^{k}X$, and showing that
the map ``integration along the suspension coordinates'' gives a
simplicial homotopy equivalence
\[
Z^{n}(\Sigma^{k}X\times \Delta^{\bullet})\to
Z^{n-k}(X\times \Delta^{\bullet}).
\]
This latter equivalence is a consequence of Corollary~\ref{thm:22}.  
\end{pf}

\section{Picard categories and Anderson duality}
\label{sec:picard-categ-anders}

\subsection{Anderson duality}
\label{sec:anderson-duality}

In this appendix we will work entirely with spectra, and we will use
the symbol $\pi_{k}X$ to denote the $k^{\text{th}}$ homotopy group of
the spectrum.  In case $X$ is the suspension spectrum of a space $M$
this group is the $k^{\text{th}}$ {\em stable} homotopy of $M$.

To define the Anderson dual of the sphere, first note that the functors
\begin{align*}
X\mapsto\hom(\pi_{\ast}^{\text{st}}X,\Q) \\
X\mapsto\hom(\pi_{\ast}^{\text{st}}X,\qz)
\end{align*}
satisfies the Eilenberg-Steenrod axioms and so define cohomology
theories which are represented by spectra.  In the first case this
cohomology theory is just ordinary cohomology with coefficients in the
rational numbers, and the representing spectrum is the
Eilenberg-MacLane spectrum $H\Q$.  It the second case the spectrum is
known as the {\em Brown-Comenetz dual of the sphere}, and denoted
$I$~\cite{BC}.  There is a natural map
\begin{equation}\label{eq:t-67}
H\Q\to I
\end{equation}
representing the transformation
\[
\hom(\pi_{\ast}^{\text{st}}M,\Q)\to
\hom(\pi_{\ast}^{\text{st}}M,\qz).
\]
\begin{defin}
The {\em Anderson dual of the sphere}, $\tilde I$, is the homotopy
fiber of the map~\eqref{eq:t-67}.
\end{defin}

By definition there is a long exact sequence
\begin{multline*}
\cdots\hom(\pi_{n-1}X,\qz)\to
\tilde I^{n}(X) \\ \to \hom(\pi_{n}X,\Q)
\to
\hom(\pi_{n}X,\qz)\cdots,
\end{multline*}
from which one can extract a short exact sequence
\begin{equation}\label{eq:t-68}
\ext\left(\pi_{n-1}X,\Z \right)\rightarrowtail
\tilde I^{n}(X)\twoheadrightarrow
\hom\left(\pi_{n}X,\Z \right).
\end{equation}
The sequence~\eqref{eq:t-68} always splits, but not canonically.  

By means of~\eqref{eq:t-68}, the group $\hom\left(\pi_{n}X,\Z \right)$
gives an algebraic approximation to $\tilde I^{n}X$.  In the next two
sections we will refine this to an algebraic description of all of
$\tilde I^{n}(X)$ (Corollary~\ref{thm:t-27}).

For a spectrum $E$, let $\tilde I(E)$ denote the function spectrum of
maps from $E$ to $\tilde I$.  Since
\[
\left[X,\tilde I(E) \right] = 
\left[X\wedge E,\tilde I \right],
\]
there is a splittable short exact sequence
\[
\ext\left(\pi_{n-1}X\wedge E,\Z \right)\rightarrowtail
\tilde I^{n}(X\wedge E)\twoheadrightarrow
\hom\left(\pi_{n}X\wedge E,\Z \right).
\]
If $E=H\Z$, then $\tilde I(E)=H\Z$, and the above sequence is simply
the universal coefficient sequence.  One can also check that $\tilde
I(H\Z/2)=\Sigma^{-1}H\Z/2$, that the Anderson dual of mod $2$
reduction
\[
H\Z\to H\Z/2
\]
is the Bockstein
\[
\Sigma^{-1}H\Z/2\to H\Z,
\]
and that the effect of Anderson duality on the Steenrod algebra
\[
\left[H\Z/2,H\Z/2 \right]_{\ast}
\]
is given by the canonical anti-automorphism $\chi$.

\begin{rem}
If $X$ is a pointed space, we will abbreviate $\tilde
I\left(\Sigma^{\infty}X \right)$ to $ \tilde I(X)$.
\end{rem}

\subsection{Picard categories}
\label{sec:picard-categories}

\begin{defin}
A {\em Picard category} is a groupoid equipped with the structure of a
symmetric monoidal category in which each object is invertible.
\end{defin}

A functor between Picard categories is a functor $f$ together with a
natural transformation
\[
f(a\otimes b)\to f(a)\otimes f(b)
\]
which is compatible with the commutativity and associativity laws.
The collection of Picard categories forms a $2$-category%
\footnote{A $2$-category is a category in which the morphisms form a
category}.%
If $\mathcal C$ is a Picard category, then the set $\pi_{0}\mathcal C$
of isomorphism classes of objects in $\mathcal C$ is an abelian group.
For every $x\in\mathcal C$, the map
\begin{align*}
\Aut\left(e \right) &\to \Aut(e\otimes x)\approx \Aut(x) \\
f &\mapsto f\otimes\text{Id}
\end{align*}
is an isomorphism ($e$ is the unit for $\otimes$).  This group is also
automatically abelian, and is denoted $\pi_{1}\mathcal C$.

For each $a\in\mathcal C$, the symmetry of $\otimes$ gives an
isomorphism 
\[
\left(a\otimes a\to a\otimes a \right)\in\Aut(a\otimes
a)\approx\Aut(e)=\pi_{1}\mathcal C
\]
which squares to $1$.  Write $\epsilon(a)$ for this element.  The
invariant $\epsilon$ descends to a homomorphism 
\[
\pi_{0}\mathcal C\otimes\Z/2
\to \pi_{1}\mathcal C,
\]
known as the {\em $k$-invariant} of $\mathcal C$.  The Picard category
$\mathcal C$ is determined by this invariant, up to equivalence of
Picard categories (see Proposition~\ref{thm:t-26} below).  A Picard
category with $\epsilon=0$ is said to be a {\em strict} Picard
category.  

\begin{eg}
A map $\partial:A\to B$ of abelian groups determines a strict Picard
category with $B$ as the set of objects, and in which a morphism from
$b_{0}$ to $b_{1}$ is a map for which $\partial a=b_{1}-b_{0}$.  The
operation $+$ comes from the group structure, and the natural
transformations $\psi$ and $\alpha$ are the identity maps.  We'll
denote this groupoid with the symbol
\[
\left(A\to B \right)=
\left(A\xrightarrow{\partial}{} B \right).
\]
It is immediate from the definition that 
\begin{align*}
\pi_{0}\left(A\xrightarrow{\partial}{} B \right) &= \coker \partial \\
\pi_{1}\left(A\xrightarrow{\partial}{} B \right) &= \ker \partial.
\end{align*}
One can easily check that the Picard category $\left(A\to B \right)$
is equivalent, as a Picard category, to
$\left(\pi_{1}\xrightarrow{0}{}\pi_{0} \right)$.  Every strict Picard
category is equivalent to one of this form.
\end{eg}

\begin{eg}
If $E$ is a spectrum, then each of the fundamental groupoids
\[
\pi_{\le 1}E_{n}
\]
is a Picard category, in which the $\otimes$ structure comes from the
``loop multiplication'' map 
\[
E_{n}\times E_{n}\to E_{n}.
\]
To make this precise, first note that the space $E_{n}$ is the space
of maps of spectra
\[
S^{0}\to \Sigma^{n}E,
\]
and a choice of ``loop multiplication'' amounts to a choice of
deformation of the diagonal $S^{0}\to S^{0}\times S^{0}$ to a map
$S^{0}\to S^{0}\vee S^{0}$. 

Let $\mathcal L_{k}(\R^{n})$ denote the space of $k$-tuples of 
linearly embedded $n$-cubes in $I^{n}$.  The ``Pontryagin-Thom
collapse'' gives a map from $\mathcal L_{k}(\R^{n})$ to the space of
deformations of the diagonal (here $S^{n}$ denotes the {\em space} $S^{n}$)
\[
S^{n}\to \prod^{k}S^{n}
\]
to a map 
\[
S^{n}\to \bigvee^{k}S^{n}.
\]
Set
\[
\mathcal L_{k}=\varinjlim_{n\to\infty}L_{k}(\R^{n}).  
\]
The space $L_{k}$ is contractible, and by passing to the limit from
the above, parameterizes deformations of the iterated diagonal map of
{\em spectra}
\begin{equation}\label{eq:73}
S^{0}\to \prod^{k}S^{0}
\end{equation}
to 
\begin{equation}\label{eq:74}
S^{0}\to\bigvee^{k}S^{0.}
\end{equation}

It is easy to see that a choice of point $x\in \mathcal L_{2}$ for the
monoidal structure, a path $x\to \tau x$ ($\tau$ the transposition in
$\Sigma_{2}$) for the symmetry, and a path in $\mathcal L_{3}$ for the
associativity law give the fundamental groupoid
\[
\pi_{\le 1}E_{n}
\]
the structure of a Picard category.  This Picard category structure is
independent, up to equivalence of Picard categories, of the choice of
point $x$, and the paths.
\end{eg}

\begin{rem}\label{rem:21}
The space of all deformations of~\eqref{eq:73} to~\eqref{eq:74} is
contractible, and can be used for constructing the desired Picard
category structure.  The ``little cubes'' spaces were introduced
in~\cite{may72,boardman73:_homot}, and have many technical advantages.
As we explain in the next example, in the case of Thom spectra they
give compatible Picard category structures to the the fundamental
groupoids of both the transverse and geometric singular complexes.
\end{rem}

\begin{eg}\label{eg:1}
Suppose $X=\thom\left(B;V \right)$ is a Thom spectrum.
Let $\gsing X_{-n}^{S}$ be the simplicial set  with $k$-simplices the set
of $B$-oriented maps $E\to S\times\Delta^{k}$.  The groupoid
\[
\pi_{\le 1}\gsing X_{-n}^{S}
\] 
then has objects the $n$-dimensional $B$-oriented maps $E\to S$, and
in morphisms from $E_{0}$ to $E_{1}$  equivalence
classes of $B$-oriented maps
\[
p:E\to S\times\Delta^{1}
\]
with $p\vert_{\partial_{i}\Delta^{1}}=E_{i}/S$, $i=0,1$.  In principle
$\pi_{\le 1}\gsing X_{-n}^{S}$ is a Picard category with the the
$\otimes$-structure coming from disjoint union of manifolds.  But this
requires a little care.  Recall that a $B$-orientation of $E/S$
consists of an embedding $E\subset S\times\R^{N}\subset
S\times\R^{\infty}$, a tubular neighborhood $W\subset S\times\R^{N}$,
and a map $E\to B$ classifying $W$.  To give a $B$-orientation to the
disjoint union of two $B$-oriented maps, $E_{1}\subset S\times \R^{N}$
and $E_{2}\subset S\times \R^{M}$ we need to construct an embedding
\[
E_{1}\amalg E_{2}\subset W_{1}\amalg W_{2}\subset S\times \R^{P}
\]
for some $P$.  Regard $\R^{k}$ as the interior of $I^{k}$, and choose
a point
\[
I^{P}\amalg I^{P}\subset I^{P}\qquad P\gg0.
\]
in $\mathcal L_{2}=\varinjlim L_{2}(\R^{P})$.  Making sure $P>N,M$, we
can then use
\[
W_{1}\amalg W_{2}\subset I^{P}\amalg I^{P}\subset I^{P}.
\]
The rest of the data needed to give $E_{1}\amalg E_{2}$ a
$B$-orientation is then easily constructed.   The Pontryagin-Thom construction
gives an equivalence of simplicial sets
\[
\gsing X_{-n}^{S} \to \sing X_{-n}^{S}
\]
and an equivalence of Picard categories 
\[
\pi_{\le 1}\gsing X_{-n}^{S} \to \pi_{\le 1}\sing X_{-n}^{S}=\pi_{\le
1}X_{n}^{S}.
\]
\end{eg}

All Picard categories arise as the fundamental groupoid of a spectrum:
\begin{prop}\label{thm:t-26}
The correspondence 
\[
E\mapsto \pi_{\le 1}E_{0}
\]
is an equivalence between the $2$-category of spectra $E$ satisfying 
\[
\pi_{i}E=0\quad i\ne 0,1
\]
and the $2$-category of Picard categories.  (The collection of such
spectra $E$ is made into a $2$-category by taking as morphisms the
fundamental groupoid of the space of maps.)
\end{prop}

\begin{pf}
We merely sketch the proof.  The homotopy type of such a spectrum $E$
is determined by its $k$-invariant
\begin{equation}\label{eq:t-33}
H\pi_{0}E\to \Sigma^{2}H\pi_{1}E.
\end{equation}
Now the set of homotopy classes of maps~\eqref{eq:t-33} is naturally
isomorphic to 
\[
\hom\left(A\otimes\Z/2,B \right).
\]
One can associate to a Picard category with $k$-invariant $\epsilon$, the
spectrum $E$ with the same $k$-invariant.  Using this it is not
difficult to then check the result.
\end{pf}

\begin{rem}
Clearly, the functor
\[
E\mapsto \pi_{\le 1}E_{n-1}
\]
gives an equivalence between the $2$-category of Picard categories,
and the $2$-category of spectra $E$ whose only non-zero homotopy
groups are $\pi_{n-1}E$ and $\pi_{n}E$.
\end{rem}

\subsection{Anderson duality and functors of Picard categories}
\label{sec:anders-dual-funct}

We now relate this discussion to Anderson duality.
First a couple of easy observations.  

\begin{lem}
The maps
\begin{gather*}
[X,\Sigma^{n}\tilde I] \xrightarrow{}{} 
\left[X\langle n-1,\infty \rangle,\Sigma^{n}\tilde I
\right]\xleftarrow{}{}\left[X\langle n-1,n \rangle,\Sigma^{n}\tilde I\right]
\\
\left[X\langle n-1,n \rangle,\Sigma^{n}\tilde I\langle n-1,\infty
\rangle\right] \xrightarrow{}{} \left[X\langle n-1,0 \rangle,\Sigma^{n}\tilde
I\right]
\end{gather*}
are isomorphisms.
\end{lem}

\begin{pf}
This follows easily from the exact sequence
\[
\ext\left(\pi_{n-1}X,\Z \right)\rightarrowtail
\tilde I^{n}(X)\twoheadrightarrow
\hom\left(\pi_{n}X,\Z \right).
\]
\end{pf}

Since 
\[
\Sigma^{n}\tilde I\langle n-1,\dots,\infty \rangle\approx \Sigma^{n}H\Z,
\]
we have 
\begin{cor}\label{thm:t-24}
There is a natural isomorphism 
\[
\tilde I^{n}(X)\approx H^{n}\left(X\langle n-1,n \rangle;\Z \right).
\]
\end{cor}

Now the Picard category
\[
\pi_{\le 1}\left(\Sigma^{n}\tilde I \right)_{n-1}\approx
\pi_{\le 1}\left(\Sigma^{n}H\Z)_{n-1}\right)
=\pi_{\le 1}K(\Z,1)
\]
is canonically equivalent to
\[
\left(\Z\to 0 \right)\approx
\left(\Q\to \qz \right).  
\]
Combining this with Proposition~\ref{thm:t-26} then gives
\begin{cor}\label{thm:t-27}
The group $\tilde I^{n}(X)$ is naturally isomorphic to the group of
natural equivalences classes of functors of Picard categories
\[
\pi_{\le 1}X_{n-1}\to\left(\Q\to\qz \right).
\]
\qed
\end{cor}

\section{Manifolds with corners}\label{sec:manif-with-corn} In this
section we briefly review the basics of the theory of manifolds with
corners.  For more information
see~\cite{cerf61:_topol,douady61:_variet,cerf62:_la_diff_s} and for a
discussion in connection with cobordism see~\cite{laures00}.  In our
discussion we have followed
Melrose~\cite{melrose:_differ_analy_manif_corner}.

\subsection{$t$-manifolds}\label{sec:t-manifolds}
Let $\R^{n}_{k}\subset \R^{n}$ be the subspace
\[
\R^{n}_{k}=\left\{(x_{1},\dotsc, x_{n})\in\R^{n}\mid x_{i}\ge 0, 1\le
i\le k \right\}. 
\]

\begin{defin}
The set of points $x=(x_{1},\dots,x_{n})\in\R^{n}_{k}$ for which
exactly $j$ of $\{x_{1},\dots,x_{k}\}$ are $0$ is the
set of {\em $j$-corner points} of $\R^{n}_{k}$
\end{defin}

\begin{defin}
A function $f$ on an open subset $U\subset\R^{n}_{k}$ is {\em smooth}
if and only if $f$ extends to a $C^{\infty}$ function on a
neighborhood $U'$ with
\[
U\subset U'\subset \R^{n}.
\]
A map $g$ from an open subset $U\subset\R^{n}_{k}$ to an open subset
$U'\subset \R^{m}_{l}$ is {\em smooth} if $f\circ g$ is smooth
whenever $f$ is smooth (this only need be checked when $f$ is one of the
coordinate functions).
\end{defin}

\begin{rem}
The definition of smoothness can be characterized without reference to
$U'$.  It is equivalent to requiring that $f$ be $C^{\infty}$ on
$(0,\infty)^{k}\times \R^{n-k}\cap U$ with all derivatives bounded on
all subsets of the form $K\cap U$ with $K\subset \R^{n}_{k}$ compact.
\end{rem}

\begin{defin}
A {\em chart with corners} on a space $U$ is a homeomorphism
$\phi:U\to\R^{n}_{k}$ of $U$ with an open subset of $\R^{n}_{k}$.  Two
charts $\phi_{1},\phi_{2}$ on $U$ are said to be {\em compatible} if
the function $\phi_{1}\circ\phi_{2}^{-1}$ is smooth.
\end{defin}

\begin{lem}
Suppose that $\phi_{1}$ and $\phi_{2}$ are two compatible charts on a
space $U$.  If $x\in U$, and $\phi_{1}(x)$ is a $j$-corner point, then
$\phi_{2}(x)$ is a $j$-corner point.
\qed
\end{lem}

\begin{defin}
An {\em atlas (with corners)} on a space $X$ is a collection of pairs
$(\phi_{a}, U_{a})$ with $\phi_{a}$ a chart on $U_{a}$, $U_{a}$ a
cover of $X$, and for which $\phi_{a}$ and $\phi_{a'}$ are compatible
on $U_{a}\cap U_{a'}$.  A {\em $C^{\infty}$-structure with corners} on
$X$ is a maximal atlas.  
\end{defin}

\begin{defin}
A {\em $t$-manifold}  is a paracompact Hausdorff space $X$ together
with a $C^{\infty}$-structure with corners.
\end{defin}

\begin{defin}
A {\em $j$-corner point} of a $t$-manifold $X$ is a point $x$ which in
some (hence any) chart corresponds to a $j$-corner point of $\R^{n}_{k}$.
\end{defin}

Let $U$ be an open subset of $\R^{n}_{k}$.  The {\em tangent} bundle
to $U$ is the restriction of the tangent bundle of $\R^{n}$ to $U$,
and will be denoted $TU$.  A smooth map $f:U_{1}\to U_{2}$ has a
derivative $df:TU_{1}\to TU_{2}$.  It follows that if $X$ is a
$t$-manifold, the tangent bundles to the charts in an atlas patch
together to form a vector bundle over $X$, the tangent bundle to $X$.

If $x\in X$ is a $j$-corner point, then the tangent space to $X$ at
$x$ contains $j$ distinguished hyperplanes $H_{1},\dots,H_{j}$ in
general position, equipped with orientations of $T_{x}X/H_{i}$.  The
intersection of these hyperplanes is the {\em tangent space to the
corner at $x$} and will be denoted $bT_{x}$.  In a general
$t$-manifold, these hyperplanes need not correspond to distinct
components of the space of $1$-corner points, as, for example happens
with the polar coordinate region
\begin{gather*}
0\le \theta\le \pi/2 \\
0\le r\le \sin(2\theta).
\end{gather*}
In a {\em manifold with corners} (Definition~\ref{def:1}) a global
condition is imposed, which guarantees that the tangent hyperplanes do
correspond to distinct components.

\subsection{Neat maps and manifolds with corners}

\begin{defin}
A smooth map of $t$-manifolds $f:X\to Y$ is {\em neat} if $f$ maps
$j$-corner points to $j$-corner points, and if for each $j$-corner
point, the map
\[
df: T_{x}/bT_{x} \to T_{f(x)}/bT_{f(x)}
\]
is an isomorphism.
\end{defin}

The term {\em neat} is due to Hirsch~\cite[p. 30]{hirsch94:_differ},
who considered the case of embeddings of manifolds with boundary.

\begin{rem}
The condition on tangent spaces is simply the requirement that $f$ be
transverse to the corners.
\end{rem}

\begin{eg}
The map 
\begin{align*}
I^{m}&\to I^{m+1} \\
(x) &\mapsto (x,t)
\end{align*}
is ``neat'' if and only if $t\ne 0,1$.
\end{eg}

\begin{defin}
\label{def:1}
A $t$-manifold $X$ is a {\em manifold with corners} if there
exists a neat map $X\to \Delta^{n}$ for some $n$.
\end{defin}

\begin{rem}
The existence of a neat map $X\to \Delta^{n}$ for some $n$ is equivalent
to the existence of a neat map $X\to I^{m}$ for some $m$.
\end{rem}

\begin{rem}
In this paper we have restricted our discussion to manifolds with
corners, for sake of simplicity, but in fact our results apply to neat
maps of $t$-manifolds.
\end{rem}

\begin{prop}
The product of manifolds with corners is a manifold with corners.
If $X\to M$ is a neat map, and $M$ is a manifold with corners, then so
is $X$.\qed
\end{prop}

\subsection{Normal bundles and tubular
neighborhoods}\label{sec:normal-bundles}

Suppose that $f:X\to Y$ is a neat embedding of $t$-manifolds.  The
relative {\em normal bundle} of $f$ is the vector bundle over $X$
whose fiber at $x$ is $T_{f(x)}Y/df\left(T_{x}X \right)$.  An
embedding of $t$-manifolds does not necessarily admit a tubular
neighborhood, as the example
\begin{align*}
[0,1] &\to [0,1]\times \R \\
x & \mapsto (x,\sqrt{x})
\end{align*}
shows.  However, when $f$ is a neat embedding, the fiber of the
relative normal bundle at $x$ can be identified with
$bT_{f(x)}Y/df(bT_{x}X)$, and the usual proof of the existence of a
normal bundle applies.

\begin{prop}
Suppose that $X$ is a compact $t$-manifold, and $f:X\to Y$ is a neat
embedding.  There exists a neighborhood $U$ of $f(X)$, a vector bundle
$W$ over $X$, and a diffeomorphism of $W$ with $U$ carrying the
zero section of $W$ to $f$. \qed
\end{prop}

The neighborhood $W$ is a tubular neighborhood of $f(X)$, and the derivative
of the embedding identifies $W$ with the relative normal bundle of $f$
\[
W_{x}\approx T_{f(x)}Y/df\left(T_{x}X \right)\approx
bT_{f(x)}Y/df(bT_{x}X).
\]

\begin{prop}
Let $f:X\to Y$ be a neat map of compact $t$-manifolds.  There exists
$N\gg0$ and a factorization
\[
X\to Y\times\R^{N}\to Y
\]
of $f$ through a neat embedding $X\hookrightarrow Y\times\R^{N}$.\qed
\end{prop}

\subsection{Transversality}\label{sec:transversality}

\begin{defin}
A {\em vector space with $j$-corners} is a real vector space $V$
together with $j$ codimension $1$, relatively oriented subspaces $\{H_{1},\dots,H_{j} \}$
in general position:  
\begin{thmList}
\item each of the $1$-dimensional vector spaces
$V/H_{i}$ is equipped with an orientation;
\item for all $\alpha\subseteq\{1,\dots,j \}$, $\codim H_{\alpha}=s$,
where $H_{\alpha}=\bigcap_{i\in\alpha}H_{i}$.
\end{thmList}
\end{defin}

\begin{eg}
The tangent space to each point of a $t$-manifold is a vector space
with corners.
\end{eg}

Suppose that $(V; H_{1},\dots,H_{j})$ and $(V'; H'_{1},\dots,H'_{l})$
are two vector spaces with corners, and $W$ is a vector space.
\begin{defin}
Two linear maps
\[
V\xrightarrow{}{} W\xleftarrow{}{}V'
\]
are {\em transverse} if for each $\alpha\subset \{1,\dots,j \}$,
$\beta\subset\{1,\dotsc ,l \}$ the map
\[
H_{\alpha}\oplus H'_{\beta}\to W
\]
is surjective.
\end{defin}

The space of transverse linear maps is a possibly empty open subset of
the space of all maps.

\begin{defin}
Suppose that $Z$ and $M$ are $t$-manifolds, and that $S$ is a closed
manifold.  Two maps
\[
Z\xrightarrow{f}{} S\quad\text{and}\quad
M\xrightarrow{g}{} S
\]
are {\em transverse} if for each $z\in Z$ and $m\in M$ with
$f(z)=g(m)$, the maps
\begin{gather*}
T_{z}Z \xrightarrow{Df}{} T_{s}S \xleftarrow{Dg}{} T_{m}M
\qquad s=f(z)=g(m)
\end{gather*}
are transverse.
\end{defin}

\begin{lem}
The class of ``neat'' maps is stable under composition and transverse
change of base.
\end{lem}

\section{Comparison of $\HZ(n)^{k}(S)$ and
$\chh{n}^{k}(S)$}\label{sec:comp-checkhznks-chsn}

The purpose of this appendix to prove that the differential cohomology
groups defined in \S\ref{sec:diff-funct-spac} using the
Eilenberg-MacLane spectrum coincide with the differential cohomology
groups defined in \S\ref{subsec-differential-character} using the
cochain complex.  

The proof is in two steps.  The first is to replace the differential
function space
\[
\fil_{k-n}\left(K(\Z,k);\iota_{k} \right)^{S}
\]
with a simplicial abelian group.  To describe the group, recall from
\S\ref{subsec-differential-character} the cochain complex
$\chzero^{\ast}(S)$ with
\[
\chzero^{k}(S)=C^{k}(S;\Z)\times C^{k-1}(S;\R)\times\Omega^{k}(S),
\]
and
\[
d(c,h,\omega) = (\delta c, \omega-c-\delta h, d\omega).
\]
Let $\chzerococycle^{k}(S)$ be the group of $k$-cocycles in
$\chzero^{\ast}(S)$.  We define
\[
\fil_{s}\chzero^{\ast}(S\times
\Delta^{k})=C^{k}(S\times\Delta^k;\Z)\times
C^{k-1}(S\times\Delta^k;\R) \times\fil_{s}\Omega^{k}(S\times\Delta^k)
\]
(see the discussion preceding Definition~\ref{def:7}).  The
coboundary map in $\chzero^{\ast}(S\times \Delta^{k})$ carries the
$\fil_{s}$ to $\fil_{s+1}$.  Let 
\[
\fil_{s}\chzerococycle^{k}(S\times \Delta^{m})=
\fil_{s}\chzero^{k}(S\times \Delta^{m})\cap
\chzerococycle^{k}(S\times \Delta^{m}).
\]
Thus $\fil_{s}\chzerococycle^{k}(S\times\Delta^{m})$ consists of triples
$(c,h,\omega)$ for which $c$ and $\omega$ are closed, 
\[
\delta h =\omega-c,
\]
and for which the weight filtration of $\omega$ is less than or equal
to $s$.

Associating to a $k$-simplex of $\fil_{s}\left(K(\Z,k);\iota_{k}
\right)^{S}$ its underlying differential cocycle gives a map of
simplicial sets 
\begin{equation}\label{eq:77}
\fil_{s}\left(K(\Z,k);\iota_{k}\right)^{S}\to
\fil_{s}\chzerococycle^{k}\left(S\times \Delta^{\bullet} \right).
\end{equation}

\begin{lem}\label{thm:29}
The map~\eqref{eq:77} is a weak equivalence.
\end{lem}

\begin{pf}
By definition, the square
\[
\begin{CD}
\fil_{s}\left(K(\Z,k);\iota_{k}\right)^{S} @>>> 
\fil_{s}\chzerococycle^{k}\left(S\times \Delta^{\bullet} \right) \\
@VVV @VVV \\
\sing K(\Z,k)^{S} @>>> Z^{k}(S\times \Delta^{\bullet})
\end{CD}
\]
is a pullback square.  The left vertical map is a surjective map of
simplicial abelian groups, hence a Kan fibration, so the square is in
fact homotopy Cartesian.  The bottom map is a weak equivalence by
Proposition~\ref{thm:30}.  It follows that the top map is a weak
equivalence as well.
\end{pf}

For the second step note that ``slant product with the
fundamental class of the variable simplex'' gives a map from 
the chain complex associated to the simplicial abelian group
\begin{equation}\label{eq:79}
\fil_{k-n}\chzerococycle^{k}\left(S\times \Delta^{\bullet} \right)
\end{equation}
to  
\begin{equation}\label{eq:80}
\chcocycle{n}^{k}(S)
\xleftarrow{}{}
\chcochain{n}^{k-1}(S)\dots 
\xleftarrow{}{}
\chcochain{n}^{0}(S).
\end{equation}

\begin{prop}\label{thm:31}
The map described above, from~\eqref{eq:79} to~\eqref{eq:80} is a
chain homotopy equivalence.
\end{prop}

Together Propositions~\ref{thm:29} and~\ref{thm:31} then give the
following result
\begin{prop}
The map ``slant product with the fundamental class of the variable
simplex'' gives an isomorphism
\[
\HZ(n)^{k}(S)\xrightarrow{\approx}{} \chh{n}^{k}(S). \qed
\]
\end{prop}

\begin{pf*}{Proof of Proposition~\ref{thm:31}}
Let $A_{\ast}$ be the total complex of the following bicomplex:
\begin{equation}\label{eq:20}
\xymatrix{
\ar@{}[d]|{\vdots}  &\ar@{}[d]|{\vdots}  & \\
\fil_{k-n}\chzerococycle^{k}(S\times\Delta^{2})
\ar[d]_{\sum(-1)^{i}\partial_{i}^{\ast}}
& \ar[l]_{d} \ar[d] \fil_{k-n-1}\chzero^{k-1}(S\times\Delta^{2})
& \ar[l] \dots \\
\fil_{k-n}\chzerococycle^{k}(S\times\Delta^{1}) \ar[d]
& \ar[l] \ar[d] \fil_{k-n-1}\chzero^{k-1}(S\times\Delta^{1})
& \ar[l] \dots \\
\fil_{k-n}\chzerococycle^{k}(S\times \Delta^{0}) & 
\ar[l]\fil_{k-n-1}\chzero^{k-1}(X\times\Delta^{0})
& \ar[l] \dots }
\end{equation}
By Lemma~\ref{thm:33} below, the inclusion of the leftmost column
\[
\fil_{s}\chzerococycle^{k}\left(S\times \Delta^{\bullet} \right)
\to A_{\ast}
\]
is a quasi-isomorphism.  By Lemma~\ref{thm:34} below, the inclusion
$B^{(0)}_{\ast}\subset A_{\ast}$ is also a quasi-isomorphism.  One
easily checks that map ``slant-product with the fundamental class of
$\Delta^{\bullet}$'' gives a retraction of $A_{\ast}$ to
$B^{(0)}_{\ast}$, which is therefore also a quasi-isomorphism.  This
means that the composite map
\begin{equation}\label{eq:82}
\fil_{s}\chzerococycle^{k}\left(S\times \Delta^{\bullet} \right)\to
A_{\ast}\to B^{(0)}_{\ast}
\end{equation}
is a quasi-isomorphism.  But the chain complex $B^{(0)}_{\ast}$
is exactly the chain complex~\eqref{eq:80}, and the map~\eqref{eq:82}
is given by slant product with the fundamental class of the variable
simplex.  This completes the proof.
\end{pf*}

We have used
\begin{lem}\label{thm:33}
For each $m$ and $s$ the simplicial abelian group
\[
[n]\mapsto \fil_{s}\chzero^{m}(S\times \Delta^{n})
\]
is contractible.
\end{lem}

\begin{lem}\label{thm:34}
Let $B^{(i)}_{\ast}$ be the $i^{\text{th}}$ row of \eqref{eq:20}.
The map 
\[
\Delta^{0}\to \Delta^{n}
\]
given by inclusion of any vertex induces an isomorphism of homology groups
\[
H_{\ast}(B^{(0))}_{\ast})\xrightarrow{\approx}{} 
H_{\ast}(B^{(i))}_{\ast}).
\] 
The inclusion of $B^{0}_{\ast}$ into the total complex associated
to~\eqref{eq:20} is therefore a quasi-isomorphism.
\end{lem}

\begin{pf*}{Proof of Lemma~\ref{thm:33}}  It suffices to separately
show that 
\[
[n]\mapsto C^{m}(S\times \Delta^{n}),\qquad
[n]\mapsto C^{m-1}(S\times \Delta^{n};\R)
\]
and 
\[
[n]\mapsto \fil_{s}\Omega^{m}(S\times \Delta^{n})
\]
are acyclic.  These are given by Lemmas~\ref{thm:32} and~\ref{thm:14}
below.
\end{pf*}

\begin{lem}\label{thm:32}
For any abelian group $A$, the simplicial abelian group
\[
[n]\mapsto C^{m}(S\times \Delta^{n};A)
\]
is contractible.
\end{lem}

\begin{pf}
It suffices to construct a contracting homotopy of the associated
chain complex.  The construction makes use of ``extension by zero.''
If $f:X\to Y$ is the inclusion of a subspace, the surjection map
\[
C^{\ast}Y\to C^{\ast}X
\]
has a canonical section $f_!$ given by {\em extension by zero}:  
\[
f_!(c)(z) = \begin{cases}
	 c(z') &\qquad \text{if $z=f\circ z'$} \\
	0 &\qquad\text{otherwise}.
	 \end{cases}
\]
Note that $f_{!}$ is {\em not} a map of cochain complexes.  Using the
formula
\[
\partial_{i}^{\ast} (\partial_{0})_{!}c=
\begin{cases}
c &\qquad i=0 \\
(\partial_{0})_{!}\partial_{i-1}^{\ast} c &\qquad i\ne 0	
\end{cases}
\]
one easily checks that
\[
h={\partial_0}_{!}:C^{n} (M\times\Delta^{k}) \to C^{n} (M\times\Delta^{k+1}) 
\]
is a contracting homotopy.  
\end{pf}

\begin{lem}\label{thm:14}
For each $s\ge 0$ and $t$ chain complex
\[
\fil_{s}\Omega^{t}(S\times\Delta^{\bullet})
\]
is contractible.
\end{lem}

\begin{pf}
We will exhibit a contracting homotopy.  We will use barycentric
coordinates on 
\[
\Delta^{n}=\{(t_{0},\dots,t_{n})\mid 0\le t_{i}\le 1, \sum t_{i}=1 \},
\]
and let $v_{i}$ be the vertex for which $t_{i}=1$.  
Let
\begin{align*}
p_{n}:\Delta^{n}\setminus v_{0} &\to\Delta^{n-1}\\
(t_{0},\dots,t_{n}) &\mapsto
\left(t_{1}/(1-t_{0}),\dots,t_{n}/(1-t_{0})) \right)
\end{align*}
be radial projection from $v_{0}$ to the $0^{\text{th}}$ face.  Then
\begin{align*}
p^{n}\circ\partial^{0} &= \text{Id} \\
p^{n}\circ\partial^{i} &= \partial^{i-1}\circ p^{n-1} \qquad i>0.
\end{align*}
Finally, let 
\[
g:[0,1]\to \R
\]
be a $C^{\infty}$ bump function, vanishing in a neighborhood of $0$,
taking the value $1$ in a neighborhood of $1$.  The map
\begin{align*}
h:\Omega^{t}(S\times\Delta^{m-1}) &\to \Omega^{t}(S\times\Delta^{m}) \\
\omega &\mapsto g(1-t_{0})\,(1\times p^{m})^{\ast}\omega
\end{align*}
is easily checked to preserve $\fil_{s}$ and define a contracting homotopy.
\end{pf}

\begin{pf*}{Proof of Lemma~\ref{thm:34}}.  It suffices to establish
the analogue of Lemma~\ref{thm:34} with the complex $\chzero^{\ast}$
replaced by $C^{\ast}$, $C^{\ast}\otimes\R$, and
$\fil_{s}\Omega^{\ast}$.  The cases of $C^{\ast}$ and
$C^{\ast}\otimes\R$ are immediate from the homotopy invariance of
singular cohomology.  For $\fil_{s}\Omega^{\ast}$, the spectral
sequence associated to the filtration by powers of
$\Omega^{\ast>0}(S)$ shows that the $i^{\text{th}}$-homology group of
\[
\fil_{s}\Omega^{k}(S\times\Delta^{m})_{\text{cl}}
\leftarrow\fil_{s-1}\Omega^{k-1}(S\times\Delta^{m})\leftarrow\dots
\]
is
\[
\begin{cases}H^{k-i}_{\text{DR}}(S;\R)\otimes
H^{0}_{\text{DR}}(\Delta^{m})
&\qquad i<s \\
\Omega^{n}_{\text{cl}}(S)\otimes
H^{0}_{\text{DR}}(\Delta^{m}) &\qquad
i=s,
\end{cases}
\]
so the result follows from the homotopy invariance of deRham cohomology.
\end{pf*}

The results in \S\ref{sec:diff-funct-spectr} require one other
variation on these ideas.  For a compact manifold $S$ and a graded
real vector space $\grv$, let $C_{c}^{*}(S\times\R;\grv)$ denote the
compactly supported cochains, $\Omega_{c}^{*}(S\times\R;\grv)$ the
compactly supported forms, and
$\fil_{s}\Omega_{c}^{\ast}(S\times\Delta^{k}\times\R;\grv)$ the
subspace consisting forms whose Kunneth component on
$\Delta^{k}\times\R$ has degree $\le s$.  The following are
easily verified using the techniques and results of this appendix.

\begin{cor}\label{thm:22}
The map ``slant product with the fundamental class of
$\Delta^{\bullet}$ is a chain homotopy equivalence from the chain
complex underlying
\[
Z^{k}(S\times \Delta^{\bullet};\grv)
\]
to
\[
Z^{k}(S;\grv)\leftarrow
C^{k-1}(S;\grv)\leftarrow\cdots\leftarrow
C^{0}(S;\grv)
\]
\qed
\end{cor}

\begin{cor}\label{thm:16}
``Integration over $\Delta^{\bullet}$'' is a chain homotopy
equivalence between the chain complex underlying the simplicial
abelian group
$\fil_{s}\Omega^{\ast}(S\times\Delta^{\bullet};\grv)^{k}_{\text{cl}}$
and
\[
\Omega^{\ast}(S;\grv)^{k}_{\text{cl}}\xleftarrow{}{}
\Omega^{\ast}(S;\grv)^{k-1}\xleftarrow{}{}\dots\xleftarrow{}{}
\Omega^{\ast}(S;\grv)^{k-s}.
\]
\end{cor}

\begin{cor}\label{thm:7}
The homotopy groups of $\fil_{s}\Omega^{\ast}(S\times
\Delta^{\bullet}; \grv)^{0}_{\text{cl}}$ are given by
\[
\pi_{m}\fil_{s}\Omega^{\ast}(S\times \Delta^{\bullet};
\grv)^{0}_{\text{cl}}=
\begin{cases}
H^{-m}_{\text{DR}}(S;\grv) &\qquad m<s\\
\Omega^{\ast}(S;\grv)^{-s}_{\text{cl}} &\qquad m=s\\
0&\qquad m>s.
\end{cases}
\]
\end{cor}

\begin{cor}\label{thm:23}
Let $[\Delta^{m}\times\R]$ be the product of the fundamental cycle of
$\Delta^{m}$ with the fundamental cycle of the $1$-point
compactification of $\R$.  The map from 
\[
Z^{k}_{c}(S\times\Delta^{\bullet}\times\R;\grv)
\]
to
\[
Z^{k}_{c}(S\times\R;\grv)\leftarrow
C^{k-1}_{c}(S\times\R;\grv)\leftarrow\cdots\leftarrow
C^{0}_{c}(S\times\R;\grv)
\]
sending $f\in Z^{k}(S\times \Delta^{m}\times\R;\grv)$ to $f\slant
[\Delta^{m}\times\R]$ is a chain homotopy equivalence.  
\qed
\end{cor}

\begin{cor}\label{thm:21}
``Integration over $\Delta^{\bullet}\times\R$'' is a chain homotopy
equivalence between chain complex underlying the simplicial abelian
group $\fil_{s}\Omega^{\ast}_{c}(S\times
\Delta^{\bullet}\times\R;\grv)^{k}_{\text{cl}}$ and
\[
\Omega^{\ast}(S;\grv)^{k-1}_{\text{cl}}\xleftarrow{}{}
\Omega^{\ast}(S;\grv)^{k-2}\xleftarrow{}{}\dots\xleftarrow{}{}
\Omega^{\ast}(S;\grv)^{k-s-1}.
\]
\end{cor}

\begin{cor}\label{thm:24}
The maps ``integration over $\R$'' and ``slant product with the
fundamental class of $\bar\R$'' give simplicial homotopy equivalences
\begin{gather*}
\fil_{s}\Omega_{c}^{\ast}\left(S\times\Delta^{\bullet}\times\R;\grv
\right)^{k}_{\text{cl}} \to
\fil_{s-1}\Omega^{\ast}\left(S\times\Delta^{\bullet};\grv
\right)^{k-1}_{\text{cl}}  \\
Z_{c}^{\ast}\left(S\times\Delta^{\bullet}\times\R;\grv
\right)^{k} \to
Z^{\ast}\left(S\times\Delta^{\bullet};\grv
\right)^{k-1}.
\end{gather*}
\end{cor}

\begin{pf}
We will do the case of forms.  The result for singular cocycles is
similar.  It suffices to show that the map of underlying chain
complexes is a chain homotopy equivalence.  Now in the sequence
\begin{multline*}
\fil_{s}\Omega^{\ast}_{c}(S\times\Delta^{\bullet}\times\R;
\grv)^{k}_{\text{cl}} \xrightarrow{\int_{\R}}{}
\fil_{s-1}\Omega^{\ast}(S\times\Delta^{\bullet};
\grv)^{k-1}_{\text{cl}} \\ \xrightarrow{\int_{\Delta^{\bullet}}}{}
\left\{\Omega^{\ast}(S;\grv)^{k-1}_{\text{cl}}\xleftarrow{}{}
\Omega^{\ast}(S;\grv)^{k-2}\xleftarrow{}{}\dots\xleftarrow{}{}
\Omega^{\ast}(S;\grv)^{k-s-1} \right\},
\end{multline*}
The second map and the composition are chain homotopy equivalences by
Corollaries~\ref{thm:16} and~\ref{thm:21} respectively.  It follows
that the first map is also a chain homotopy equivalence.
\end{pf}

\section{Integer Wu-classes for $\spin$-bundles}\label{sec:integer-wu-classes}

Let $S$ be a space, and $V$ and oriented vector bundle over $S$.  Our
aim in this appendix is to associate to a $\spin$-structure on $V$ an
integer lift of the total Wu-class of $V$, and to describe the
dependence of this integer lift on the choice of $\spin$-structure.
We begin by describing 
a family of integer cohomology characteristic classes
$\nu^{\spin}_{k}(V)$ for $\spin$-bundles $V$, whose mod $2$-reduction
are the Wu-classes of the underlying vector bundle.  The classes
$\nu^{\spin}_{k}(V)$ are zero if $k\not\equiv 0$ mod $4$.
The $\nu^{\spin}_{4k}$ satisfy the Cartan formula
\[
\nu^{\spin}_{4n}\left(V\oplus W \right)=
\sum_{i+j=n}\nu^{\spin}_{4k}(V)\, \nu^{\spin}_{4j}(W),
\]
making the total $\spin$ Wu-class
\[
\nu^{\ast}_{t}= 1 + \nu_{4}^{\spin}+\nu_{8}^{\spin}+\dots
\]
into a stable exponential characteristic class.

The characteristic classes $\nu_{k}^{\spin}$, while very natural,
don't quite constitute integral lifts of the Wu-classes.  They are
integer {\em cohomology classes} lifting the Wu-classes, defined by
cocycles up to arbitrary coboundary, while integral lift is
represented by a cocycle chosen up to the coboundary of {\em twice} a
cochain.  We don't know of a natural way of making this choice in
general, so instead we work with an arbitrary one.  Fortunately, the
formula for the effect of a change of $\spin$-structure does not
depend on this choice of lift.  
We have decided to describe the characteristic classes
$\nu_{k}^{\spin}$ partly because they seem to be interesting in their own right,
and partly as a first step toward finding a natural choice of integer
lift of the total Wu-class of $\spin$-bundles.

\subsection{Spin Wu-classes}
\subsubsection{Wu-classes}
Let $V$ be a real $n$-dimensional vector vector bundle over a space
$X$, and write $U\in H^{n}(V,V-\{0 \};\Z/2)$ for the Thom class.  
The Wu-classes $\nu_{i}=\nu_{i}(V)\in H^{i}(X;\Z/2)$ are defined by
the identity
\[
\nu_{t}(V)\,U = \chi(\Sq_{t})(U)
\]
in which 
\[
\Sq_{t}=1+\Sq_{1}+\Sq_{2}+\dots
\]
is the total mod $2$ Steenrod operation, $\chi$ the canonical
anti-automorphism (antipode) of the Steenrod algebra, 
\[
\nu_{t}=1+\nu_{1}+\nu_{2}+\dots
\]
the total Wu-class of $V$.  One checks from the definition that the
Wu-classes satisfy

\begin{textList}
\item (Cartan formula) 
\[
\nu_{t}(V\oplus W)= \nu_{t}(V)\,\nu_{t}(W);
\]

\item If $V$ is a real line bundle with $w_{1}(V)=\alpha$ then 
\[
\nu_{t}(V)= \sum \alpha^{2^{n}-1}.
\]
\end{textList}
Indeed, the first part follows from the fact the $\chi$ is compatible
with the coproduct, and the second part follows from the formula
\[
\chi(\Sq_{t})(x)= \sum x^{{2^{n}}}.
\]
These two properties characterize the Wu-classes as the mod $2$ cohomology stable
exponential characteristic class with characteristic series 
\[
1+x+x^{3}+\dots+x^{2^{n}-1}+\dots\in\Z/2\LL x\RR = H^{\ast}(\rp;\Z/2).
\]

\subsubsection{Complex Wu classes}
\label{sec:complex-wu-classes}

Now consider complex analogue of the above.  We define
\[
\nu_{t}^{\C}=1+\nu_{1}^{\C} + \nu_{2}^{\C}\dots
\in H^{\ast}(\bu;\Z)
\]
to be the stable exponential characteristic class with characteristic series
\[
1+x+\dots+ x^{2^{n}-1}+\dots\in \Z\LL x\RR=H^{\ast}(\cp;\Z).
\]
So for a complex line bundle $L$, with first Chern class $x$,
\[
\nu_{t}^{\C}(L)=\sum_{n\ge 0}x^{2^{n}-1}.
\]
One easily checks that
\[
\nu^{\C}(L) \equiv \nu(L)\mod 2.
\]
so the classes $\nu_{t}^{\C}$ is a stable exponential integer characteristic
class lifting the total Wu class in case the vector bundle $V$ has a
complex structure.

\subsubsection{Wu-classes for Spin bundles}
\label{sec:wu-classes-spin}

Because of the presence of torsion in $H^{\ast}(\bspin;\Z)$,
specifying a stable exponential characteristic class 
\[
\nu_{t}^{\spin}=1+\nu_{4}^{\spin}+\nu_{8}^{\spin}+\dots\in H^{\ast}(\bspin;\Z)
\]
bundles is a little subtle.  Things are simplified somewhat by the
fact that the torsion in $H^{\ast}(\bspin)$ is all of order $2$.  
It implies that
\[
\begin{CD}
H^{\ast}(\bspin;\Z) @>>> H^{\ast}(\bspin;\Z)/\text{torsion} \\
@VVV @VVV \\
H^{\ast}(\bspin;\Z)\otimes\Z/2 @>>> (H^{\ast}(\bspin;\Z)/\text{torsion})\otimes\Z/2
\end{CD}
\]
is a pullback square and that $H^{\ast}(\bspin;\Z)\otimes\Z/2$ is 
just the kernel of 
\[
\Sq^{1}: H^{\ast}(\bspin;\Z/2)\to H^{\ast+1}(\bspin;\Z/2).
\]
One also knows that $H^{\ast}(\bspin;\Z)/\text{torsion}$
is a summand of $H^{\ast}(\bsu;\Z)$.  So to specify an element of
$H^{\ast}(\bspin;\Z)$ one needs to give $a\in H^{\ast}(\bspin;\Z/2)$
and $b\in H^{\ast}(\bsu;\Z)$ with the properties
\begin{textList}
\item\label{item-sq-1} $\Sq^{1}(a)=0$;
\item The image of $b$ in $H^{\ast}(\bsu;\Q)$ is in the image of $H^{\ast}(\bspin;\Q)$;
\item $b\equiv a\mod 2\in H^{\ast}(\bsu;\Z/2)=H^{\ast}(\bsu;\Z)\otimes\Z/2$.
\end{textList}

We take $a$ to be (the restriction of) the total Wu-class $\nu_{t}\in
H^{\ast}(\bspin;\Z/2)$.  Property~\textItemref{item-sq-1} is then given by
the following (well-known) result.

\begin{lem}\label{thm:3}
If $V$ is a $\spin$-bundle then $\nu_{k}(V)=0$ if $n\not\equiv 0\mod
4$, and $\Sq^{1}\nu_{4k}=0$.
\end{lem}

\begin{pf}
Since $V$ is a spin bundle $\Sq^{k}U=w_{k}U=0$ for $k\le 3$.  Using
the Adem relations
\begin{align*}
\Sq^{1}\Sq^{2k} &=\Sq^{2k+1} \\
\Sq^{2}\Sq^{4k} &=\Sq^{4k+2} +\Sq^{4k+1}\Sq^{1}\\
&=\Sq^{4k+2} +\Sq^{1}\Sq^{4k}\Sq^{1}\\
\end{align*}
we calculate
\begin{align*}
\nu_{2k+1}U=\chi(\Sq^{2k+1})U &=\chi(\Sq^{1}\Sq^{2k})\\
&=\chi(Sq^{2k})\chi(Sq^{1})U \\
&=\chi(Sq^{2k})Sq^{1}U=0.
\end{align*}
and 
\begin{align*}
\nu_{4k+2}U=\chi(\Sq^{4k+2})U &=\chi(\Sq^{2}\Sq^{4k+1})U+\chi(\Sq^{1}\Sq^{4k}\Sq^{1})\\
&= \chi(Sq^{4k})\chi(Sq^{2})U + \chi(\Sq^{1})\chi(\Sq^{4k})\chi(\Sq^{1})U \\
&= \chi(Sq^{4k})\,Sq^{2}\,U + \Sq^{1}\chi(\Sq^{4k})\,\Sq^{1}\,U \\
&= 0.
\end{align*}
This gives the first assertion.  For the second,
use the Adem relation
\[
\Sq^{2}\Sq^{4k-1}=\Sq^{4k}\Sq^{1}
\]
and calculate
\begin{align*}
\Sq^{1}\nu_{4k}U = \Sq^{1}\chi(\Sq^{4k})U &= \chi(\Sq^{4k}\Sq^{1})U \\
&= \chi(\Sq^{2}\Sq^{4k-1})U \\
&= \chi(\Sq^{4k-1})\Sq^{2}U =0.\\
\end{align*}
\end{pf}

We would like to take for $b\in H^{\ast}(\bsu;\Z)$ the restriction of
$\nu_{t}^{\C}$.  Unfortunately that class is not in the image of
$H^{\ast}(\bspin)$.  Indeed, $\nu_{t}^{\C}$ is the stable exponential characteristic
class with characteristic series
\[
f(x)= 1+x+x^{3}+\dots+x^{2^{n}-1}+\dots.
\]
which is not symmetric.  But if we symmetrize $f$ by setting
\[
g(x)= \sqrt{f(x)f(-x)}= 1 - \frac{x^2}{2} - \frac{9\,x^4}{8} - \frac{17\,x^6}{16}+\dots
\in \Z[\tfrac12]\LL x\RR
\]
then $g(x)$, being even, will be the characteristic series of a stable exponential
characteristic class
\[
\chi\in H^{\ast}(\bso;\Z[\tfrac12]).
\]

\begin{lem}\label{thm:11}
The restriction of $\chi$ to $H^{\ast}(\bsu;\Z[\tfrac12])$ lies in $H^{\ast}(\bsu;\Z)$.
\end{lem}

\begin{pf}
We will use the result
of~\cite{hopkins95:_topol_witten,ando01:_ellip_witten} that a stable
exponential characteristic class for $\su$-bundles is determined by
its value on
\[
(1-L_{1})(1-L_{2})
\]
which can be any power series 
\[
h(x,y)\in H^{\ast}(\cp\times\cp)
\]
satisfying
\begin{textList}
\item $h(0,0)=1$;
\item $h(x,y)=h(y,x)$;
\item $h(y,z)h(x,y+z)=h(x+y,z)h(x,y)$.
\end{textList}
We will call $h$ the {\em $\su$-characteristic series}.  The
$\su$-characteristic series of $\chi$ is
\[
\delta g(x,y) :=  \frac{g(x+y)}{g(x)g(y)}.
\]
The lemma  will follow once we show that $\delta
g(x,y)$ has integer coefficients.  

Now 
\[
g(x,y)^{2}= \delta f(x,y)\, \delta f(-x,-y).
\]
Since  the power series $\sqrt{1+4 x}$ has integer coefficients, it
suffices to show that 
\[
\delta f(x,y)\, \delta f(-x,-y)
\]
is a square in $\Z/4\LL x\RR$.  For this, start with 
\begin{align*}
f(x)-f(-x) &= 2 x f(x^{2}) &\\
&\equiv 2 x \,f(x)^{2} &\mod 4\,\Z\LL x\RR \\
&\equiv 2 x \,f(x)\,f(-x) &\mod 4\,\Z\LL x\RR
\end{align*}
to conclude that 
\[
f(x)\equiv f(-x)\left(1+2 x f(x) \right)\mod 4\,\Z\LL x\RR.
\]
Write 
\[
e(x)=1+2 x f(x) = 1 + 2 \sum_{n\ge 0} x^{2^{n}}.
\]
Then 
\[
e(x+y)\equiv e(x)e(y) \mod 4\,\Z\LL x\RR
\]
(in fact $e(x)\equiv e^{2x}$ mod $4\,\Z\LL x\RR$).  This implies
\[
\delta f(x,y) \equiv \delta f(-x,-y) \mod 4\,\Z\LL x\RR,
\]
and so 
\begin{equation}\label{eq:2}
\delta f(x,y)\, \delta f(-x,-y) = \delta f(x,y)^{2} \mod 4\,\Z\LL x\RR.
\end{equation}
\end{pf}

Continuing with the construction of $\nu_{t}^{\spin}$, we now take the
class $b$ to be the stable exponential characteristic class with
$\su$-characteristic series $g(x)$.    It remains to check that
$b\equiv a\in H^{\ast}(\bsu;\Z/2)$, or, in terms of
$\su$-characteristic series, that
\[
\delta g(x,y)\equiv \delta f(x,y) \mod 2\, \Z\LL x\RR.
\]
But this identity obviously holds after squaring both sides, which suffices, 
since $\Z/2\LL x\RR$ is an integral domain over $\Z/2$.

To summarize, we have constructed a stable exponential characteristic
class $\nu_{t}^{\spin}$ for $\spin$ bundles, with values in integer
cohomology, whose mod $2$-reduction is the total Wu-class.
Rationally, in terms of Pontryagin classes, the first few are
\begin{align*}
\nu^{\spin}_{4} &= -\frac{p_{1}}{2} \\
\nu^{\spin}_{8} &= \frac{20\,p_{2}-9\,p_{1}^{2}}{8} \\
\nu^{\spin}_{12} &= \frac{-80\, p_{3} + 60\,p_{1}p_{2}-17\,p_{1}^{3}}{16}.
\end{align*}
It frequently comes up in geometric applications that one wishes to
express the value of a stable exponential characteristic class of the normal bundle
in terms characteristic classes of the tangent bundle.  We therefore also record
\begin{align*}
\nu^{\spin}_{4}(-T) &= \frac{p_{1}}{2} \\
\nu^{\spin}_{8}(-T) &= \frac{-20\, p_{2}+11\,p_{1}^{2}}{8} \\
\nu^{\spin}_{12}(-T) &= \frac{80\, p_3-100\, p_2 p_1+37\, p_1^3}{16},
\end{align*}
where in these expressions the Pontryagin classes are those of $T$ (and not $-T$).

\subsection{Integral Wu-structures and change of $\spin$-structure}
\label{sec:change-spin-struct}

Suppose $S$ is a space, and $\nu\in Z^{k}(S;\Z/2)$ a cocycle.  Recall
from \S\ref{subsec-wustruct} that the category $\hcat_{\nu}^{k}(S)$ of {\em
integer lifts} of $\nu$ is the category whose objects are cocycles
$x\in Z^{k}(S;\Z)$ whose reduction modulo $2$ is $\nu$.  A morphism
from $x$ to $y$ in $\hcat_{\nu}^{k}(S)$ is a $(k-1)$-cochain $c\in
C^{k-1}(S;\Z)$ with the property that
\[
2\,\delta c=y-x.
\]
We identify two morphisms $c$ and $c'$ if they differ by a
$(k-1)$-cocycle.  The set $H_{\nu}^{k}(S)$ of isomorphism classes of
objects in $\hcat_{\nu}^{k}(S)$ is a torsor for $H^{k}(S;\Z)$.  We
write the action of $b\in H^{k}(S)$ on $x\in H^{k}_{\nu}(S)$ as
\[
x\mapsto x + (2)\,b.
\]

There is a differential analogue of this notion when $S$ is a
manifold.   We define a {\em differential integral lift} of $\nu$ to
be a differential cocycle 
\[
x=(c,h,\omega)\in \chcocycle{k}^{k}(S)
\]
with the property that $c\equiv \nu$ modulo $2$.  The set of
differential integral lifts of $\nu$ forms a category
$\chcat_{\nu}^{k}(S)$ in which a morphism from $x$ to $y$ is an element
$c\in\chcochain{k}^{k-1}/\chcocycle{k}^{k-1}$ with the property  that
\[
2\,\delta c=y-x.
\]
The set of isomorphism classes in $\chcat_{\nu}^{k}(S)$ is denoted
$\chh{k}^{k}_{\nu}(S).$  It is a torsor for the differential
cohomology group $\chh{k}^{k}(S)$.   As above, we
write the action of $b\in \chh{k}^{k}(S)$ on $x\in \chh{k}^{k}_{\nu}(S)$ as
\[
x\mapsto x + (2)\,b.
\]

Suppose now that $V$ is an oriented vector bundle over $S$, and $\nu$
a cocycle representing the total Wu-class of $V$.  By
Proposition~\ref{thm:3}, if $w_{2}(V)=0$ the category
$H_{\nu}^{\ast}(S)$ is non-empty, and an integer lift of the total Wu
class can be associated to a $\spin$-structure.  We wish to describe
the effect of a change of $\spin$-structure on these integral lifts.
In case $S$ is a manifold, and $V$ is equipped with a connection, we
can associate a differential integral lift of the total Wu-class to
the $\spin$-structure.  We are also interested in the effect of a
change of $\spin$-structure on these differential integral lifts.

Choose cocycles $z_{2k}\in Z^{2k}(\bso;\Z/2)$ representing the
universal Wu-classes.  To associate an integer lift of the total
Wu-class to a $\spin$-structure we must choose cocycles $\zbar_{2k}\in
Z^{2k}(\bspin;\Z)$ which reduce modulo $2$ to the restriction of the
$z_{2k}$ cocycles to $\bspin$.  The classes $\zbar$ could be taken to
represent the classes $\nu^{\spin}$ constructed in the previous
section, though this is not necessary.  We represent the cocycles
$z_{2k}$ and $\zbar_{2k}$ by maps to Eilenberg-MacLane spaces,
resulting in a diagram
\[
\begin{CD}
\bspin @> \zbar=\prod\zbar_{2k}>> \prod_{k\ge 2}K(\Z,2k) \\
@VVV @VVV \\
\bso @>> z=\prod z_{2k}> \prod_{k\ge 1} K(\Z/2,2k)
\end{CD}
\]
Choose a map $S\to\bso$ classifying $V$,
$s_{0}:S\to \bspin$ a lift
corresponding to a $\spin$-structure on $V$, and $\alpha:S\to
K(\Z/2,1)$ a $1$-cocycle on $S$.  Let $s_{1}$ be the $\spin$-structure
on $V$ gotten by changing $s_{0}$ by $\alpha$.  We want a formula
relating $\zbar(s_{0})$  and $\zbar(s_{1})$, or, more precisely, a
formula for the cohomology class represented by 
\[
\frac{\zbar(s_{1})-\zbar(s_{1})}{2}.
\]

We first translate this problem into one involving the long exact
sequence of homotopy groups of a fibration.  Write $P_{0}\to
B_{0}$ for the fibration
\begin{equation}\label{eq:3}
\bspin^{S}\to\bso^{S}
\end{equation}
and $P_{1}\to B_{1}$ for
\begin{equation}\label{eq:5}
\left(\prod_{k\ge 1}K(\Z,2k) \right)^{S}
\to
\left(\prod_{k\ge 1}K(\Z/2,2k) \right)^{S}.
\end{equation}
Then~\eqref{eq:3} and~\eqref{eq:5} are principal bundles with
structure groups
\[
G_{0}=K(\Z/2,1)^{S}\quad\text{and}\quad
\left(\prod_{k\ge 1}K(\Z,2k) \right)^{S}
\]
respectively.  We give $P_{0}$ the basepoint corresponding to $s_{0}$,
and $B_{0}$ the basepoint corresponding to $V$, and $P_{1}$ and
$B_{1}$ the basepoints corresponding to $\zbar(s_{0})$ and $z(V)$.
These choices identify the fibers over $s_{0}$ and
$\zbar(s_{0})$ with $G_{0}$ and $G_{1}$, and lead to a map of pointed fiber sequences
\[
\begin{CD}
G_{0} @>>> P_{0} @>>> B_{0}\\
@VVV @VVV @VVV \\
G_{1} @>>> P_{1} @>>> B_{1}
\end{CD}
\]
derived from $(\zbar,z)$.  Our problem is to describe the map
\[
\pi_{0}G_{0}=H^{1}(S;\Z/2)\to
\pi_{0}G_{1}=\prod_{k\ge 1}H^{2k}(S;\Z).
\]
In general this can be difficult; but for elements in the image of
$\pi_{1}B_{0}\to\pi_{0}G_{0}$ the answer is given by the composite
$\pi_{1}B_{0}\to\pi_{1}B_{1}\to\pi_{0}G_{0}$.
In~\cite{milnor63:_spin} (stated in the language of
$\spin$-structures) Milnor shows that $\pi_{1}B_{0}\to \pi_{0}G_{0}$
is surjective.  So in our case, every element of $\pi_{0}G_{0}$ is in
the image of $\pi_{1}B_{0}$.  This leads to our desired formula.

We now translate this discussion back into the language of
$\spin$-structures and characteristic classes.  The element of
$\pi_{0}G_{0}=H^{1}(S;\Z/2)$ is the cohomology class represented by
$\alpha$.  To lift this to an element of $\pi_{1}B_{0}$ is equivalent
to finding an oriented stable vector bundle $W$ over 
\[
S\times S^{1}
\]
satisfying
\begin{gather*}
W\vert_{S\times\{1 \}}=V \\
w_{2}(W)=\alpha\cdot U
\end{gather*}
where $U\in H^{1}(S^{1};\Z/2)$ is the generator.  This is easily
done.  We take
\[
W=V\oplus (1-L_{\alpha})\otimes(1-H)
\]
in which $L_{\alpha}$ is the real line bundle whose first
Stiefel-Whitney class is represented by $\alpha$, and $H$ is the
non-trivial real line bundle over $S^{1}$.   (This proves Milnor's result
on the surjectivity of $\pi_{1}B_{0}\to\pi_{0}G_{0}$.)   The image of
this class in $\pi_{1}B_{1}$ is the total Wu-class of $W$, which is 
\begin{equation}\label{eq:6}
\nu_{t}(V)\,\frac{\nu(L_{\alpha}\otimes H)}{\nu(L_{\alpha})\nu(H)}.
\end{equation}
Writing
\begin{gather*}
\nu_{t}(V)= 1+\nu_{1}+\nu_{2}+\dots \\
w_{1}(L_{\alpha})=\alpha\\
w_{1}(H)=\epsilon
\end{gather*}
and remembering that $\epsilon^{2}=0$, one evaluates~\eqref{eq:6} to be 
\begin{align*}
\nu_{t}(V)\,\frac{\sum_{n\ge 0}
(\alpha+\epsilon)^{2^{n}-1}}{(1+\epsilon)\,\sum_{n\ge 0}
\alpha^{2^{n}-1}} &=
\nu_{t}(V)\frac{\sum \alpha^{2^{n}-1}+\epsilon\,\left(\sum
\alpha^{2^{n}-1} \right)^{2}}{(1+\epsilon)\,\sum \alpha^{2^{n}-1}} \\
&=
\nu_{t}(V)\, \frac{1+\epsilon\,\sum \alpha^{2^{n}-1}}{(1+\epsilon)} \\
&= \nu_{t}(V)\,\left( 1+\epsilon\,\sum_{n\ge 1} \alpha^{2^{n}-1}\right).
\end{align*}
The image of this class under $\pi_{1}B_{1}\to\pi_{0}G_{1}$ is
computed by taking the slant product with the fundamental class of
$S^{1}$ and then applying the Bockstein homomorphism.  This leads to 
\begin{equation}\label{eq:8}
\beta\left(\nu_{t}(V)\,\sum_{n\ge 1}\alpha^{2^{n}-1} \right)
=\nubar_{t}(V)\,\beta \left(\sum_{n\ge 1}\alpha^{2^{n}-1}\right) \in \prod H^{2n}(S;\Z). 
\end{equation}
This proves

\begin{prop}\label{thm:2}
Suppose that $V$ is a vector bundle over a space $S$, $s$ a
$\spin$-structure on $V$,  $\alpha\in H^{1}(S;\Z/2)$.   Let $\nubar_{t}$
be any integer lift of the restriction of the total Wu class to
$\bspin$.  Then
\begin{align*}
\nubar_{t}(s+\alpha) &=\nubar(s)+(2)\,\beta \left(\nu_{t}(V)\sum_{n\ge
1}\alpha^{2^{n}-1}\right) \\
&= \nubar(s)+(2)\,\nu_{t}(V)\sum_{n\ge
1}\beta(\alpha^{2^{n}-1})
\end{align*}
\end{prop}
In this expression, the factor ``$(2)$'' is formal.  It is written to
indicate the action of $\prod H^{2n}(S;\Z)$ on the set of integral
Wu-structures.  The number $2$ serves as a reminder that the action of
$x$ on the cohomology class $w$ underlying an integral Wu-structure is
$w\mapsto w+2 x$.

An analogous discussion, using differential cohomology, leads to the
following 

\begin{prop}\label{thm:4}
Let $S$ be a manifold and $s:S\to \bspin$, represent a
$\spin$-structure on a stable oriented vector bundle $V$ with
connection $\nabla$.  Write $\nubar(s,\nabla)\in \prod_{k\ge
0}\chcat^{2k}_{\nu}(S)$ for the twisted differential cocycle associated to
$s$, $\nabla$ and the cocycle $\nubar$ (see
\S\ref{sec:char-class}). If $\alpha\in Z^{1}(S;\Z/2)$, then
\[
\nubar(s+\alpha,\nabla)=
\nubar(s,\nabla) + (2)\, \beta\left(\nu(V)\sum_{k\ge 1}\alpha^{2^{k}-1} \right)
\]
where, again, the factor ``$(2)$'' is formal, indicating the action of
$\prod \chh{2k}^{2k}(S)$ on the set of isomorphism classes of
$\nu$-twisted differential cocycles, and $\beta$ denotes the map
\[
\prod H^{2k-1}(S;\Z/2)\to \prod H^{2k-1}(S;\rz)\subset
\prod \chh{2k}^{2k}(S).
\]
\end{prop}


\begin{thebibliography}{10}

\bibitem{Ad:SHGH}
J.~F. Adams, \emph{Stable homotopy and generalised homology}, University of
  Chicago Press, Chicago, 1974.

\bibitem{alvarez85:_cohom}
Orlando Alvarez, \emph{Cohomology and field theory}, Symposium on anomalies,
  geometry, topology (Chicago, Ill., 1985), World Sci. Publishing, Singapore,
  1985, pp.~3--21. \MR{87j:81081}

\bibitem{alvarez01:_beyon_ellip_genus}
Orlando Alvarez and I.~M. Singer, \emph{Beyond the {E}lliptic {G}enus},
  hep-th/0104199, 2001.

\bibitem{anderson69:_univer_k}
D.~W. Anderson, \emph{Universal coefficient theorems for $k$-theory},
  mimeographed notes, Univ. California, Berkeley, Calif., 1969.

\bibitem{anderson66:_su_ko_kervair}
D.~W. Anderson, E.~H. Brown, Jr., and F.~P. Peterson, \emph{${\rm
  {s}{u}}$-cobordism, ${\rm {k}{o}}$-characteristic numbers, and the {K}ervaire
  invariant}, Ann. of Math. (2) \textbf{83} (1966), 54--67. \MR{32 \#6470}

\bibitem{ando01:_ellip_witten}
M.~Ando, M.~J. Hopkins, and N.~P. Strickland, \emph{Elliptic spectra, the
  {W}itten genus and the theorem of the cube}, Invent. Math. \textbf{146}
  (2001), no.~3, 595--687. \MR{2002g:55009}

\bibitem{aps3}
M.~F. Atiyah, V.~K. Patodi, and I.~M. Singer, \emph{Spectral asymmetry and
  {R}iemannian geometry. {I}{I}{I}}, Math. Proc. Cambridge Philos. Soc.
  \textbf{79} (1976), no.~1, 71--99. \MR{53 \#1655c}

\bibitem{atiyah71:_rieman}
Michael~F. Atiyah, \emph{Riemann surfaces and spin structures}, Ann. Sci.
  \'Ecole Norm. Sup. (4) \textbf{4} (1971), 47--62.

\bibitem{boardman73:_homot}
J.~M. Boardman and R.~M. Vogt, \emph{Homotopy invariant algebraic structures on
  topological spaces}, Springer-Verlag, Berlin, 1973, Lecture Notes in
  Mathematics, Vol. 347. \MR{54 \#8623a}

\bibitem{browder69:_kervair}
W.~Browder, \emph{The {Kervaire} invariant of framed manifolds and its
  generalization}, Annals of Mathematics \textbf{90} (1969), 157--186.

\bibitem{BC}
E.~H. Brown and M.~Comenetz, \emph{Pontrjagin duality for generalized homology
  and cohomology theories}, American Journal of Mathematics \textbf{98} (1976),
  1--27.

\bibitem{brown65:_note_kervair}
Edgar~H. Brown, Jr., \emph{Note on an invariant of {K}ervaire}, Michigan Math.
  J. \textbf{12} (1965), 23--24. \MR{32 \#8351}

\bibitem{brown71:_kervair}
\bysame, \emph{The {K}ervaire invariant of a manifold}, Algebraic topology
  (Proc. Sympos. Pure Math., Vol. XXII, Univ. Wisconsin, Madison, Wis., 1970),
  Amer. Math. Soc., Providence, R. I., 1971, pp.~65--71. \MR{48 \#1249}

\bibitem{brown72:kervair}
\bysame, \emph{Generalizations of the {K}ervaire invariant}, Ann. of Math. (2)
  \textbf{95} (1972), 368--383. \MR{45 \#2719}

\bibitem{brown65:_kervair}
Edgar~H. Brown, Jr. and Franklin~P. Peterson, \emph{The {K}ervaire invariant of
  $(8k+2)$-manifolds}, Bull. Amer. Math. Soc. \textbf{71} (1965), 190--193.
  \MR{30 \#584}

\bibitem{brown66:_kervair}
\bysame, \emph{The {K}ervaire invariant of $(8k+2)$-manifolds}, Amer. J. Math.
  \textbf{88} (1966), 815--826. \MR{34 \#6786}

\bibitem{brylinski93:_loop}
Jean-Luc Brylinski, \emph{Loop spaces, characteristic classes and geometric
  quantization}, Birkh\"auser Boston Inc., Boston, MA, 1993. \MR{94b:57030}

\bibitem{cerf61:_topol}
Jean Cerf, \emph{Topologie de certains espaces de plongements}, Bull. Soc.
  Math. France \textbf{89} (1961), 227--380. \MR{25 \#3543}

\bibitem{cerf62:_la_diff_s}
\bysame, \emph{La nullit\'e de $\pi \sb{0}\,({\rm {d}iff}\,{S}\sp{3})$.
  {T}h\'eor\`emes de fibration des espaces de plongements. {A}pplications},
  S\'eminaire Henri Cartan, 1962/63, Exp. 8, Secr\'etariat math\'ematique,
  Paris, 1962/1963, p.~13. \MR{28 \#3444}

\bibitem{cheeger85:_differ}
Jeff Cheeger and James Simons, \emph{Differential characters and geometric
  invariants}, Geometry and topology (College Park, Md., 1983/84), Lecture
  Notes in Math., vol. 1167, Springer, Berlin, 1985, pp.~50--80.

\bibitem{chern74:_charac}
Shiing~Shen Chern and James Simons, \emph{Characteristic forms and geometric
  invariants}, Ann. of Math. (2) \textbf{99} (1974), 48--69. \MR{50 \#5811}

\bibitem{curtis71:_simpl}
Edward~B. Curtis, \emph{Simplicial homotopy theory}, Advances in Math.
  \textbf{6} (1971), 107--209 (1971). \MR{43 \#5529}

\bibitem{deligne87:_det}
P.~Deligne, \emph{Le d\'eterminant de la cohomologie}, Current trends in
  arithmetical algebraic geometry (Arcata, Calif., 1985), Amer. Math. Soc.,
  Providence, RI, 1987, pp.~93--177. \MR{89b:32038}

\bibitem{deligne71:_theor_hodge}
Pierre Deligne, \emph{Th\'eorie de {H}odge. {I}{I}}, Inst. Hautes \'Etudes Sci.
  Publ. Math. (1971), no.~40, 5--57. \MR{58 \#16653a}

\bibitem{dolan98}
Louise Dolan and Chiara~R. Nappi, \emph{A modular invariant partition function
  for the fivebrane}, Nuclear Phys. B \textbf{530} (1998), no.~3, 683--700,
  hep-th/9806016. \MR{99k:81231}

\bibitem{douady61:_variet}
Adrien Douady, \emph{Vari\'et\'es \`a bord anguleux et voisinages tubulaires},
  S\'eminaire Henri Cartan, 1961/62, Exp. 1, Secr\'etariat math\'ematique,
  Paris, 1961/1962, p.~11. \MR{28 \#3435}

\bibitem{elmendorf95:_moder}
A.~D. Elmendorf, I.~K{\v{r}}{\'\i}{\v{z}}, M.~A. Mandell, and J.~P. May,
  \emph{Modern foundations for stable homotopy theory}, Handbook of algebraic
  topology, North-Holland, Amsterdam, 1995, pp.~213--253. \MR{97d:55016}

\bibitem{elmendorf97:_rings}
A.~D. Elmendorf, I.~Kriz, M.~A. Mandell, and J.~P. May, \emph{Rings, modules,
  and algebras in stable homotopy theory}, American Mathematical Society,
  Providence, RI, 1997, With an appendix by M. Cole. \MR{97h:55006}

\bibitem{freed00:_dirac_charg_quant_gener_differ_cohom}
Daniel~S. Freed, \emph{Dirac charge quantization and generalized differential
  cohomology}, hep-th/0011220, 2000.

\bibitem{freed:_ramon_ramon_k}
Daniel~S. Freed and Michael Hopkins, \emph{On {R}amond-{R}amond fields and
  ${K}$-theory}, J. High Energy Phys. (2000), no.~5, Paper 44, 14. \MR{1 769
  477}

\bibitem{goerss99:_simpl}
Paul~G. Goerss and John~F. Jardine, \emph{Simplicial homotopy theory},
  Birkh\"auser Verlag, Basel, 1999. \MR{2001d:55012}

\bibitem{harris82:_theta}
Joe Harris, \emph{Theta-characteristics on algebraic curves}, Trans. Amer.
  Math. Soc. \textbf{271} (1982), no.~2, 611--638.

\bibitem{harvey:_from_spark_grund}
Reese Harvey and H.~Blaine~Lawson Jr, \emph{{F}rom sparks to
  grundles--differential characters}, arXiv:math.DG/0306193.

\bibitem{hirsch94:_differ}
Morris~W. Hirsch, \emph{Differential topology}, Springer-Verlag, New York,
  1994, Corrected reprint of the 1976 original. \MR{96c:57001}

\bibitem{HopGoerss:Mult}
M.~J. Hopkins and P.~Goerss, \emph{Multiplicative stable homotopy theory}, book
  in preparation.

\bibitem{hopkins95:_topol_witten}
Michael~J. Hopkins, \emph{Topological modular forms, the {W}itten genus, and
  the theorem of the cube}, Proceedings of the International Congress of
  Mathematicians, Vol.\ 1, 2 (Z\"urich, 1994) (Basel), Birkh\"auser, 1995,
  pp.~554--565. \MR{97i:11043}

\bibitem{kervaire62:_la_pontr}
Michel Kervaire, \emph{La m\'ethode de {P}ontryagin pour la classification des
  applications sur une sph\`ere}, Topologia Differenziale (Centro Internaz.
  Mat. Estivo, deg. 1 Ciclo, Urbino, 1962), Lezione 3, Edizioni Cremonese,
  Rome, 1962, p.~13.

\bibitem{kervaire60}
Michel~A. Kervaire, \emph{A manifold which does not admit any differentiable
  structure}, Comment. Math. Helv. \textbf{34} (1960), 257--270.

\bibitem{laures00}
Gerd Laures, \emph{On cobordism of manifolds with corners}, Trans. Amer. Math.
  Soc. \textbf{352} (2000), no.~12, 5667--5688 (electronic). \MR{1 781 277}

\bibitem{LMayS}
L.~G. Lewis, J.~P. May, and M.~Steinberger, \emph{Equivariant stable homotopy
  theory}, Lecture Notes in Mathematics, vol. 1213, Springer--Verlag, New York,
  1986.

\bibitem{lott:_higher}
John Lott, \emph{Higher-degree analogs of the determinant line bundle},
  http://www.math.lsa.umich.edu/$\sim$lott.

\bibitem{lott94:_r_z}
\bysame, \emph{${\bf {r}}/{\bf {z}}$ index theory}, Comm. Anal. Geom.
  \textbf{2} (1994), no.~2, 279--311. \MR{95j:58166}

\bibitem{mathai86:_super_thom}
Varghese Mathai and Daniel Quillen, \emph{Superconnections, {T}hom classes, and
  equivariant differential forms}, Topology \textbf{25} (1986), no.~1, 85--110.
  \MR{87k:58006}

\bibitem{may72}
J.~P. May, \emph{The geometry of iterated loop spaces}, Springer-Verlag,
  Berlin, 1972, Lectures Notes in Mathematics, Vol. 271. \MR{54 \#8623b}

\bibitem{may92:_simpl}
J.~Peter May, \emph{Simplicial objects in algebraic topology}, University of
  Chicago Press, Chicago, IL, 1992, Reprint of the 1967 original.
  \MR{93m:55025}

\bibitem{melrose:_differ_analy_manif_corner}
Richard Melrose, \emph{Differential analysis on manifolds with corners},
  http://www-math.mit.edu/$\sim$rbm/book.html, in preparation.

\bibitem{milgram74:_surger}
R.~James Milgram, \emph{Surgery with coefficients}, Ann. of Math. (2)
  \textbf{100} (1974), 194--248. \MR{50 \#14801}

\bibitem{milnor63:_spin}
J.~Milnor, \emph{Spin structures on manifolds}, Enseignement Math. (2)
  \textbf{9} (1963), 198--203. \MR{MR0157388 (28 \#622)}

\bibitem{milnor61}
John Milnor, \emph{A procedure for killing homotopy groups of differentiable
  manifolds.}, Proc. Sympos. Pure Math., Vol. III, American Mathematical
  Society, Providence, R.I, 1961, pp.~39--55.

\bibitem{milnor73:_symmet}
John Milnor and Dale Husemoller, \emph{Symmetric bilinear forms},
  Springer-Verlag, New York, 1973, Ergebnisse der Mathematik und ihrer
  Grenzgebiete, Band 73. \MR{58 \#22129}

\bibitem{morgan74}
John~W. Morgan and Dennis~P. Sullivan, \emph{The transversality characteristic
  class and linking cycles in surgery theory}, Ann. of Math. (2) \textbf{99}
  (1974), 463--544. \MR{50 \#3240}

\bibitem{mumford71:_theta}
David Mumford, \emph{Theta characteristics of an algebraic curve}, Ann. Sci.
  \'Ecole Norm. Sup. (4) \textbf{4} (1971), 181--192.

\bibitem{narasimhan61:_exist}
M.~S. Narasimhan and S.~Ramanan, \emph{Existence of universal connections},
  Amer. J. Math. \textbf{83} (1961), 563--572. \MR{24 \#A3597}

\bibitem{narasimhan63:_exist}
\bysame, \emph{Existence of universal connections. {I}{I}}, Amer. J. Math.
  \textbf{85} (1963), 223--231. \MR{27 \#1904}

\bibitem{pontryagin38_2}
L.~S. Pontryagin, \emph{A classification of continuous transformations of a
  complex into a sphere. 1}, Doklady Akad. Nauk SSSR (N.S.) \textbf{19} (1938),
  361--363.

\bibitem{pontryagin38_1}
\bysame, \emph{A classification of continuous transformations of a complex into
  a sphere. 2}, Doklady Akad. Nauk SSSR (N.S.) \textbf{19} (1938), 147--149.

\bibitem{pontryagin50:_homot}
\bysame, \emph{Homotopy classification of mappings of an $(n+2)$-dimensional
  sphere on an $n$-dimensional one}, Doklady Akad. Nauk SSSR (N.S.) \textbf{19}
  (1950), 957--959, Russian.

\bibitem{pontryagin59:_smoot}
\bysame, \emph{Smooth manifolds and their applications in homotopy theory},
  American Mathematical Society Translations, Ser. 2, Vol. 11, American
  Mathematical Society, Providence, R.I., 1959, pp.~1--114.

\bibitem{quillen85:_super_chern}
Daniel Quillen, \emph{Superconnections and the {C}hern character}, Topology
  \textbf{24} (1985), no.~1, 89--95. \MR{86m:58010}

\bibitem{reidemeister72:_einfueh_topol}
Kurt Reidemeister, \emph{Einf\"uhrung in die kombinatorische {T}opologie},
  Wissenschaftliche Buchgesellschaft, Darmstadt, 1972, Unver\"anderter
  reprografischer Nachdruck der Ausgabe Braunschweig 1951. \MR{49 \#9827}

\bibitem{Spanier}
E.~H. Spanier, \emph{Algebraic topology}, McGraw--Hill, New York, 1966.

\bibitem{spanier19:_algeb}
Edwin~H. Spanier, \emph{Algebraic topology}, Springer-Verlag, New York, 19??,
  Corrected reprint of the 1966 original. \MR{96a:55001}

\bibitem{Thom}
R.~Thom, \emph{Quelques propri\'et\'es globales des vari\'et\'es
  differentiables}, Commentarii Mathematici Helvetici \textbf{28} (1954),
  17--86.

\bibitem{whitehead50}
George~W. Whitehead, \emph{The $(n+2)\sp {\rm nd}$ homotopy group of the
  $n$-sphere}, Ann. of Math. (2) \textbf{52} (1950), 245--247.

\bibitem{witten97:_five_m}
Edward Witten, \emph{Five-brane effective action in ${M}$-theory}, J. Geom.
  Phys. \textbf{22} (1997), no.~2, 103--133, hep-th/9610234.

\bibitem{witten97:_m}
\bysame, \emph{On flux quantization in ${M}$-theory and the effective action},
  J. Geom. Phys. \textbf{22} (1997), no.~1, 1--13, hep-th/9609122.

\bibitem{witten01:_overv_k}
\bysame, \emph{Overview of ${K}$-theory applied to strings}, Internat. J.
  Modern Phys. A \textbf{16} (2001), no.~5, 693--706. \MR{1 827 946}

\bibitem{yosimura75:_univer_cw}
Zen-ichi Yosimura, \emph{Universal coefficient sequences for cohomology
  theories of ${\rm {c}{w}}$-spectra}, Osaka J. Math. \textbf{12} (1975),
  no.~2, 305--323. \MR{52 \#9212}

\end{thebibliography}

\providecommand{\bysame}{\leavevmode\hbox to3em{\hrulefill}\thinspace}
\providecommand{\MR}{\relax\ifhmode\unskip\space\fi MR }
\providecommand{\MRhref}[2]{%
  \href{http://www.ams.org/mathscinet-getitem?mr=#1}{#2}
}
\providecommand{\href}[2]{#2}

\end{document}